\documentclass[10pt]{article}

\usepackage{xr}
\usepackage[doi=false,isbn=false,url=false,
	hyperref=auto,
	sorting=nyt,
	maxnames=10,
	maxcitenames=3,
	backend=biber,
	texencoding=auto,
	giveninits=true,
	block=space,
]{biblatex}

\DefineBibliographyStrings{english}{%
	andothers = {et al}
} 
 
\DeclareFieldFormat{eprint:arxiv}{%
	\href{https://arxiv.org/abs/#1}{\emph{arXiv:\allowbreak #1}}}

\DeclareFieldFormat{eprint:mrnumber}{%
	\href{http://www.ams.org/mathscinet-getitem?mr=MR#1}{MR\allowbreak #1}}

\DeclareFieldFormat{eprint:jstor}{%
	\href{http://www.jstor.org/stable/#1}{JSTOR\allowbreak #1}}

\DeclareFieldFormat{eprint:onlineshown}{%
	\em Available at \allowbreak \upshape \ttfamily \href{https://#1}{#1}}

\DeclareFieldFormat{eprint:onlinehidden}{%
	\em Available \href{https://#1}{online}}

\DeclareFieldFormat{eprint:inprep}{%
	\emph{In preparation}}

\DeclareFieldFormat{eprint:manual}{%
	\emph{#1}}

\DeclareFieldFormat{eprint:onarxiv}{%
	\emph{#1}}

\DeclareFieldFormat{eprint:toappear}{%
	To appear in \emph{#1}}

\DeclareFieldFormat{eprint:accepted}{%
	\emph{Accepted at {#1}; to appear}}

\DeclareSourcemap{
	\maps[datatype=bibtex]{
		\map{
			\step[fieldsource=mrnumber,		fieldtarget=eprint, final]
			\step[fieldset=eprinttype,		fieldvalue=mrnumber]
		}
		\map{
			\step[fieldsource=arxiv,		fieldtarget=eprint, final]
			\step[fieldset=eprinttype,		fieldvalue=arxiv]
		}
		\map{
			\step[fieldsource=jstor,		fieldtarget=eprint, final]
			\step[fieldset=eprinttype,		fieldvalue=jstor]
		}
		\map{
			\step[fieldsource=onlineshown,	fieldtarget=eprint, final]
			\step[fieldset=eprinttype,		fieldvalue=onlineshown]
		}
		\map{
			\step[fieldsource=onlinehidden,		fieldtarget=eprint, final]
			\step[fieldset=eprinttype,		fieldvalue=onlinehidden]
		}
		\map{
			\step[fieldsource=inprep,		fieldtarget=eprint, final]
			\step[fieldset=eprinttype,		fieldvalue=inprep]
		}
		\map{
			\step[fieldsource=manual,		fieldtarget=eprint, final]
			\step[fieldset=eprinttype,		fieldvalue=manual]
		}
		\map{
			\step[fieldsource=onarxiv,		fieldtarget=eprint, final]
			\step[fieldset=eprinttype,		fieldvalue=onarxiv]
		}
		\map{
			\step[fieldsource=toappear,		fieldtarget=eprint, final]
			\step[fieldset=eprinttype,		fieldvalue=toappear]
		}
		\map{
			\step[fieldsource=accepted,		fieldtarget=eprint, final]
			\step[fieldset=eprinttype,		fieldvalue=accepted]
		}
	}
}

\DeclareFieldFormat{doi}{%
	\href{https://doi.org/#1}{DOI}%
}

\DeclareFieldFormat{url}{%
	\href{https://#1}{\texttt{#1}}%
}

\DeclareFieldFormat{comments}{%
	\emph{#1}%
} 
 
\DeclareFieldFormat
	[article,inbook,incollection,inproceedings,patent,thesis,unpublished]
	{title}{{#1\isdot}}
\DeclareFieldFormat
	[article,book,inbook,incollection,inproceedings,patent,thesis,unpublished, online]
	{date}{\mkbibparens{#1}}
\DeclareFieldFormat
	[article,book,inbook,incollection,inproceedings,patent,thesis,unpublished]
	{volume}{\mkbibbold{#1}}
\DeclareFieldFormat
	{pages}{\mkbibparens{#1}}

\DeclareFieldFormat{edition}{%
	\ifinteger{#1}
	{\mkbibordedition{#1}~\bibstring{edition}\addcomma}
	{#1~\bibstring{edition}\addcomma}}

\renewbibmacro*{volume+number+eid}{%
	\printfield{volume}%
\setunit*{\adddot}%
	\printfield{number}%
\setunit{\space}%
	\printfield{eid}%
}

\DeclareBibliographyDriver{article}{%
	\usebibmacro{bibindex}%
	\usebibmacro{begentry}%
	\usebibmacro{author}%
\setunit{\addspace}%
	\usebibmacro{date}%
\newunit\newblock
	\usebibmacro{title}%
\newunit\newblock
	\printfield{journaltitle}\addperiod%
\setunit*{\addspace}%
	\usebibmacro{volume+number+eid}%
\setunit*{\addspace}\newblock
	\printfield{pages}
\setunit*{\addspace}%\newblock
	\usebibmacro{eprint}%
\setunit*{\addspace}%
	\printfield{doi}%
	\usebibmacro{finentry}%
}

\DeclareBibliographyDriver{online}{%
	\usebibmacro{bibindex}%
	\usebibmacro{begentry}%
	\usebibmacro{author}%
\setunit{\addspace}%
	\usebibmacro{date}%
\newunit\newblock
	\usebibmacro{title}%
\newunit\newblock
	\usebibmacro{eprint}%
	\usebibmacro{finentry}%
}

\DeclareBibliographyDriver{book}{%
	\usebibmacro{bibindex}%
	\usebibmacro{begentry}%
	\usebibmacro{author}%
\setunit{\addspace}\newblock
	\usebibmacro{date}%
\newunit\newblock
	\usebibmacro{title}%
\setunit*{\addspace}\newblock
	\printfield{edition}%
\newunit\newblock
	\printlist{publisher}%\addcomma
\setunit*{\addspace}%
	\printfield{volume}%
\setunit*{\addspace}\newblock%
	\usebibmacro{eprint}%
\setunit*{\addspace}
	\printfield{doi}
	\usebibmacro{finentry}%
}

\DeclareBibliographyDriver{thesis}{%
	\usebibmacro{bibindex}%
	\usebibmacro{begentry}%
	\usebibmacro{author}%
\setunit{\addspace}\newblock
	\usebibmacro{date}%
\newunit\newblock
	\usebibmacro{title}.%
\setunit*{\addspace}\newblock
	\space Thesis, \printlist{institution}
	\usebibmacro{eprint}%
\setunit*{\addspace}
	\printfield{doi}%
	\usebibmacro{finentry}%
}

\DeclareBibliographyDriver{inproceedings}{%
	\usebibmacro{bibindex}%
	\usebibmacro{begentry}%
	\usebibmacro{author}%
\setunit{\addspace}%
	\usebibmacro{date}%
\newunit\newblock
	\usebibmacro{title}%
\newunit\newblock
	\printfield{booktitle}\addcomma%
\setunit*{\addspace}%
	\printfield{series}\addcomma%
\setunit*{\addspace}\newblock
	\usebibmacro{editor}\addcomma
\setunit*{\addspace}\newblock
	\printlist{publisher}
\setunit*{\addspace}%
	\printfield{volume}%
\setunit*{\addspace}%
	\printfield{pages}
\setunit*{\addspace}%\newblock
	\usebibmacro{eprint}%
\setunit*{\addspace}
	\printfield{doi}
	\usebibmacro{finentry}%
}

\DeclareBibliographyDriver{inbook}{%
	\usebibmacro{bibindex}%
	\usebibmacro{begentry}%
	\usebibmacro{author}%
\setunit{\addspace}%
	\usebibmacro{date}%
\newunit\newblock
	\usebibmacro{title}%
\newunit\newblock
	\printfield{booktitle}\addcomma%
\setunit*{\addspace}%
	\printfield{series}\addcomma%
\setunit*{\addspace}\newblock
	\usebibmacro{editor}\addcomma
\setunit*{\addspace}\newblock
	\printlist{publisher}%\addcomma%
\setunit*{\addspace}%
	\printfield{volume}%
\setunit*{\addspace}%
	\printfield{pages}
\setunit*{\addspace}%\newblock
	\usebibmacro{eprint}%
\setunit*{\addspace}
	\printfield{doi}
	\usebibmacro{finentry}%
}

\DeclareBibliographyDriver{incollection}{%
	\usebibmacro{bibindex}%
	\usebibmacro{begentry}%
	\usebibmacro{author}%
\setunit{\addspace}%
	\usebibmacro{date}%
\newunit\newblock
	\usebibmacro{title}%
\newunit\newblock
	\printfield{booktitle}\addcomma%
\setunit*{\addspace}%
	\printfield{series}\addcomma%
\setunit*{\addspace}\newblock
	\printlist{publisher}%\addcomma%
\setunit*{\addspace}%
	\printfield{volume}%
\setunit*{\addspace}%
	\printfield{pages}
\setunit*{\addspace}%\newblock
	\usebibmacro{eprint}%
\setunit*{\addspace}
	\printfield{doi}
	\usebibmacro{finentry}%
}

\AtEveryBibitem{%
	\clearfield{day}%
	\clearfield{month}%
} 
 
\DeclareNameAlias{sortname}{given-family}
\DeclareNameAlias{default}{given-family}

\newcommand{\acknofootnote}{%
	\blfootnote{%
		\hspace*{1em}%
		Jonathan Hermon%
	\hfill%
		Sam Olesker-Taylor%
		\hspace*{1em}%
	\\%
		\href{mailto:jhermon@math.ubc.ca}{jhermon@math.ubc.ca},
		\href{http://www.math.ubc.ca/~jhermon/}{math.ubc.ca/$\sim$jhermon/}%
	\hfill%
		\href{mailto:oleskertaylor.sam@gmail.com}{oleskertaylor.sam@gmail.com},
		\href{https://sites.google.com/view/sam-ot/}{sites.google.com/view/sam-ot/}%
	\\%
		University of British Columbia, Vancouver, Canada%
	\hfill%
		Department of Statistics, University of Warwick, UK
	\\%
		Supported by EPSRC EP/L018896/1 and an NSERC Grant%
	\hfill%
		Supported by EPSRC Grants 1885554 and EP/N004566/1%
	\par\smallskip\par
	\centering%
		The vast majority of this work was undertaken whilst both authors were at the University of Cambridge%
	}
}

\newcommand*{\mm}{\ensuremath{L}}

\newcommand*{\MM}{\ensuremath{\Gamma}}

\newcommand*{\tent}{\tau}
\newcommand*{\sent}{\sigma}
\newcommand*{\tpro}{\mathfrak t}
\newcommand*{\spro}{\mathfrak s}

\newcommand{\printtoc}[1]{%
	\ifthenelse%
		{\equal{#1}{1}}%
		{\sffamily\boldmath\tableofcontents\unboldmath\normalfont}%
		{\newpage\small\sffamily\boldmath\tableofcontents\unboldmath\normalfont\normalsize}%
	}

\newcommand{\nextresult}{%
	\setcounter{introthm}{\value{introthm}}
	\setcounter{introcor}{\value{introthm}}
	\setcounter{introconj}{\value{introthm}}
	\setcounter{introdefn}{\value{introthm}}
	\setcounter{intrormkT}{\value{introthm}}
}

\usepackage[utf8]{inputenc}

\usepackage{empheq}

\usepackage{manyfoot}
\SetFootnoteHook{\hspace*{-1.8em}}
\DeclareNewFootnote{bl}[gobble]
\newcommand{\blfootnote}[1]{\footnotebl{\sffamily#1}}

\usepackage{graphicx}
\graphicspath{{../Figures/}}
\usepackage[labelsep=period]{caption}

\usepackage{chngcntr}
\counterwithin{figure}{section}

\usepackage[UKenglish]{babel}
\usepackage{amsmath, amsthm, amssymb, mathrsfs, bm, mathtools}
\usepackage{wasysym}
\usepackage{mathtools}

\allowdisplaybreaks

\usepackage{cases}

\usepackage{titlesec}
\usepackage{titletoc}

\titleformat{\title}
{\sffamily \Huge \bfseries \scshape \boldmath}{}{}{}
\titleformat{\section}
{\sffamily \Large \bfseries \scshape \boldmath}{\thesection}{1em}{}
\titleformat{\subsection}
{\sffamily \large \bfseries \scshape \boldmath}{\thesubsection}{1em}{}
\titleformat{\subsubsection}
{\sffamily \normalsize \bfseries \scshape \boldmath}{\thesubsubsection}{1em}{}
\titleformat{\paragraph}
{\sffamily \normalsize \bfseries \scshape \boldmath}{\theparagraph}{1em}{}
\titleformat{\subparagraph}[runin]
{\sffamily \normalsize \bfseries \scshape \boldmath}{\thesubparagraph}{1em}{}

\def\IfAmpersandUseAlign#1#2&#3\EndIfAmpersandUseAlign
{%
	\if\relax\detokenize{#3}\relax
	\begin{equation*}%
	#1%
	\end{equation*}%
	\else
	\begin{align*}%
	#1%
	\end{align*}%
	\fi
}
\def\[#1\]%
{%
	\IfAmpersandUseAlign{#1}#1&\EndIfAmpersandUseAlign
}

\usepackage{eqparbox}
\newcommand{\eqmathsbox}[3][\mathrel]{%
	#1{\eqmakebox[#2]{$\displaystyle#3$}}%
}
\usepackage{scrextend}

\usepackage{enumitem}

\newcommand{\numberingroman}{%
	\renewcommand{\labelenumi}{(\roman{enumi})}%
	\renewcommand{\theenumi}{(\roman{enumi})}%
}

\setlist[description]{%
	topsep		= 0pt,		% space before start / after end of list
	noitemsep,				% space between items
	font		= {\mdseries\itshape},	% set the label font
}

\newcommand{\ggr}{G}
\newcommand{\gab}{G^\ab}
\newcommand{\gcom}{G^\com}

\newcommand{\UU}{U}
\newcommand{\HH}{H}

\newcommand{\hab}{\HH_{p,d}^\ab}
\newcommand{\hcom}{\HH_{p,d}^\com}

\usepackage[dvipsnames,hyperref]{xcolor}

\newcommand{\bcdot}{\ensuremath{\bm{\cdot}}}
\let\mod\relax
\DeclareMathOperator{\mod}{\, mod}

\DeclareMathOperator{\lcm}{lcm}

\let\div\relax
\DeclareMathOperator{\div}{div}

\newcommand{\Quad}[1]{
	\mathchoice
	{\quad\text{#1}\quad}%	\displaystyle
	{\text{ #1 }}%			\textsyle
	{\text{ #1 }}%			\scriptstyle
	{\text{ #1 }}%			\scriptscriptstyle
}

\newcommand{\whp}{\text{whp}\xspace}

\newcommand{\forallZ}{\text{for all $Z$}\xspace}

\newcommand{\Qforall}{\Quad{for all}}
\newcommand{\Qfor}{\Quad{for}}
\newcommand{\Qand}{\Quad{and}}
\newcommand{\Qwhere}{\Quad{where}}
\newcommand{\Qwhen}{\Quad{when}}

\newcommand{\id}{\mathsf{id}}

\newcommand{\loc}{\mathrm{loc}}
\newcommand{\glo}{\mathrm{glo}}
\newcommand{\typ}{\mathsf{typ}}

\newcommand{\cq}{\coloneqq}
\renewcommand{\epsilon}{\varepsilon}

\newcommand{\eps}{\epsilon}

\usepackage{xifthen}

\newcommand{\binomt}[2]{ \textstyle \binom{#1}{#2} \displaystyle }

\newcommand{\maxt}[1]{ \textstyle \max_{#1} \displaystyle }

\newcommand{\mint}[1]{ \textstyle \min_{#1} \displaystyle }

\usepackage{calc}
\newlength{\halfplusheight}
\newcommand{\MAX}[1]{\mathop{\raisebox{\halfplusheight}{\(\displaystyle\max_{#1}\)}}}
\newcommand{\MIN}[1]{\mathop{\raisebox{\halfplusheight}{\(\displaystyle\min_{#1}\)}}}
\newcommand{\LIM}[1]{\mathop{\raisebox{\halfplusheight}{\(\displaystyle\lim_{#1}\)}}}
\newcommand{\LIMSUP}[1]{\mathop{\raisebox{\halfplusheight}{\(\displaystyle\limsup_{#1}\)}}}
\newcommand{\LIMINF}[1]{\mathop{\raisebox{\halfplusheight}{\(\displaystyle\liminf_{#1}\)}}}

\DeclareMathOperator*{\sumTT}{\textstyle\sum}
\newcommand{\sumT}[2][]{
	\ifthenelse{\isempty{#1}}
	{\sumTT_{#2}}
	{\sumTT_{#2}^{#1}}
}
\newcommand{\sumt}[2][]{
	\ifthenelse{\isempty{#1}}
	{\textstyle \sum_{#2}      \displaystyle}
	{\textstyle \sum_{#2}^{#1} \displaystyle}
}
\newcommand{\sumd}[2][]{
	\ifthenelse{\isempty{#1}}
	{\displaystyle \sum_{#2}}
	{\displaystyle \sum_{#2}^{#1}}
}

\newcommand{\intt}[2][]{
	\ifthenelse{\isempty{#1}}
	{\textstyle \int_{#2}      \displaystyle}
	{\textstyle \int_{#2}^{#1} \displaystyle}
}

\newcommand{\prodt}[2][]{
	\ifthenelse{\isempty{#1}}
	{\textstyle \prod_{#2}      \displaystyle}
	{\textstyle \prod_{#2}^{#1} \displaystyle}
}
\newcommand{\prodd}[2][]{
	\ifthenelse{\isempty{#1}}
	{\prod_{#2}}
	{\prod_{#2}^{#1}}
}

\usepackage{xspace}

\let\originalexp\exp
\let\exp\relax

\DeclareRobustCommand{\exp} [1]{\originalexp(#1)}
\newcommand{\expb} [1]{\originalexp\bigl( #1 \bigr)}

\newcommand{\abs}  [1]{| #1 |}
\newcommand{\absb} [1]{\bigl| #1 \bigr|}

\newcommand{\norm}  [1]{\lVert #1 \rVert}
\newcommand{\normb} [1]{\big\lVert #1 \bigr\rVert}

\newcommand{\rbr} [1]{ ( #1 ) }
\newcommand{\rbb} [1]{\bigl( #1 \bigr)}

\newcommand{\sbr} [1]{ [ #1 ] }
\newcommand{\sbb} [1]{\bigl[ #1 \bigr]}

\newcommand{\bra} [1]{ \{ #1 \} }
\newcommand{\brb} [1]{\bigl\{ #1 \bigr\}}

\newcommand{\brbb}[1]{\biggl\{ #1 \biggr\}}

\DeclareMathOperator{\diam}{diam}

\newcommand{\mix}{\mathrm{mix}}

\newcommand{\rel}{\mathrm{rel}}

\newcommand{\st}{{ \ \mathrm{st} \ }}

\newcommand{\ab} {\mathrm{ab}}
\newcommand{\com}{\mathrm{com}}

\DeclareMathOperator{\step}{step}

\newcommand{\Unif}{\mathrm{Unif}}
\newcommand{\iid}{\mathrm{iid}}

\newcommand{\Bin}{\mathrm{Bin}}

\newcommand{\tmix}{t_\mix}

\newcommand{\trel}{t_\rel}

\newcommand{\Ninn}{{N\in\mathbb{N}}}

\newcommand{\floor}[1]{\lfloor #1 \rfloor}
\newcommand{\floorb}[1]{\bigl \lfloor #1 \bigr \rfloor}

\newcommand{\ceil}[1]{\lceil #1 \rceil}

\newcommand{\midb}{\bigm\vert}

\newcommand{\one}  [1]{\bm1( #1 )}

\newcommand{\tv} [1]{\lVert #1 \rVert_{\mathrm{TV}}}
\newcommand{\tvb}[1]{\bigl\lVert #1 \bigr\rVert_{\mathrm{TV}}}
\newcommand{\logk}[1][]{
	\ifthenelse{\equal{}{#1}}
	{\log k}
	{(\log k)^{#1}}
}

\newcommand{\logn}[1][]{
	\ifthenelse{\equal{}{#1}}
	{\log n}
	{(\log n)^{#1}}
}

\newcommand{\logm}[1][]{
	\ifthenelse{\equal{}{#1}}
	{\log m}
	{(\log m)^{#1}}
}

\newcommand{\loglogn}[1][]{
	\ifthenelse{\equal{}{#1}}
	{\log\log n}
	{(\log\log n)^{#1}}
}

\newcommand{\prt}[2][]{
	\ifthenelse{\equal{}{#1}}
	{\mathbb{P}(#2)}
	{\mathbb{P}_{#1}(#2)}
}
\newcommand{\pr}[2][]{
	\mathchoice
	{\ifthenelse{\isempty{#1}}
		{\mathbb{P}\bigl(#2\bigr)}
		{\mathbb{P}_{#1}\bigl(#2\bigr)}}%	\displaystyle
	{\ifthenelse{\isempty{#1}}
		{\mathbb{P}(#2)}
		{\mathbb{P}_{#1}(#2)}}%	\textsyle
	{\ifthenelse{\isempty{#1}}
		{\mathbb{P}(#2)}
		{\mathbb{P}_{#1}(#2)}}%	\scriptstyle
	{\ifthenelse{\isempty{#1}}
		{\mathbb{P}(#2)}
		{\mathbb{P}_{#1}(#2)}}%	\scriptscriptstyle
}
\newcommand{\prb}[2][]{
	\ifthenelse{\equal{}{#1}}
	{\mathbb{P}\bigl( #2 \bigr)}
	{\mathbb{P}_{#1}\bigl( #2 \bigr)}
}
\newcommand{\prB}[2][]{
	\ifthenelse{\equal{}{#1}}
	{\mathbb{P}\Bigl( #2 \Bigr)}
	{\mathbb{P}_{#1}\Bigl( #2 \Bigr)}
}
\newcommand{\prbb}[2][]{
	\ifthenelse{\equal{}{#1}}
	{\mathbb{P}\biggl( #2 \biggr)}
	{\mathbb{P}_{#1}\biggl( #2 \biggr)}
}
\newcommand{\prBB}[2][]{
	\ifthenelse{\equal{}{#1}}
	{\mathbb{P}\Biggl( #2 \Biggr)}
	{\mathbb{P}_{#1}\Biggl( #2 \Biggr)}
}
\newcommand{\prs}[2][]{
	\ifthenelse{\equal{}{#1}}
	{\mathbb{P}\left( #2 \right)}
	{\mathbb{P}_{#1}\left( #2 \right)}
}

\newcommand{\qr}[2][]{
	\mathchoice
	{\ifthenelse{\isempty{#1}}
		{\mathbb{Q}\bigl(#2\bigr)}
		{\mathbb{Q}_{#1}\bigl(#2\bigr)}}%	\displaystyle
	{\ifthenelse{\isempty{#1}}
		{\mathbb{Q}(#2)}
		{\mathbb{Q}_{#1}(#2)}}%	\textsyle
	{\ifthenelse{\isempty{#1}}
		{\mathbb{Q}(#2)}
		{\mathbb{Q}_{#1}(#2)}}%	\scriptstyle
	{\ifthenelse{\isempty{#1}}
		{\mathbb{Q}(#2)}
		{\mathbb{Q}_{#1}(#2)}}%	\scriptscriptstyle
}
\newcommand{\qrb}[2][]{
	\ifthenelse{\equal{}{#1}}
	{\mathbb{Q}\bigl( #2 \bigr)}
	{\mathbb{Q}_{#1}\bigl( #2 \bigr)}
}
\newcommand{\qrB}[2][]{
	\ifthenelse{\equal{}{#1}}
	{\mathbb{Q}\Bigl( #2 \Bigr)}
	{\mathbb{Q}_{#1}\Bigl( #2 \Bigr)}
}
\newcommand{\qrbb}[2][]{
	\ifthenelse{\equal{}{#1}}
	{\mathbb{Q}\biggl( #2 \biggr)}
	{\mathbb{Q}_{#1}\biggl( #2 \biggr)}
}
\newcommand{\qrBB}[2][]{
	\ifthenelse{\equal{}{#1}}
	{\mathbb{Q}\Biggl( #2 \Biggr)}
	{\mathbb{Q}_{#1}\Biggl( #2 \Biggr)}
}
\newcommand{\qrs}[2][]{
	\ifthenelse{\equal{}{#1}}
	{\mathbb{Q}\left( #2 \right)}
	{\mathbb{Q}_{#1}\left( #2 \right)}
}

\newcommand{\ext}[2][]{
\ifthenelse{\equal{}{#1}}
{\mathbb{E}(#2)}
{\mathbb{E}_{#1}(#2)}
}
\newcommand{\ex}[2][]{
	\mathchoice
	{\ifthenelse{\isempty{#1}}
		{\mathbb{E}\bigl(#2\bigr)}
		{\mathbb{E}_{#1}\bigl(#2\bigr)}}%	\displaystyle
	{\ifthenelse{\isempty{#1}}
		{\mathbb{E}(#2)}
		{\mathbb{E}_{#1}(#2)}}%	\textsyle
	{\ifthenelse{\isempty{#1}}
		{\mathbb{E}(#2)}
		{\mathbb{E}_{#1}(#2)}}%	\scriptstyle
	{\ifthenelse{\isempty{#1}}
		{\mathbb{E}(#2)}
		{\mathbb{E}_{#1}(#2)}}%	\scriptscriptstyle
}
\newcommand{\exb}[2][]{
	\ifthenelse{\equal{}{#1}}
	{\mathbb{E}\bigl( #2 \bigr)}
	{\mathbb{E}_{#1}\bigr( #2 \bigr)}
}
\newcommand{\exB}[2][]{
	\ifthenelse{\equal{}{#1}}
	{\mathbb{E}\Bigl( #2 \Bigr)}
	{\mathbb{E}_{#1}\Bigl( #2 \Bigr)}
}
\newcommand{\exbb}[2][]{
	\ifthenelse{\equal{}{#1}}
	{\mathbb{E}\biggl( #2 \biggr)}
	{\mathbb{E}_{#1}\biggl( #2 \biggr)}
}
\newcommand{\exBB}[2][]{
	\ifthenelse{\equal{}{#1}}
	{\mathbb{E}\Biggl( #2 \Biggr)}
	{\mathbb{E}_{#1}\Biggl( #2 \Biggr)}
}

\newcommand{\fx}[2][]{
	\ifthenelse{\equal{}{#1}}
	{\mathbb{F}(#2)}
	{\mathbb{F}_{#1}(#2)}
}
\newcommand{\fxb}[2][]{
	\ifthenelse{\equal{}{#1}}
	{\mathbb{F}\bigl( #2 \bigr)}
	{\mathbb{F}_{#1}\bigr( #2 \bigr)}
}
\newcommand{\fxB}[2][]{
	\ifthenelse{\equal{}{#1}}
	{\mathbb{F}\Bigl( #2 \Bigr)}
	{\mathbb{F}_{#1}\Bigl( #2 \Bigr)}
}
\newcommand{\fxbb}[2][]{
	\ifthenelse{\equal{}{#1}}
	{\mathbb{F}\biggl( #2 \biggr)}
	{\mathbb{F}_{#1}\biggl( #2 \biggr)}
}
\newcommand{\fxBB}[2][]{
	\ifthenelse{\equal{}{#1}}
	{\mathbb{F}\Biggl( #2 \Biggr)}
	{\mathbb{F}_{#1}\Biggl( #2 \Biggr)}
}

\newcommand{\Var}[1]{\mathbb{V}\mathrm{ar}(#1)}
\newcommand{\Varb}[2][]{
	\ifthenelse{\equal{}{#1}}
	{\mathbb{V}\mathrm{ar} \bigl(#2\bigr)}
	{\mathbb{V}\mathrm{ar}_{#1} \bigl(#2\bigr)}
}
\newcommand{\VAR}[2][]{
	\ifthenelse{\equal{}{#1}}
	{\mathrm{Var}(#2)}
	{\mathrm{Var}_{#1}(#2)}
}

\newcommand{\Oh}  [1]{\mathcal{O}( #1 )}

\newcommand{\oh}  [1]{o( #1 )}

\newcommand{\mbe}{\mathbb{E}}

\newcommand{\mbn}{\mathbb{N}}

\newcommand{\mbp}{\mathbb{P}}

\newcommand{\mbr}{\mathbb{R}}

\newcommand{\mbz}{\mathbb{Z}}

\newcommand{\mca}{\mathcal{A}}

\newcommand{\mce}{\mathcal{E}}

\newcommand{\mch}{\mathcal{H}}

\newcommand{\mcw}{\mathcal{W}}

\newcommand{\mfgcd}{\mathfrak{g}}

\newcommand{\toinf}[1]{\ensuremath{#1\to\infty}}
\newcommand{\asinf}[1]{\text{as \ensuremath{#1\to\infty}}}

\newcommand{\Ninf}{{N\to\infty}}

\makeatletter
\newenvironment{subtheorem}[1]{%
	\def\subtheoremcounter{#1}%
	\refstepcounter{#1}%
	\protected@edef\theparentnumber{\csname the#1\endcsname}%
	\setcounter{parentnumber}{\value{#1}}%
	\setcounter{#1}{0}%
	\expandafter\def\csname the#1\endcsname{\theparentnumber\alph{#1}}%
	\expandafter\def\csname theH#1\endcsname{thm.\theparentnumber\alph{#1}}%
	\unskip\ignorespaces
}{%
	\setcounter{\subtheoremcounter}{\value{parentnumber}}%
	\ignorespacesafterend
}
\makeatother
\newcounter{parentnumber}

\makeatletter
\newenvironment{subtheorem-num}[1]{%
	\def\subtheoremcounter{#1}%
	\refstepcounter{#1}%
	\protected@edef\theparentnumber{\csname the#1\endcsname}%
	\setcounter{parentnumber}{\value{#1}}%
	\setcounter{#1}{0}%
	\expandafter\def\csname the#1\endcsname{\theparentnumber.\arabic{#1}}%
	\expandafter\def\csname theH#1\endcsname{thm.\theparentnumber.\arabic{#1}}%
	\unskip\ignorespaces
}{%
	\setcounter{\subtheoremcounter}{\value{parentnumber}}%
	\ignorespacesafterend
}
\makeatother
\usepackage[colorlinks, pdfauthor = {Jonathan\ Hermon\ and\ Sam\ Olesker-Taylor}, pdfusetitle]{hyperref}
\hypersetup{
	colorlinks,
	linkcolor	= {blue},
	filecolor	= {red},
	urlcolor	= {blue},
	citecolor	= {ForestGreen}
}

\newcommand{\qedtriangle}{\renewcommand{\qedsymbol}{\ensuremath{\triangle}}}

\usepackage[capitalise]{cleveref}

\crefformat{equation}{\textup{(#2#1#3)}}
\crefmultiformat{equation}{\textup{(#2#1#3}}{\textup{, #2#1#3)}}{\textup{, #2#1#3}}{\textup{, #2#1#3)}}
\crefrangeformat{equation}{\textup{(#3#1#4--#5#2#6)}}

\newenvironment{Proof}[1][\proofname]{%
	\proof[\upshape\bfseries\sffamily\boldmath#1]
}{\endproof}

\usepackage{etoolbox}
\usepackage{needspace}

\newtheoremstyle{sfsl}
{1\baselineskip}		% Space above
{1\baselineskip}		% Space below
{\slshape}				% Theorem body font
{}						% Indent amount
{\bfseries\sffamily}	% Theorem head font
{.}						% Punctuation after theorem head
{0.5em}					% Space after theorem head
{\thmname{#1}\thmnumber{ #2}\thmnote{ {\mdseries(#3)}}}
\newtheoremstyle{sfup}
{1\baselineskip}		% Space above
{1\baselineskip}		% Space below
{\upshape}				% Theorem body font
{}						% Indent amount
{\bfseries\sffamily}	% Theorem head font
{.}						% Punctuation after theorem head
{0.5em}					% Space after theorem head
{\thmname{#1}\thmnumber{ #2}\thmnote{ {\mdseries(#3)}}}
\theoremstyle{sfsl}

\newtheorem*{thm*}{Theorem}
\newtheorem{thm} {Theorem}[section]
\crefname{thm}{Theorem}{Theorems}

\newtheorem*{introthm*}{Theorem}
\newtheorem{introthm}{Theorem}

\crefname{introthm}{Theorem}{Theorems}

\newtheorem*{cor*}{Corollary}
\newtheorem{cor} [thm]{Corollary}
\crefname{cor}{Corollary}{Corollaries}

\newtheorem*{introcor*}{Corollary}
\newtheorem{introcor}{Corollary}

\crefname{introcor}{Corollary}{Corollaries}

\newtheorem*{introconj*}{Conjecture}
\newtheorem{introconj}{Conjecture}

\crefname{introconj}{Conjecture}{Conjectures}

\newtheorem*{introques*}{Question}

\crefname{introques}{Question}{Questions}

\newtheorem*{lem*}    {Lemma}
\newtheorem{lem} [thm]{Lemma}
\crefname{lem}{Lemma}{Lemmas}

\newtheorem*{introlem*}{Lemma}

\crefname{introlem}{Lemma}{Lemmas}

\newtheorem*{prop*}    {Proposition}
\newtheorem{prop} [thm]{Proposition}
\crefname{prop}{Proposition}{Propositions}

\newtheorem*{clm*}    {Claim}

\crefname{clm}{Claim}{Claims}

\newtheorem*{defn*}    {Definition}
\newtheorem{defn} [thm]{Definition}
\crefname{defn}{Definition}{Definitions}

\newtheorem*{introdefn*}{Definition}
\newtheorem{introdefn}{Definition}

\crefname{introdefn}{Definition}{Definitions}

\newtheorem{innercustomgenericsl}{\customgenericnamesl}
\providecommand{\customgenericnamesl}{}
\newcommand{\newcustomtheoremsl}[2]{%
	\newenvironment{#1}[1]
	{%
		\renewcommand\customgenericnamesl{#2}%
		\renewcommand\theinnercustomgenericsl{##1}%
		\innercustomgenericsl
	}
	{\endinnercustomgenericsl}
}

\newcustomtheoremsl{customthm}{Theorem}
\newcustomtheoremsl{customprop}{Proposition}
\newcustomtheoremsl{customlemma}{Lemma}
\newcustomtheoremsl{customrmk}{Remark}
\newcustomtheoremsl{customconj}{Conjecture}
\newcustomtheoremsl{customdefn}{Definition}
\newcustomtheoremsl{customquestion}{Question}
\newcustomtheoremsl{customopenproblem}{Open Problem}
\newcustomtheoremsl{customhyp}{Hypothesis}

\newtheorem*{conj*}   {Conjecture}

\crefname{conj}{Conjecture}{Conjectures}

\newenvironment{conj-ind*}
	{\begin{quote}\textsf{\textbf{Conjecture.}}\slshape}
	{\end{quote}}
\newenvironment{conj-ind}
	{\begin{quote}\vspace{-\glueexpr\baselineskip+\topsep}\begin{customconj}}
	{\end{customconj}\end{quote}}

\newenvironment{question-ind*}
	{\begin{quote}\textsf{\textbf{Question.}}\slshape}
	{\end{quote}}
\newenvironment{question-ind}
	{\begin{quote}\vspace{-\glueexpr\baselineskip+\topsep}\begin{customquestion}}
	{\end{customquestion}\end{quote}}

\newenvironment{openproblem-ind*}
	{\begin{quote}\textsf{\textbf{Open Problem.}}\slshape}
	{\end{quote}}
\newenvironment{openproblem-ind}
	{\begin{quote}\vspace{-\glueexpr\baselineskip+\topsep}\begin{customopenproblem}}
	{\end{customopenproblem}\end{quote}}

\newtheorem*{hypothesis*}{Hypothesis}

\newtheorem*{hyp*}{Hypothesis}
\newtheorem{hyp}{Hypothesis}

\crefname{hyp}{Hypothesis}{Hypotheses}

\newtheorem*{rmk*}{Remark}

\theoremstyle{sfup}

\newtheorem{innercustomgenericup}{\customgenericnameup}
\providecommand{\customgenericnameup}{}
\newcommand{\newcustomtheoremup}[2]{%
	\newenvironment{#1}[1]
	{%
		\renewcommand\customgenericnameup{#2}%
		\renewcommand\theinnercustomgenericup{##1}%
		\innercustomgenericup
	}
	{\endinnercustomgenericup}
}

\crefname{exm} {Example}{Examples}
\crefname{exmT}{Example}{Examples}

\newtheorem{rmkT}[thm]{Remark}
	\newenvironment{rmkt}
	{\pushQED{\qed}\renewcommand{\qedsymbol}{\ensuremath{\triangle}}\rmkT}
	{\popQED\endrmkT}
\crefname{rmk} {Remark}{Remarks}
\crefname{rmkT}{Remark}{Remarks}

\newtheorem*{rmkTT}{Remark}
\newenvironment{rmkt*}
	{\pushQED{\qed}\renewcommand{\qedsymbol}{\ensuremath{\triangle}}\rmkTT}
	{\popQED\endrmkTT}

\newcustomtheoremup{customrmkT}{Remark}

\crefname{rmks} {Remarks}{Remarks}
\crefname{rmksT}{Remarks}{Remarks}

\newtheorem*{rmks*} {Remarks}
\newtheorem*{rmksTT}{Remarks}
\newenvironment{rmkst*}
	{\pushQED{\qed}\renewcommand{\qedsymbol}{\ensuremath{\triangle}}\rmksTT}
	{\popQED\endrmksTT}

\newtheorem{intrormkT}{Remark}
	\newenvironment{intrormkt}
	{\pushQED{\qed}\renewcommand{\qedsymbol}{\ensuremath{\triangle}}\intrormkT}
	{\popQED\endintrormkT}

\crefname{intrormk} {Remark}{Remarks}
\crefname{intrormkT}{Remark}{Remarks}

\newtheorem*{intrormk*} {Remark}
\newtheorem*{intrormkTT}{Remark}
\newenvironment{intrormkt*}
	{\pushQED{\qed}\renewcommand{\qedsymbol}{\ensuremath{\triangle}}\intrormkTT}
	{\popQED\endintrormkTT}

\newtheorem*{exm*} {Example}
\newtheorem*{exmTT}{Lemma}
	\newenvironment{exmt*}
	{\pushQED{\qed}\renewcommand{\qedsymbol}{\ensuremath{\triangle}}\exmTT}
	{\popQED\endexmTT}

\newtheorem*{note*} {Note}
\newtheorem*{noteTT}{Note}
	\newenvironment{notet*}
	{\pushQED{\qed}\renewcommand{\qedsymbol}{\ensuremath{\triangle}}\noteTT}
	{\popQED\endnoteTT}

\catcode`\^^A=14

\usepackage{xparse}

\newcounter{mixedsubequations}

\ExplSyntaxOn

\int_new:N \g_mixedsubeq_int

\NewDocumentEnvironment{mixedsubequations}{o}
{
	\IfNoValueTF { #1 }
	{
		\addtocounter{equation}{-\g_mixedsubeq_int}
		\stepcounter{mixedsubequations}
	}
	{
		\int_gset:Nn \g_mixedsubeq_int { \clist_count:n { #1 } }
		\clist_map_inline:nn { #1 }
		{
			\refstepcounter{equation}\label{##1}
		}
		\addtocounter{equation}{-\g_mixedsubeq_int}
		\setcounter{mixedsubequations}{1}
	}
	\domixedsubequations
}
{\ignorespacesafterend}

\NewDocumentCommand{\domixedsubequations}{}
{
	\cs_set:Npx \theequation
	{
		\exp_not:o { \theequation }
		\exp_not:n { \alph{mixedsubequations} }
	}
	\ignorespaces
}
\ExplSyntaxOff

\makeatletter
\let\save@mathaccent\mathaccent
\newcommand*\if@single[3]{%
  \setbox0\hbox{${\mathaccent"0362{#1}}^H$}%
  \setbox2\hbox{${\mathaccent"0362{\kern0pt#1}}^H$}%
  \ifdim\ht0=\ht2 #3\else #2\fi
  }
\newcommand*\rel@kern[1]{\kern#1\dimexpr\macc@kerna}
\newcommand*\widebar[1]{\@ifnextchar^{{\wide@bar{#1}{0}}}{\wide@bar{#1}{1}}}
\newcommand*\wide@bar[2]{\if@single{#1}{\wide@bar@{#1}{#2}{1}}{\wide@bar@{#1}{#2}{2}}}
\newcommand*\wide@bar@[3]{%
  \begingroup
  \def\mathaccent##1##2{%
    \let\mathaccent\save@mathaccent
    \if#32 \let\macc@nucleus\first@char \fi
    \setbox\z@\hbox{$\macc@style{\macc@nucleus}_{}$}%
    \setbox\tw@\hbox{$\macc@style{\macc@nucleus}{}_{}$}%
    \dimen@\wd\tw@
    \advance\dimen@-\wd\z@
    \divide\dimen@ 3
    \@tempdima\wd\tw@
    \advance\@tempdima-\scriptspace
    \divide\@tempdima 10
    \advance\dimen@-\@tempdima
    \ifdim\dimen@>\z@ \dimen@0pt\fi
    \rel@kern{0.6}\kern-\dimen@
    \if#31
      \overline{\rel@kern{-0.6}\kern\dimen@\macc@nucleus\rel@kern{0.4}\kern\dimen@}%
      \advance\dimen@0.4\dimexpr\macc@kerna
      \let\final@kern#2%
      \ifdim\dimen@<\z@ \let\final@kern1\fi
      \if\final@kern1 \kern-\dimen@\fi
    \else
      \overline{\rel@kern{-0.6}\kern\dimen@#1}%
    \fi
  }%
  \macc@depth\@ne
  \let\math@bgroup\@empty \let\math@egroup\macc@set@skewchar
  \mathsurround\z@ \frozen@everymath{\mathgroup\macc@group\relax}%
  \macc@set@skewchar\relax
  \let\mathaccentV\macc@nested@a
  \if#31
    \macc@nested@a\relax111{#1}%
  \else
    \def\gobble@till@marker##1\endmarker{}%
    \futurelet\first@char\gobble@till@marker#1\endmarker
    \ifcat\noexpand\first@char A\else
      \def\first@char{}%
    \fi
    \macc@nested@a\relax111{\first@char}%
  \fi
  \endgroup
}
\makeatother

\numberwithin{equation}{section}

\title{\sffamily Cutoff for Almost All Random Walks on Abelian Groups}

\author{\sffamily Jonathan Hermon\quad Sam Olesker-Taylor}

\date{}
\usepackage{geometry}
\begin{document}

%	\noindent
%odd side: \the\oddsidemargin\\
%even side: \the\evensidemargin\\
%top: \the\topmargin\\
%%bottom: \the\bottommargin\\
%text height: \the\textheight\\
%paper height: \the\paperheight\\
%paper width: \the\paperwidth\\
%text width: \the\textwidth\\
%line width: \the\linewidth\\
%
%\newlength{\mytextsize}
%\makeatletter
%\show\f@size
%\setlength{\mytextsize}{\f@size pt}
%\makeatother
%\the\mytextsize

\maketitle

\acknofootnote

\vspace{-6ex}

\renewcommand{\abstractname}{\sffamily Abstract}
\begin{abstract}
Consider the random Cayley graph of a finite group $G$ with respect to $k$ generators chosen uniformly at random, with $1 \ll \log k \ll \log |G|$; denote it $G_k$. A conjecture of Aldous and Diaconis~\cite{AD:conjecture} asserts, for $k \gg \log |G|$, that the random walk on this graph exhibits cutoff. Further, the cutoff time should be a function only of $k$ and $|G|$, to sub-leading order.

This was verified for all Abelian groups in the '90s. We extend the conjecture to $1 \ll k \lesssim \log |G|$. We establish cutoff for \emph{all} Abelian groups under the condition $k - d(G) \gg 1$, where $d(G)$ is the minimal size of a generating subset of $G$, which is almost optimal. The cutoff time is described (abstractly) in terms of the entropy of random walk on $\mathbb Z^k$. This abstract definition allows us to deduce that the cutoff time can be written as a function only of $k$ and $|G|$ when $d(G) \ll \log |G|$ and $k - d(G) \asymp k \gg 1$; this is not the case when $d(G) \asymp \log |G| \asymp k$. For certain regimes of $k$, we find the limit profile of the convergence to equilibrium.

Wilson~\cite{W:rws-hypercube} conjectured that $\mathbb Z_2^d$ gives rise to the slowest mixing time for $G_k$ amongst all groups of size at most $2^d$. We give a partial answer, verifying the conjecture for nilpotent groups. This is obtained via a comparison result of independent interest between the mixing times of nilpotent $G$ and a corresponding Abelian group $\widebar G$, namely the direct sum of the Abelian quotients in the lower central series of $G$. We use this to refine a celebrated result of Alon and Roichman~\cite{AR:cayley-expanders}: we show for nilpotent $G$ that $G_k$ is an expander provided $k - d(\widebar G) \gtrsim \log \abs G$. As another consequence, we establish cutoff for nilpotent groups with relatively small commutator subgroup, including high-dimensional special groups, such as Heisenberg groups.

The aforementioned results all hold with high probability over the random Cayley graph~$G_k$.
\end{abstract}

\small
\begin{quote}
\begin{description}
	\item [Keywords:]
	cutoff, mixing times, random walk, random Cayley graphs, entropy
	
	\item [MSC 2020 subject classifications:]
	05C48, 05C80, 05C81; 20D15; 60B15, 60J27, 60K37
\end{description}
\end{quote}
\normalsize

%05C12: Distance in graphs
%05C48: Expander graphs
%05C80:	Random graphs
%05C81:	Random walks on graphs
%20D15: Nilpotent groups, p-groups
%60B15: Probability measures on groups or semigroups, Frouier transforms, factorization
%60C05: Combinatorial probability
%60J10: Markov chains (discrete-time Markov processes on discrete state spaces)
%60J27:	Continuous-time Markov processes on discrete state spaces
%60K37:	Processes in random environments

\vspace{2ex}
\numberingroman
\renewcommand{\mm}{\ensuremath{\gamma}}

\setcounter{tocdepth}{1}
\vfill
\printtoc{\value{tocdepth}}
\vspace*{2ex}

\newpage
\section{Introduction and Statement of Results}
\label{sec-p2:intro}

\subsection{Motivation, Brief Overview of Results and Notation}

\subsubsection{Motivating Conjectures of \citeauthor{AD:conjecture} and \citeauthor{W:rws-hypercube}}

We analyse properties of the random walk (abbreviated \textit{RW}) on a \textit{Cayley graph} of a finite group.
The generators of this graph are chosen independently and uniformly at random.
Precise definitions are given in \S\ref{sec-p2:intro:cayley-def}. For now, let $G$ be a finite group, let $k$ be an integer (allowed to depend on $G$) and denote by $G_k$ the Cayley graph of $G$ with respect to $k$ independently and uniformly random generators.
We consider values of $k$ with $1 \ll \log k \ll \log \abs G$ for which $G_k$ is connected with high probability (abbreviated \textit{\whp}), ie with probability tending to 1 as $\abs G$ grows.

\medskip

Since pioneering work of Erd\H{o}s, it has been  understood that the typical behaviour of \emph{random} objects in some class can shed valuable light on the class as a whole.
Thus, when considering some class of combinatorial objects, it is natural to ask questions such as the following.
\begin{itemize}[noitemsep, topsep = \smallskipamount, label = \bcdot] \slshape
	\item What does a typical object in this class `look like'?
	\item If an object is chosen uniformly at random, which properties hold with high probability?
\end{itemize}
\textcite{AD:conjecture} applied this philosophy to the study of random walks on groups.

\smallskip

\textcite{AD:conjecture,AD:shuff-stop} coined the phrase \textit{cutoff phenomenon}:
	this occurs when the total variation distance (TV) between the law of the RW and its invariant distribution drops abruptly from close to $1$ to close to $0$ in a time-interval of smaller order than the mixing time.
The material in this article is motivated by a conjecture of theirs regarding `universality of cutoff' for the RW on the random Cayley graph $G_k$. It is given in \cite[Page~40]{AD:conjecture}, which is an extended version of \cite{AD:shuff-stop}.

\begin{introconj*}[\citeauthor{AD:conjecture}, \citeyear{AD:conjecture}]
	For any group $G$,
	if $k \gg \log \abs G$ and $\log k \ll \log \abs G$, then
	the random walk on $G_k$ exhibits cutoff \whp.
%	the random walk on the uniform simple Cayley graph of degree $k$ exhibits cutoff \whp.
	Further, the cutoff time, to leading order, is independent of the algebraic structure of the group: it can be written as a function only of~$k$~and~$\abs G$.
\end{introconj*}

This conjecture spawned a large body of work, including \cite{D:phd,DH:enumeration-rws,H:cutoff-cayley->,H:cutoff-cayley-survey,H:cutoff-cayley-<,R:random-random-walks,W:rws-hypercube}; see \S\ref{sec-p2:intro:previous-work}.
It has been established in the Abelian set-up
%as well as the general set-up when $\log k \gg \log \log \abs G$,
by \citeauthor{DH:enumeration-rws} \cite{DH:enumeration-rws,H:cutoff-cayley->}; see \S\ref{sec-p2:intro:previous-work:ad-conj}, and \cref{res-p2:conc-rmks:roichman} where we give a short proof.
Save \cite{H:cutoff-cayley-<} which considers the cyclic group $\mbz_p$ for prime $p$ and \cite{W:rws-hypercube} which considers $\mbz_2^d$ (which enforces $k \ge d = \log_2 \abs G$),
%Save \cite{H:cutoff-cayley-<,W:rws-hypercube},
focus has been on $k \gg \log \abs G$.

We establish cutoff for all Abelian groups when $1 \ll k \lesssim \log \abs G$ under almost optimal conditions in terms of group-generation.
We also give simple conditions under which the cutoff time is independent of the algebraic structure of the group.
%See \cref{res-p2:intro:tv}.
%See \cref{def-p2:intro:ent-time,res-p2:intro:tv}.

\medskip

The second part of this article is motivated by a conjecture of \citeauthor{W:rws-hypercube}.
\textcite{W:rws-hypercube} established cutoff for the RW on $G_k$ when $G = \mbz_2^d$ and then conjectured that $\mbz_2^d$ is the slowest amongst all groups of size at most $2^d$, asymptotically \asinf d; see \cite[Theorem~1 and Conjecture~7]{W:rws-hypercube}.

\begin{introconj*}[\citeauthor{W:rws-hypercube}, \citeyear{W:rws-hypercube}]
	For
		all diverging $d$ and $n$ with $n \le 2^d$
	and
		all groups $G$ of size $n$,
	if $k - \log_2 n \gg 1$ and $\log k \ll \log n$, then
	\(
		\tmix(\eps, G_k) / \tmix(\eps', H_k) \le 1 + \oh1 \ \whp
	\)
	for all $\eps, \eps' \in (0,1)$
	where $H \cq \mbz_2^d$---ie, the mixing time for $G_k$ is at most that of $H_k$ \whp up to smaller order terms.
\end{introconj*}

We establish a comparison between the mixing times for nilpotent and Abelian groups,
of which
%see \cref{res-p2:intro:comp:nil-abe}.
Wilson's conjecture in the nilpotent set-up is an immediate consequence.
%of this general comparison result.
%see \cref{res-p2:intro:comp:wilson}.
Additionally, we apply our nilpotent--Abelian comparison theorem to establish cutoff for various examples of non-Abelian groups, including $p$-groups with `small' commutator subgroup and Heisenberg groups.
%see \cref{res-p2:intro:comp:cor}.

\subsubsection{Brief Overview of Results}

Our focus is on mixing properties of the RW on the random Cayley graph $G_k$.
We consider the limit as $n \cq \abs G \to \infty$ under the assumption that $1 \ll \log k \ll \log \abs G$.
The condition $1 \ll \log k \ll \log \abs G$ is necessary for cutoff on $G_k$ for all nilpotent $G$; see \cref{rmk-p2:intro:cutoff:nec-con} below.

\begin{quote}
	We establish cutoff when $G$ is any Abelian group, requiring only $k - d(G) \gg 1$, where $d(G)$ is the minimal size of a generating subset of $G$.
	We show that the leading order term in the cutoff time is independent of the algebraic structure of $G$ when $d(G) \ll \log \abs G$ and $k - d(G) \asymp k$, ie it depends only on $k$ and $\abs G$.
	It is the time at which the entropy of RW on $\mbz^k$ is $\log \abs G$.
	This extends the Aldous--Diaconis conjecture to $1 \ll k \lesssim \log \abs G$.
	For certain $k$, we find the limit profile of the convergence to equilibrium.
	
	We deduce Wilson's conjecture in the Abelian set-up, as a consequence of our cutoff results.
	We then extend this to the nilpotent set-up via the following result, which is of independent interest:
		to a nilpotent group $G$, we associate an Abelian group $\widebar G$ of the same size,
		which is the direct sum of the Abelian quotients in the lower central series of $G$,
		and show that
%		the RW on $G$ mixes at least as fast as on $\widebar G$, ie
		\(
			\tmix(G_k) / \tmix(\widebar G_k) \le 1 + \oh1 \ \whp
		\)
		(provided $k - d(\widebar G) \gg 1$).
	
	We give examples where this comparison is tight:
		we establish cutoff \whp for the RW on $G_k$ where $G$ is a nilpotent group with a relatively small commutator subgroup.
		Examples of such groups include high-dimensional extra special or Heisenberg groups.
	
	Lastly, we show that the random Cayley graph of a nilpotent group $G$ is an \textit{expander} \whp whenever $k \gtrsim \log \abs G$ and $k - d(\widebar G) \asymp k$.
	(If $G$ is Abelian, then $G = \widebar G$.)
\end{quote}

Introduced by \textcite{AD:conjecture}, there has been a great deal of research into these random Cayley graphs.
Motivation for this model and an overview of historical work is given~in~\S\ref{sec-p2:intro:previous-work}.

\subsubsection{Notation and Terminology}

Cayley graphs can be either directed or undirected; we emphasise this by writing $G_k^+$ and $G_k^-$, respectively.
When we write $G_k$ or $G^\pm_k$, this means ``either $G^-_k$ or $G^+_k$'', corresponding to the undirected, respectively directed, graphs with generators chosen independently and uniformly at~random.

Conditional on being simple, $G^+_k$ is uniformly distributed over the set of all simple degree-$k$ Cayley graphs. Up to a slightly adjusted definition of \textit{simple} for undirected Cayley graphs, our results hold with $G_k$ replaced by a uniformly chosen simple Cayley graph of degree $k$; see \S\ref{sec-p2:intro:rmks:typ-simp}.% for more details.

Our results are for sequences $(G_N)_\Ninn$ of finite groups with \toinf{\abs{G_N}} \asinf N.
For ease of presentation, we write statements like ``let $G$ be a group'' instead of ``let $(G_N)_\Ninn$ be a sequence of groups''.
Likewise, the quantities $d(G)$ and, of course, $k$ appearing in the statements below all correspond to sequences, which need not be fixed (or bounded) unless we explicitly say so.
In the same vein, an event holds \textit{with high probability} (abbreviated \textit{\whp}) if its probability tends~to~1.

We use standard asymptotic notation:
	``$\ll$'' or ``$\oh{\cdot}$'' means ``of smaller order'';
	``$\lesssim$'' or $\Oh{\cdot}$'' means ``of order at most'';
	``$\asymp$'' means ``of the same order'';
	``$\eqsim$'' means ``asymptotically equivalent''.%, ie differ by smaller order terms.
%We now give our results in summarised form.
%More refined statements are given later.

\subsection{Statements of Main Results}
\label{sec-p2:intro:res}

We analyse mixing in the \textit{total variation} (abbreviated \textit{TV}) distance.
The uniform distribution on $G$, denoted $\pi_G$, is invariant for the RW.
Let $S = (S(t))_{t\ge0}$ denote the RW on $G_k$; its law is denoted $\pr[G_k]{S(t) \in \cdot}$.
For $t \ge 0$,
denote the TV distance between the law of $S(t)$ and $\pi_G$ by
\[
	d_{G_k}(t)
\cq
	\tvb{ \pr[G_k]{ S(t) \in \cdot } - \pi_G }
\cq
	\MAX{A \subseteq G}
	\absb{ \pr[G_k]{ S(t) \in A } - \abs A / \abs G }.
\]
Throughout, unless explicitly specified otherwise, we use continuous time: $t \ge 0$ means $t \in [0,\infty)$.

\subsubsection{Cutoff for All Abelian Groups}
\label{sec-p2:intro:res:tv}

\nextresult

We use standard notation and definitions for \textit{mixing} and \textit{cutoff}; see, eg, \cite[\S 4 and \S 18]{LPW:markov-mixing}.

\begin{introdefn*}
	A sequence $(X^N)_\Ninn$ of Markov chains is said to exhibit \textit{cutoff} when, in a short time-interval, known as the \textit{cutoff window}, the TV distance of the distribution of the chain from equilibrium drops from close to $1$ to close to $0$, or more precisely if there exists $(t_N)_\Ninn$ with
	\[
		\LIMSUP{\Ninf} d_N\rbb{ t_N (1 - \eps) } = 1
	\Qand
		\LIMSUP{\Ninf} d_N\rbb{ t_N (1 + \eps) } = 0
	\Qforall
		\eps \in (0,1),
	\]
	where $d_N(\cdot)$ is the TV distance of $X^N(\cdot)$ from its equilibrium distribution for each $\Ninn$.
	
	We say that a RW on a sequence of random graphs $(H_N)_\Ninn$ \textit{exhibits cutoff around time $(t_N)_\Ninn$ \whp} if,
	for all fixed $\eps$, in the limit $\Ninf$,
	the TV distance at time $(1 + \eps) t_N$ converges in distribution to $0$ and at time $(1 - \eps) t_N$ to $1$,
	where the randomness is over the random graph~$H_N$.
\end{introdefn*}

To extend the Aldous--Diaconis conjecture to $1 \ll k \lesssim \log \abs G$, one needs additional assumptions.
For an Abelian group $G$, write $d(G)$ for the minimal size of a generating set of $G$.
%; eg, $d(\mbz_2^d) = d$ and $d(\mbz_n) = 1$, and in general $d(G) \le \log_2 \abs G$.
If $k < d(G)$, then the group cannot be generated by any choice of generators.
\textcite{P:gen-abelian} shows that the expected number of independent, uniform generators required to generate the group is at most $d(G) + 3$.
(That is, if $Z_1, Z_2, ... \sim^\iid \Unif(G)$ and $\kappa \in \mbn$ is minimal with $\langle Z_1, ..., Z_\kappa \rangle = G$, then $d(G) \le \ex{\kappa} \le d(G) + 3$.)
Thus $k - d(G) \gg 1$ is always sufficient for $G$ to be generated by $\bra{Z_1^\pm, ..., Z_k^\pm}$ \whp (by Markov's inequality); we assume this throughout.
In many cases, $k - d(G) \gg 1$ is necessary to generate the group \whp, and so this assumption cannot be removed. For a characterisation of these cases and related discussion, see \cite[\cref{res-p5:gen:dichotomy}]{HOt:rcg:abe:extra}.
%(This is all for Abelian groups.)
The condition $k - d(G) \asymp k$ is particularly relevant for the Aldous--Diaconis conjecture; see \cref{rmk-p2:intro:cutoff:ad}.

\smallskip

We use an \emph{entropic method}, which involves defining \emph{entropic times}; see \S\ref{sec-p2:intro:previous-work:generic-ent} for a high-level description of the method and \S\ref{sec-p2:cutoff1:ent:method} for the specific application.
The main idea is to use an auxiliary process $W$ to generate the walk $S$; one then studies the entropy of the process $W$.
Write $Z = [Z_1, ..., Z_k]$ for the multiset of generators of the Cayley graph; then $G_k$ corresponds to choosing $Z_1, ..., Z_k \sim^\iid \Unif(G)$.
	Here, $W_i(t)$ is, for each $i$, the number of times generator $Z_i$ has been applied minus the number of times $Z_i^{-1}$ has been applied;
	$W$ is a rate-1 RW on $\mbz^k$.
	Then, $S(t) = W(t) \bcdot Z$ when the group is Abelian.
	(This auxiliary process $W$ is key even when studying nilpotent groups.)

For undirected graphs, $W$ is the usual simple RW (abbreviated \textit{SRW}):
	a coordinate is selected uniformly at random and incremented/decremented by 1 each with probability $\tfrac12$.
For directed graphs, inverses are never applied, so a step of $W$ is as follows:
	a coordinate is selected uniformly at random and incremented by 1;
	we term this the \textit{directed} RW (abbreviated~\textit{DRW}).

The \textit{entropic times} are defined in terms of the entropy of this auxiliary RW $W$.

\nextresult

\begin{subtheorem-num}{intrormkT}
	\label{rmk-p2:intro:cutoff}

\begin{introdefn}%[Entropic Time]
\label{def-p2:intro:ent-time}
	For $\gamma \in \mbn \cup \bra{\infty}$,
	let $\tent^\pm_\gamma \cq \tent^\pm_\gamma(k,G)$ be the time at which the entropy of rate-1 RW (ie, SRW or DRW, as appropriate) on $\mbz_\gamma^k$ is $\log \abs{G/\gamma G}$,
	where $\gamma G \cq \bra{ \gamma g \mid g \in G }$;
	we use the convention, $\mbz_\infty \cq \mbz$ and $\infty G \cq \abs G G = \bra{\id}$.
	Set
	\(
		\tent^\pm_*
	\cq
		\tent^\pm_*(k,G)
	\cq
		\maxt{\gamma \in \mbn \cup \bra{\infty}}
		\tent^\pm_\gamma(k,G).
	\)
\end{introdefn}

We establish cutoff for all Abelian groups, under almost optimal conditions on $k$ in terms of~$G$.
This gives an affirmative answer for Abelian groups in a strong sense to the primary part of the conjecture (occurrence of cutoff) of \textcite{AD:conjecture,AD:shuff-stop} as well as the informal question asked by \textcite{D:group-rep}; we discuss the secondary part (time depending only on $k$ and $\abs G$) in \cref{rmk-p2:intro:cutoff:ad}.

Cutoff has already been established for Abelian groups when $k \gg \log \abs G$ with $\log k \ll \log \abs G$, as mentioned above; see \S\ref{sec-p2:intro:previous-work:ad-conj}.
We thus restrict our statements to $1 \ll k \lesssim \log \abs G$.
For $1 \ll k \lesssim \log \abs G$, only two groups had been considered previously:
	$\mbz_2^d$ in \cite{W:rws-hypercube}
and
	$\mbz_p$ with $p$ prime in \cite{H:cutoff-cayley-<}.
Recall that $1 \ll \log k \ll \log \abs G$ is necessary for cutoff for nilpotent $G$, eg Abelian $G$;
%Recall that the condition $1 \ll \log k \ll \log \abs G$ is necessary for cutoff; 
% on $G_k^\pm$ for all $G$.
see \cref{rmk-p2:intro:cutoff:nec-con}.
More refined statements are given in \cref{res-p2:cutoff1:res,res-p2:cutoff2:res,res-p2:cutoff3:res}. %; see also \cref{hyp-p2:cutoff1,hyp-p2:cutoff2,hyp-p2:cutoff3}.

\begin{introthm}%[Cutoff]
\label{res-p2:intro:tv}
	Let $G$ be an Abelian group and $k$ an integer with $1 \ll k \lesssim \log \abs G$.
	Suppose that $k - d(G) \gg 1$.
	Then, the RW on $G^\pm_k$ exhibits cutoff at time
	\(
		\tent^\pm_*(k,G)
	\)
	\whp.
%	Moreover, the following~hold.
	Further, if $k - d(G) \asymp k$ and $d(G) \ll \log \abs G$, then $\tent_*(k, G) \eqsim \tent_\infty(k, \abs G)$, which depends only on $k$ and $\abs G$.
	
	Moreover, the following asymptotic relations regarding the entropic times hold.
	\begin{itemize}[noitemsep, topsep = \smallskipamount, label = \bcdot]
%		\item 
%		If $k - d(G) \asymp k$ and $d(G) \ll \log \abs G$, then $\tent_*(k, G) \eqsim \tent_\infty(k, \abs G)$.
		
		\item 
		If $k \ll \log \abs G$, then $\tent_\infty(k, \abs G) \eqsim k \abs G^{2/k} / (2 \pi e)$.
		
		\item 
		If $k - d(G) \asymp k \asymp \log \abs G$, then $\tent_*(k, G) \asymp k \abs G^{2/k} \asymp k$.
		
		\item 
		If $k > d(G)$, then $k \abs G^{2/k} \lesssim \tent_*(k, G) \lesssim k \abs G^{2/k} \log k$.
	\end{itemize}
%	Moreover, if
%		$k - d(G) \asymp k$
%	and
%		$d(G) \ll \log \abs G$,
%	then
%	\(
%		\tent_*(k,G)
%	\eqsim
%		\tent_\infty(k, \abs G)
%	\eqsim
%		k \abs G^{2/k} / (2 \pi e).
%	\)
%	If $k > d(G)$, then
%	\(
%		\tent_*(k, G)
%	\lesssim
%		k \abs G^{2/k} \log k.
%	\)
%	If $k - d(G) \asymp \log \abs G$, then
%	\(
%		\tent_*(k, G)
%	\asymp
%		k \abs G^{2/k}.
%	\)
\end{introthm}

We now give some remarks on this theorem. Further remarks are deferred to \S\ref{sec-p2:intro:res:rmks}.

\begin{intrormkt}
\label{rmk-p2:intro:cutoff:ad}
	\cref{res-p2:intro:tv} establishes cutoff for all Abelian groups, under the mild (almost necessary) condition $k - d(G) \gg 1$,
	verifying the primary part of the Aldous--Diaconis conjecture.
	Further, the secondary part is partially verified, too:
		the cutoff time depends only on $k$ and $\abs G$, up to smaller order terms, when $k - d(G) \asymp k$ and $d(G) \ll \log \abs G$.
	Cases with $k - d(G) \ll k$ or $d(G) \asymp \log \abs G$ need not satisfy this, however.
	Eg, if $k \eqsim 2 \log(4^r)$, then the groups $\mbz_2^{2r}$ and $\mbz_4^r$ give rise to mixing times which differ by a constant factor; see \cite{W:rws-hypercube} also.
	For a counterexample with $1 \ll k \ll \log \abs G$,
%	so necessarily $1 \ll k - d(G) \ll k$,
	see \cite[Proposition~\ref{res-p5:p:t0a} and \cref{res-p5:p:res}]{HOt:rcg:abe:extra} where $\mbz_p^d$, with $p$ prime, is studied.
%	both of $p$ and $d$ allowed to diverge.
\end{intrormkt}

\begin{intrormkt}
%\label{rmk-p2:intro:cutoff:profile}
	For certain regimes of $k$,
	we find the limit profile of the convergence to equilibrium:
		we define entropic times $\tent_\alpha$ and show that
		\(
			d_{G_k}(\tent_\alpha) \to^\mbp \Psi(\alpha),
		\)
		where $\Psi$ is the standard Gaussian tail; see \cref{def-p2:cutoff1:ent-times,res-p2:cutoff1:ent-times,res-p2:cutoff1:res}.
	This holds for any Abelian group if, for example,
		$k - d(G) \asymp k $ and $1 \ll k \ll \log \abs G / \log \log \log \abs G$
	or
		$k - d(G) \gg 1$ and $1 \ll k \ll \sqrt{\log \abs G / \log \log \log \abs G}$.
	The result holds for any $1 \ll k \ll \log \abs G$ under some constraints on the group.
	In \cite[\cref{res-p5:intro:cutoff}]{HOt:rcg:abe:extra}, we show the same for $k \asymp \log \abs G$, again with some constraints on~$G$.
\end{intrormkt}

\begin{intrormkt}
\label{rmk-p2:intro:cutoff:wilson-abe}
From the abstract entropic definition,
%it is not difficult to see that
$\mbz_2^d$ is the slowest amongst Abelian groups:
\[
	\max\brb{ \tent_*(k,G) \midb G \text{ Abelian group with } \abs G \le 2^d }
=
	\tent_*(k,\mbz_2^d).
\]
This verifies Wilson's conjecture in the Abelian set-up;
the general nilpotent set-up comes later.
%the more general nilpotent set-up is considered later.
	%
\end{intrormkt}

\begin{intrormkt}
\label{rmk-p2:intro:cutoff:mod-l2}
The entropic time $\tent_\infty(k, G)$ arises naturally; see \S\ref{sec-p2:cutoff1:outline} for an outline. In essence,~we~want
\[
	\mcw_t
\cq
	\brb{ w \in \mbz^k \mid \pr{W(t) = w} \ll 1/\abs G}
=
	\brb{ w \in \mbz^k \midb - \log \pr{W(t) = w} - \log \abs G \gg 1 }
\]
to satisfy $\pr{W(t) \in \mcw_t} = 1 - \oh1$.
We thus want the entropy of $W(t)$ to be at least $\log \abs G$.

The arisal of the entropic times $\tent_\mm$ ($\mm \ne \infty$) is more delicate.
We outline this in \S\ref{sec-p2:cutoff2:outline}.
\end{intrormkt}

\end{subtheorem-num}

Cutoff in $L_2$, instead of TV (ie $L_1$), can also be analysed.
For time $t \ge 0$, define
\[
	d_{G_k}^{(2)}(t)
\cq
	\normb{ \pr[G_k]{S(t) \in \cdot} - \pi_G }_{2, \pi_G}
\cq
	\rbb{ \abs G^{-1} \sumt{g \in G} \rbb{ \abs G \, \prt[G_k]{ S(t) = g } - 1 }^2 }^{1/2}.
\]
Mixing and cutoff can then be defined with respect to $L_2$ analogously to TV ($L_1$) distance.

It turns out that $L_2$ mixing time may be a constant, or even more, larger than the TV.
Similar considerations to those in \cref{rmk-p2:intro:cutoff:mod-l2} suggest that for the $L_2$ mixing the key condition is
\(
	\pr{ W(2t) = 0 }
\ll
	1/\abs G.
\)
This leads us to a conjecture for the $L_2$ mixing time, which we state informally now.
We elaborate briefly on where the proof would differ, compared with TV, in~\S\ref{sec-p2:conc-rmks:l2-mix}.
	%
%\end{intrormkt}

\begin{introconj}
\label{conj-p2:intro:l2}
	For $\mm \in \mbz \cup \bra{\infty}$,
	let $\tilde \tent^\pm_\mm \cq \tilde \tent^\pm_\mm(k,G)$ be the time $t$ at which the return probability for RW on $\mbz_\mm^k$ at time $2t$ is $\abs{G/\mm G}^{-1}$.
	Set $\tilde \tent^\pm_*(k,G) \cq \max_{\mm \in \mbn} \tilde \tent^\pm_\mm(k,G)$.
	Then, under similar conditions to those of \cref{res-p2:intro:tv}, \whp, the RW on $G_k$ exhibits cutoff in the $L_2$ metric at time~$\tilde \tent^\pm_*(k,G)$.
\end{introconj}

\nextresult

We also consider cutoff in \textit{separation distance}.
For time $t \ge 0$, define
\[
	s_{G_k}(t)
\cq
	\MAX{g \in G} \:
	\brb{ 1 - \abs G \, \pr[G_k]{ S(t) = g } }.
\]
Mixing and cutoff can then be defined with respect to separation distance analogously to TV.

It is standard that, under reversibility, the TV and separation mixing times differ by up to a factor 2; see, eg, \cite[Lemmas~6.16 and~6.17]{LPW:markov-mixing}.
However, \textcite[Theorem~1.1]{HLP:tv-sep-cutoff} showed that TV and separation cutoff are not equivalent, and that neither one implies the other.

We show that separation cutoff occurs \whp in a certain regime
%$k - d(G) \asymp k \gtrsim \log \abs G$,
and, moreover, that the cutoff time is the same, up to subleading order terms, as for TV.

A more refined statement is given in \cref{res-p2:sep:res}. %; see also \cref{hyp-p2:sep}.

\begin{introthm}%[Separation Cutoff]
\label{res-p2:intro:sep}
	Let $G$ be an Abelian group and $k$ an integer.
	Suppose that $1 \ll \log k \ll \log \abs G$ and
	\(
		k - d(G)
	\gg
		\max\bra{ (\tfrac1k \log \abs G)^2, \: (\log \abs G)^{1/2} }.
	\)
%	eg, $k - d(G) \asymp k \asymp \log \abs G$.
	Then, the RW on $G_k$ exhibits cutoff in separation distance at time
	\(
		\tent_*(k,G)
	\)
	\whp.
\end{introthm}

\begin{intrormkt}
%\label{rmk-p2:intro:sep}
	The conditions hold if $k \gtrsim (\log \abs G)^{3/4}$, $\log k \ll \log \abs G$ and $k - d(G) \gg (\log \abs G)^{1/2}$.
	Analogously to \cref{rmk-p2:intro:cutoff:wilson-abe},
	the slowest amongst Abelian groups for separation mixing is $\mbz_2^d$.
\end{intrormkt}

%\redb{make sure that the conditions are up to date}

\subsubsection{Comparison of Mixing Times Between Different Groups}
\label{sec-p2:intro:res:comp}

\nextresult

The previous results concerned cutoff.
The next results are of a slightly different flavour.
They concern \textit{nilpotent} groups:
	these are groups $G$ whose \textit{lower central series}, ie the sequence $(G_{(\ell)})_{\ell \ge 0}$ defined by $G_{(0)} \cq G$ and $G_{(\ell)} \cq [G_{(\ell-1)}, G]$ for $\ell \ge 1$, stabilises at the trivial group.
The results compare the mixing times between different groups; these mixing times are random.
%We make this concept price in the following definition.

\begin{introdefn*}
%\label{def-p2:intro:comp:mix-comparison}
	For $\eps \in (0,1)$ and a Cayley graph $H$,
	write
	\(
		\tmix(\eps, H)
	\cq
		\inf\bra{ t \ge 0 \mid d_H(t) \le \eps }.
	\)
	
%	For sequences $T \cq (T_N)_\Ninn$ of times and $H \cq (H_N)_\Ninn$ of random Cayley graphs,
%	say that
%	\textit{$\tmix(H) / T \le 1 + \oh1$ \whp}
%	if
%	there exists a non-random sequences $(\delta_N)_\Ninn$ with $\lim_N \delta_N = 0$ such that,
%	for all $\eps \in (0,1)$,
%	we have
%	\(
%		\lim_N
%		\pr{ \tmix(\eps, H_N) / T_N \le 1 + \delta_N }
%	=
%		1.
%	\)
	
	For two sequences $H \cq (H_N)_\Ninn$ and $H' \cq (H'_N)_\Ninn$ of random Cayley graphs,
	say that
	\textit{$\tmix(H) / \tmix(H') \le 1 + \oh1$ \whp}
	if
	there exist non-random sequences $(\gamma_N)_\Ninn$ and $(\delta_N)_\Ninn$ with $\lim_N \delta_N = 0$ such that,
	for all $\eps, \eps' \in (0,1)$,
	we have
	\[
		\LIM{\Ninf}
		\pr{ \tmix(\eps, H_N) \le (1 + \delta_N) \gamma_N }
	=
		1
	=
		\LIM{\Ninf}
		\pr{ (1 - \delta_N) \gamma_N \le \tmix(\eps', H'_N) }.
	\]
\end{introdefn*}

We establish Wilson's conjecture in the nilpotent set-up, as the following theorem describes.

\begin{introthm}%[Wilson's Conjecture for Nilpotent Groups]
\label{res-p2:intro:comp:wilson}
	For
		all diverging $d$ and $n$ with $n \le 2^d$
	and
		all nilpotent groups $G$ of size $n$,
	if $k - \log_2 n \gg 1$ and $\log k \ll \log n$, then
	\(
		\tmix(G_k) / \tmix(H_k) \le 1 + \oh1 \ \whp
	\)
	where $H \cq \mbz_2^d$.
\end{introthm}

As noted in \cref{rmk-p2:intro:cutoff:wilson-abe}, for Abelian groups this follows from our cutoff result
%(specifically the implicit upper bound on mixing)
and the abstract entropic definition of the cutoff time $\tent_*(k,G)$ for Abelian $G$.
The extension to nilpotent groups is then established by \cref{res-p2:intro:comp:nil-abe} below, which is of independent interest.
It is quite significantly stronger than Wilson's conjecture in the nilpotent set-up.
We can use it to establish cutoff for a class of nilpotent groups with `small commutator subgroup'; see \cref{res-p2:intro:comp:cor}.

\nextresult

\begin{subtheorem-num}{intrormkT}
	\label{rmk-p2:intro:comp}

\begin{subtheorem-num}{introcor}
	\label{res-p2:intro:comp:cor}

\begin{introthm}
\label{res-p2:intro:comp:nil-abe}
	Let $G$ be a nilpotent group.
	Set
	$\widebar G \cq \oplus_1^L \: (G_{(\ell-1)}/G_{(\ell)})$
	where $(G_{(\ell)})_{\ell\ge0}$ is the lower central series of $G$ and $L \cq \min\bra{ \ell \ge 0 \mid G_{(\ell)} = \bra{\id} }$.
	Suppose that
		$1 \ll \log k \ll \log \abs G$
	and
		$k - d(\widebar G) \gg 1$.
	Then,
	\(
		\tmix(G_k) / \tmix(\widebar G_k) \le 1 + \oh1
	\)~%
	\whp.
\end{introthm}

%This result is the key ingredient in the proof of Wilson's conjecture in the nilpotent case. Further, we demonstrate that this result is tight enough to establish cutoff for some nilpotent groups.
The quotients $G_{(\ell-1)} / G_{(\ell)} = G_{(\ell-1)} / [G_{(\ell-1)}, G]$ are all Abelian, by definition of the commutator.
The result says that the mixing time of $G_k$ is at least as fast as its Abelian counterpart $\widebar G_k$.
For a group $G$, denote its \textit{commutator} subgroup $\gcom \cq [G, G]$ and its \textit{Abelianisation}~$\gab \cq G/\gcom$.

\begin{introcor}
\label{res-p2:intro:comp:cor:gen}
	Let $G$ be a finite, non-Abelian, nilpotent group and $k$ such that $1 \ll \log k \ll \log \abs G$.
	\begin{itemize}[noitemsep, topsep = \smallskipamount, label = \bcdot]
		\item 
		If $k \eqmathsbox{abe-com}{\lesssim} \log \abs \ggr$, then
		suppose that $k \gg d(\widebar{[G,G]}) \log \abs{[G,G]}$ and $k - d(\gab) \gg d(\widebar{[G,G]})$.
		
		\item 
		If $k \eqmathsbox{abe-com}{\gg} \log \abs \ggr$, then
		suppose only that $\log \abs{[G,G]} \ll \log \abs \ggr$.
	\end{itemize}
	Then, the RW on $G_k$ exhibits cutoff at $\tent_*(k, \gab)$ \whp.
\end{introcor}

For step-2 nilpotent groups, $[G,G]$ is Abelian and hence $\overline{[G,G]} = [G,G]$.
The above corollary is thus particularly applicable for these groups.
A prime example of such groups is \textit{special groups} with small commutator subgroup.
For a prime $p$, a $p$-group is \textit{special} if it is step-2 and has centre $Z(G)$, Frattini subgroup $\Phi(G)$ and commutator subgroup $[G,G]$ all equal and elementary Abelian (ie isomorphic to $\mbz_p^s$ for some $s$).
In this case, also $\gab \cong \mbz_p^r$ where $r \cq \ell - s$ and $\ell \cq \log_p \abs G$.

We can relax the conditions on $k$ using this particular form of the Abelianisation and commutator subgroup.
The time at which the entropy of RW on $\mbz_p^k$ reaches $\log(p^r) = \log \abs{\mbz_p^r}$ is $\tent_p(k, \mbz_p^r)$

\begin{introcor}
\label{res-p2:intro:comp:cor:special}
	Let $p$ be prime, $G$ be a non-Abelian, special $p$-group and $k$ be such that $1 \ll \log k \ll \log \abs G$.
	Let $r \cq \log_p \abs \gab$, $s \cq \log_p \abs \gcom$ and $\ell \cq r + s = \log_p \abs G$.
	Suppose that $k \ge \ell$.
	\begin{itemize}[noitemsep, topsep = \smallskipamount, label = \bcdot]
		\item 
		If $k \eqmathsbox{abe-com}{\lesssim} \log \abs \ggr$, then
		suppose that $k \gg s^2 \log p$ and $k - r \gg s$.
		
		\item 
		If $k \eqmathsbox{abe-com}{\gg} \log \abs \ggr$, then
		suppose only that $s \ll r$.
	\end{itemize}
	Then, the RW on $G_k$ exhibits cutoff at
	$\tent_*(k, \gab) = \tent_p(k, \mbz_p^r)$
	\whp
	conditional that $G_k$ is connected.
	If $(k - r) p \gg 1$, then $G_k$ is connected \whp.
	If $k - r \asymp k$ and $p \gg 1$, then $\tent_p(k, \mbz_p^r) \eqsim \tent_\infty(k, \mbz_p^r)$.
\end{introcor}

Special groups are ubiquitous amongst $p$-groups of a given size in a precise, quantitative sense.
Hence, \cref{res-p2:intro:comp:cor:special} is applicable to many groups.
See \cref{rmk-p2:comp:cor:p-groups:hist} for a precise statement as well as some asymptotic expressions.
\textcite{S:enumerating-p-groups} gives, for given $(p, \ell, s)$, a simple, explicit description of all special groups of size $p^\ell$ whose commutator subgroup is of size $p^s$.

\smallskip

\textit{Extra special} groups satisfy $\gcom \cong \mbz_p$ (so $d(\gcom) = 1$) and $\abs G = p^{2d-3}$ for some integer $d \ge 3$.
For given $d$ and $p \ne 2$, up to isomorphism there are only two extra special groups.
One of these is the \textit{Heisenberg group}, which can be defined for $p$ not prime also.
For (not necessarily prime) $m, d \in \mbn$, the Heisenberg group $\HH_{m,d}$
	is the set triples $(x,y,z) \in \mbz_m^{d-2} \times \mbz_m^{d-2} \times \mbz_m$
with
%	multiplication~defined~by
	\[
		(x,y,z) \circ (x',y',z')
	\cq
		(x+x', y+y', z+z' + x \cdot y'),
	\]
	where $x \cdot y'$ is the usual dot product for vectors in $\mbz_m^{d-2}$.
We have
\(
	\HH_{m,d}^\ab
\cong
	\mbz_m^{2d-4}
\)
and
\(
	\HH_{m,d}^\com
\cong
	\mbz_m.
\)

For $p$ prime, $\HH_{p,d}$ with $d \gg 1$ falls into the class analysed in \cref{res-p2:intro:comp:cor:special} with $r = 2d-4$ and $s = 1$.
The following corollary thus focusses on $\HH_{m,d}$ with $m$ not (necessarily) prime.
Note that $\tent_\infty(k, \mbz_m^r)$ is the time at which the entropy of RW on $\mbz^k = \mbz_\infty^k$ reaches $\log(m^r) = \log \abs{\mbz_m^r}$.

\begin{introcor}
\label{res-p2:intro:comp:cor:heis}
	Let $m, d \in \mbn$ with $d \gg 1$.
	Suppose that $k - 2d \gg 1$, $k \gg \log m$ and $\log k \ll d \log m \asymp \log \abs{\HH_{m,d}}$.
	Then, \whp, the RW on $(\HH_{m,d})_k$ exhibits cutoff at $\tent_*(k, \HH_{m,d}^\ab \cong \mbz_m^{2d-4})$.
	If additionally $k - 2d \asymp k$ and $m \gg 1$, then
%	\(
%		\tent_*(k, \mbz_m^{2d-4})
%	\eqsim
%		\tent_m(k, \mbz_m^{2d-4})
%	\eqsim
%		\tent_\infty(k, \mbz_m^{2d-4})
%	\eqsim
%		k m^{2(2d-4)/k} / (2 \pi e).
%	\)
	\(
		\tent_*(k, \mbz_m^{2d-4})
	\eqsim
		\tent_\infty(k, \mbz_m^{2d-4}).
	\)
\end{introcor}

If $m$ is fixed (and thus $d \gg 1$), then the condition $k \gg \log m$ is absorbed into $k \gg 1$.
Thus, this corollary handles arbitrary Heisenberg groups $H_{m,d}$ with $m$ fixed and $k - 2d \gg 1$.

%\medskip

We now give some remarks on \cref{res-p2:intro:comp:nil-abe}
%,res-p2:intro:comp:cor}.
and Corollaries~\ref{res-p2:intro:comp:cor:gen}--\ref{res-p2:intro:comp:cor:heis}.

\begin{intrormkt}
	The bounds on $\tent_*(k,\widebar G)$, for Abelian $\widebar G$, described in \cref{res-p2:intro:tv} complement the upper bound $\tmix(G_k) \le \tmix(\widebar G_k)$ to give explicit bounds on $\tmix(G_k)$ which hold \whp.
\end{intrormkt}

\begin{intrormkt}
	In the course of proving this theorem, we prove an \emph{exact} relation between the $L_2$ mixing time for the RWs on $G_k$ and $\widebar G_k$, namely
	\(
		\ex{ d_{G_k}^{(2)}(t) } \le \ex{ d_{\widebar G_k}^{(2)}(t) }.
	\)
\end{intrormkt}

%\begin{intrormkt}
%	\color{purple}
%\label{rmk-p2:intro:comp:p-groups}
%	%
%%\cref{res-p2:intro:comp:cor} applies to many groups.
%%One class of such groups is \textit{special groups} $G$ satisfying $\log \abs \gcom \ll \log \abs G$.
%%A \textit{special group} is a step-2 $p$-group whose centre, Frattini subgroup and commutator subgroup are all equal and are elementary Abelian, ie isomorphic to $\mbz_p^s$ for some~$s$.
%%%(Heisenberg groups $H_{p,d}$ are an example of such groups.)
%%
%In the following precise sense, a large number of $p$-groups are actually special groups.
%In his classical work \cite{H:enumerating-p-groups}, \citeauthor{H:enumerating-p-groups} gave upper and lower bounds on the number groups of size $p^\ell$ for a prime $p$.
%The upper bound was later refined by \textcite{S:enumerating-p-groups}.
%A small variant of \citeauthor{H:enumerating-p-groups}'s argument shows that the logarithm of the number of such groups is dominated by special groups.
%We give further details, citing the explicit bounds derived, in \cref{rmk-p2:comp:cor:p-groups:hist}.
%
%As well as abstractly bounding the number of certain groups, \textcite{S:enumerating-p-groups} gives a simple and explicit construction of the relevant special groups.
%	%
%\end{intrormkt}

%\medskip

As explained below,
it is natural to conjecture that \cref{res-p2:intro:comp:nil-abe} does not require $G$ to be nilpotent.
The definition of the Abelian group $\widebar G$ corresponding to $G$ required $G$ to be nilpotent.
We extend this definition to allow general group $G$.
(The definitions are equivalent if $G$ is nilpotent.)

The following conjecture extends \cref{res-p2:intro:comp:nil-abe}; it contains, as a special case, Wilson's conjecture.

\begin{introconj}
\label{res-p2:intro:comp:nil-abe:conj}
	Let $G$ be a group, $(G_{(\ell)})_{\ell\ge0}$ be its lower central series and $L \cq \min\bra{ \ell \ge 0 \mid G_{(\ell)} = \bra{\id} }.$
	Let the prime decomposition of $\abs{G_L}$ be $\abs{G_L} = \prodt[r]{1} p_j$.~Set
	\(
		\widebar G
	\cq
		\rbr{ \oplus_1^L (G_{(\ell-1)}/G_{(\ell)}) }
	\oplus
		\rbr{ \oplus_1^r \: \mbz_{p_j} }.
	\)
	Suppose that $1 \ll \log k \ll \log \abs G$ and $k - d(\widebar G) \gg 1$.
	Then,
	\(
		\tmix(G_k) / \tmix(\widebar G_k) \le 1 + \oh1
	\)~%
	\whp.
\end{introconj}

We are showing in \cref{res-p2:intro:comp:nil-abe}, for nilpotent groups, that being non-Abelian can only speed up mixing.
Finite nilpotent groups are intuitively thought of as `almost Abelian'; this is partially because two elements having co-prime orders must commute.
%(Abelian groups are nilpotent groups of step~1: $G_{(0)} = G$ and $G_{(1)} = \bra{\id}$.)
Removing the nilpotent property should only mean the group is `farther from Abelian', and thus is expected to speed up mixing.

%Given \cref{res-p2:intro:comp:nil-abe}, we conjecture an analogous improvement to \cref{res-p2:intro:comp:cor:gen}.

\end{subtheorem-num} % corollaries

\end{subtheorem-num} % remarks

\subsubsection{Expander Graphs of Nilpotent Groups}

\nextresult

Our last result considers the expansion properties of the random Cayley~graph.

\begin{introdefn}
\label{def-p2:intro:isoperimetric}
	The \textit{isoperimetric constant}
%	, or \textit{Cheeger}, \textit{constant}
	of a finite $d$-regular graph $G = (V,E)$ is defined as
	\[
		\Phi_*
	\cq
		\MIN{1 \le \abs S \le \frac12 \abs V}
		\Phi(S)
	\Qwhere
		\Phi(S)
	\cq
%		\tfrac1{d \abs S}
		\absb{ \sbr{ \bra{ a, b } \in E \midb a \in S, \, b \in S^c } }
	\big/
		(d \abs S).
	\]
%	(The numerator counts the number of edges between $S$ and $S^c$, with repetitions for multigraphs.)
\end{introdefn}

%\begin{introdefn}
%	Consider a reversible Markov chain with eigenvalues $1 = \lambda_1 \ge \lambda_2 \ge \cdots \ge \lambda_n \ge -1$.
%%	 of its transition matrix.
%%	
%	The (\textit{absolute}) \textit{spectral gap} %$\gamma_*$ is defined by
%	\(
%		\gamma_*
%	\cq
%		\mint{i \ne 1} \bra{1 - \abs{\lambda_i}};
%%	=
%%		1 - \max\bra{ \abs{\lambda_2}, \: \abs{\lambda_n} };
%	\)
%	the \textit{relaxation time} %$\trel^*$ is defined by
%	\(
%		\trel^* \cq 1/\gamma_*.
%	\)
%	
%	The (\textit{absolute}) \textit{spectral gap} or \textit{relaxation time of a graph}, is that of the RW on the graph.
%\end{introdefn}

%The underlying group will be a nilpotent group, as in the previous part.
%In our cutoff results, we needed conditions like $k - d(G) \gg 1$.
%We need similar conditions here.

\begin{introthm}
\label{res-p2:intro:gap}
	Let $G$ be a nilpotent group.
	Set
	$\widebar G \cq \oplus_1^L \: (G_{(\ell-1)}/G_{(\ell)})$
	where $(G_{(\ell)})_{\ell\ge0}$ is the lower central series of $G$ and $L \cq \min\bra{ \ell \ge 0 \mid G_{(\ell)} = \bra{\id} }$.
	Then, for all $c > 0$, there exists a $c' > 0$ so that
	if $k - d(\widebar G) \ge c \log \abs G$, then
	$\Phi_*(G_k) \ge c'$ \whp.	
%	Suppose that $k - d(\widebar G) \gtrsim \log \abs G$.
%	Then $\Phi_*(G_k) \asymp 1$~\whp.
\end{introthm}

\begin{intrormkt}
\label{rmk-p2:intro:gap}
	This theorem is already known when $k - \log_2 \abs G \asymp k$, without the nilpotent restriction: it is a celebrated result of \textcite{AR:cayley-expanders}.
	For them, the constant $c'$ vanished as $k$ got closer to $\log_2 \abs G$.
	Our result removes this when $d(\overline G)$ is not close to $\log_2 \abs G$, eg $d(\overline G) \le 0.99 \log_2 \abs G$.
%	It thus suffices to consider only $k \asymp \log \abs G$.
\end{intrormkt}

%	Note also that the relaxation time is defined above for reversible Markov chains.
%	Hence we restrict attention to the \emph{undirected} random Cayley graphs $G_k^-$.

\subsubsection{Further Remarks on \cref{res-p2:intro:tv}}
\label{sec-p2:intro:res:rmks}

Here we make some remarks on \cref{res-p2:intro:tv} in addition to the three in \S\ref{sec-p2:intro:res:tv}.

\setcounter{intrormkT}{0}

\begin{subtheorem-num}{intrormkT}
	\addtocounter{intrormkT}{4}

\begin{intrormkt}
\label{rmk-p2:intro:cutoff:nec-con}
This article establishes cutoff in a variety of set-ups, but always in the regime $1 \ll \log k \ll \log \abs G$.
This leaves the regimes $k \asymp 1$ and $\log k \asymp \log \abs G$, for which there is no cutoff for any choice of generators:
	when $k \asymp 1$, this holds whenever the group is nilpotent;
	when $\log k \asymp \log \abs G$, this holds for all groups.
The former result is due to \textcite{DSc:growth-rw}; we give a short exposition of this in \cite[\S\ref{sec-p5:const-k}]{HOt:rcg:abe:extra}.
We prove the latter in \cref{res-p2:conc-rmks:roichman} below; the mixing time is order~1.
\textcite[Theorems~3.3.1 and~3.4.7]{D:phd} establishes a more general result for $\log k \asymp \log \abs G$.
\end{intrormkt}

\begin{intrormkt}
%\label{rmk-p2:intro:cutoff:lifting}
	Our approach lifts the walk $S$ from the Abelian Cayley graph $G(Z)$ to a walk $W$ on the free Abelian group with $k = \abs Z$ generators.
	Note that the walk $W$ is independent of $Z$, ie of \emph{which} $k$ generators are used.
	We study the lifted walk $W$, in particular its entropic profile, before projecting back from $W$ to $S$.
	This gives us a candidate mixing time; see \S\ref{sec-p2:intro:previous-work:generic-ent} and \S\ref{sec-p2:cutoff1:ent:method}.
\end{intrormkt}

%\smallskip

\begin{intrormkt}
%\label{rmk-p2:intro:cutoff:references}
	%
The theorem is established via two distinct approaches:
% in \S\ref{sec-p2:cutoff1} and \S\ref{sec-p2:cutoff2}, respectively.
	the former applies for $k$ not growing too rapidly;
	the second can be seen as a refinement of the first, optimised for larger $k$, where the first breaks down.
We combine the two approaches
%in \S\ref{sec-p2:cutoff3}
to analyse an interim regime of $k$.

We separate the exposition of the approaches:
	they are given in \S\ref{sec-p2:cutoff1}, \S\ref{sec-p2:cutoff2} and \S\ref{sec-p2:cutoff3}, respectively.
%	the first in \S\ref{sec-p2:cutoff1} and the second in \S\ref{sec-p2:cutoff2};
%	the third, combining the first and second, is in \S\ref{sec-p2:cutoff3}.
In the first two a concept of \textit{entropic times} is defined; see \S\ref{sec-p2:cutoff1:ent:def} and \S\ref{sec-p2:cutoff2:ent:def}.
A precise statement for each approach is given; see \S\ref{sec-p2:cutoff1:res}, \S\ref{sec-p2:cutoff2:res} and \S\ref{sec-p2:cutoff3:res}, respectively.
In summary, \cref{res-p2:intro:tv} is a direct consequence of
	Propositions~\ref{res-p2:cutoff1:ent-times} and~\ref{res-p2:cutoff2:ent:eval}
and
	\cref{res-p2:cutoff1:res,res-p2:cutoff2:res,res-p2:cutoff3:res}. %; see also \cref{hyp-p2:cutoff1,hyp-p2:cutoff2,hyp-p2:cutoff3}.
\end{intrormkt}

\end{subtheorem-num}

\subsection{Historic Overview}
\label{sec-p2:intro:previous-work}

In this subsection, we give a fairly comprehensive account of previous work on mixing and cutoff for random walk on random Cayley graphs; we compare our results with existing ones.
The occurrence of cutoff in particular has received a great deal of attention over the years.
We also mention, where relevant, other results which we have proved in companion papers;
see also \S\ref{sec-p2:intro:rmks:advert}.

%\subsubsection{Motivation: Random Cayley Graphs and Cutoff for Random Walks}
%\label{sec-p2:intro:previous-work:motivation}
%
%% Create counter `admotiv`. Set this as follows:
%%	`1` if this is Paper1, as then prints `Theorem A`;
%%	otherwise it prints `[CITATION_NUMBER, Theorem A]`
%
%%\newcounter{admotiv}
%%\setcounter{admotiv}{2}
%%
%%\input{anc/short_sections/CRC_ADMotivation_model.tex}
%%
%%\smallskip
%%
%%\input{anc/short_sections/CRC_ADMotivation_cutoff.tex}
%
%In their seminal paper, \textcite{AD:conjecture,AD:shuff-stop} considered random walks on \emph{random} Cayley graphs.
%\textcite{D:cayley:private} gave the following (paraphrased) motivation.
%\begin{quote} \slshape
%	Erd\H{o}s, when considering classes of mathematical objects, often combinatorial or graph theoretic, would often ask,
%	``What does a typical object in this class `look like'?''
%	If an object is chosen uniformly at random, are there natural properties which hold \whp?
%	
%	It is then natural to ask,
%	``How does a typical random walk on a group behave?''
%\end{quote}
%This lead \textcite{AD:conjecture,AD:shuff-stop} to consider the set of all Cayley graphs of a given group $G$ with $k$ generators.
%Drawing such a Cayley graph uniformly at random corresponds to choosing generators $Z_1, ..., Z_k \sim^\iid \Unif(G)$, conditional on giving rise to a simple graph; see \S\ref{sec-p2:intro:rmks:typ-simp}.

\subsubsection{Universal Cutoff: The Aldous--Diaconis Conjecture}
\label{sec-p2:intro:previous-work:ad-conj}

\textcite[Page~40]{AD:conjecture} stated their conjecture for $k \gg \log \abs G$.
A more refined version is given by \textcite[Conjectures~3.1.2 and~3.4.5]{D:phd}; see also \cite{H:cutoff-cayley->,R:random-random-walks}.
An informal, more general, variant was reiterated by \citeauthor{D:group-rep} in \cite[Chapter~4G, Question~8]{D:group-rep};
he gave some related open questions recently in \cite[\S 5]{D:some-things-learned}.
Towards the conjecture, an upper bound, valid for arbitrary groups, was established by \textcite[Theorem~1]{DH:enumeration-rws} and later \textcite[Theorems~1 and 2]{R:random-random-walks}, who simplified their argument.
A matching lower bound, valid only for Abelian groups, was given by \textcite[Theorem~3]{H:cutoff-cayley->}; see also \textcite[Theorem~5]{H:cutoff-cayley-survey}.
\textcite[Theorem~4]{DH:enumeration-rws} modify the proof of \cite[Theorem~3]{H:cutoff-cayley->} to extend the lower bound from Abelian groups to some families of groups with irreducible representations of bounded degree.
Combined, this established the Aldous--Diaconis conjecture for Abelian groups and such groups with low degree irreducible representations.
Moreover, the cutoff time was determined explicitly:
it is at
\[
	T(k, \abs G)
\cq
	\log \abs G / \log(k / \log \abs G)
=
	\tfrac{\rho}{\rho-1} \log_k \abs G
\Qwhere
	\text{$\rho$ is defined by $k = \rbr{ \log \abs G }^\rho$}.
\]
(To have $k \gg \log \abs G$, one needs $\rho - 1 \gg 1/\log \log \abs G$.)
See also \textcite{D:phd,H:cutoff-cayley-survey}.

There is a trivial diameter-based lower bound of $\log_k \abs G$.
If $\rho \gg 1$, ie $k$ is super-polylogarithmic in $\abs G$,
then
\(
	T(k, \abs G) \eqsim \log_k \abs G.
\)
Thus, cutoff is established for all groups~for~such~$k$.

In \cite[\cref{res-p1:intro:ext:comp}]{HOt:rcg:matrix},
using the group $\UU_{m,d}$ of $d \times d$ unit upper triangular matrices with entries in $\mbz_m$,
%we disprove the conjecture:
we disprove the part of the conjecture concerning the independence of the cutoff time from the algebraic structure of the group:
	if $d \ge 3$ is fixed and $k = (\log \abs{\UU_{m,d}})^{1 + 1/d}$,
%	if $k = \floor{\log \abs{\UU_{m,d}}}^2$ and $d = 3$,
	then there is cutoff at $\tfrac2d T(k, \abs{\UU_{m,d}})$.
In fact, $T(k,\abs{\UU_{m,d}})$ does not even capture the correct order:
	letting $d \to \infty$ sufficiently slowly, we still have $k = (\log \abs{\UU_{m,d}})^{1 + 1/d} \gg \log \abs{\UU_{m,d}}$ and the cutoff time is still shown to be $\tfrac2d T(k, \abs{\UU_{m,d}})$, which is $\oh{ T(k, \abs{\UU_{m,d}}) }$.
%	letting $d \to \infty$ sufficiently slowly, $k$ can be chosen so that $k \gg \log \abs{\UU_{m,d}}$ and there is cutoff at a time of smaller order than $T(k, \abs{\UU_{m,d}})$.

\smallskip

There has been a little investigation into the regime $1 \ll k \lesssim \log \abs G$, but with much less success.
\textcite[Theorem~4]{H:cutoff-cayley->} showed that the mixing time must be super-polylogarithmic, unlike for $k \gg \log \abs G$.
\textcite[Theorem~1]{W:rws-hypercube} established cutoff for $\mbz_2^d$; this naturally requires $k \ge d = \log_2 \abs G$.
Regarding $1 \ll k \ll \log \abs G$, a breakthrough came recently when \textcite[Theorem~1.7]{H:cutoff-cayley-<} established cutoff for $\mbz_p$ with $1 \ll k \le \log p / \log\log p$ and $p$ a (diverging) prime.
The techniques were specialised to their respective cases;
we consider arbitrary Abelian groups.

\smallskip

Relatedly, \textcite{H:random-rws-arbitrary} analysed the regime in which $k$ is just above the $\log_2 n$ threshold.
%, the minimum threshold for general results.
For $k = a \log_2 n$ with $a > 1$,
he established an upper bound on the mixing time of $a \log( a/(a-1) ) \log_2 n$ whp; see \cite[Theorem~1]{H:random-rws-arbitrary}.
However, this is quite far from tight when $a$ is large:
\[
	a \log\rbb{ a/(a-1) } \ge 1
\Quad{for all}
	a > 1,
\Quad{yet}
	T(a \log_2 n, n) / \log_2 n \eqsim 1/\log_2 a \to 0
\Quad{as}
	a \to \infty.
\]
He also analysed $k = \log_2 n + f(n)$ with $1 \ll f(n) \ll \log n$; see \cite[Theorem~2]{H:random-rws-arbitrary} particularly.

\subsubsection{Comparison of Mixing Times}

In the direction of comparison of mixing times, there has been much less work.
The only work of note (of which we are aware) is by \textcite{P:rw-few-gen}.
There, he studies universal mixing bounds (ie ones valid for all groups), but his bounds are not tight; they are always at least a constant factor away from those conjectured by \textcite{W:rws-hypercube} (and by us above).
%We are able to prove Wilson's conjecture \emph{in a specific (small) regime}: $k - \log_2 \abs G = \oh{\log \abs G}$.
%We do this by improving Pak's arguments in this regime.
%Once we have $k \ge a \log_2 \abs G$ for some constant $a > 1$, while our method does improve that of \textcite{P:rw-few-gen}, it does not match the conjectured bound.

\smallskip

A related universal bound in which $\mbz_2^d$ is the worst case is given by \textcite{P:gen-group-notes}.
Let $\varphi_k(G) \cq \pr{ G_k \ \text{is connected} }$, ie the probability that the group $G$ is generated by $k$ uniformly chosen generators.
Then, \textcite[Lecture~1, Theorem~6]{P:gen-group-notes} proves that
if $\abs G \le 2^d$ then
\(
	\varphi_k(G) \ge \varphi(\mbz_2^d)
\)
for~all~$k$.

\subsubsection{Random Walks on Upper Triangular Matrix Groups}

The study of random walks on Heisenberg groups and other groups of upper triangular matrices has a rich history.
We give a detailed historical account in \cite[\S\ref{sec-p1:intro:previous-work:rws-heisenberg}]{HOt:rcg:matrix}.

As noted above, in \cite{HOt:rcg:matrix} we study $d \times d$ unit upper triangular matrices with entries in $\mbz_m$.
By viewing the Heisenberg group as $d \times d$ matrices (see \S\ref{sec-p1:intro:not-term}), these $d \times d$ unit upper triangular matrices can be seen as a supergroup of the $d$-dimensional Heisenberg group.

\subsubsection{Expander Graphs for Nilpotent Groups}

A celebrated result of \textcite[Corollary~1]{AR:cayley-expanders} asserts that, for any finite group $G$, the random Cayley graph with at least $C_\eps \log\abs G$ random generators is \whp an $\eps$-expander, provided $C_\eps$ is a sufficiently large in terms of $\eps$.
(A graph is an \textit{$\eps$-expander} if its isoperimetric constant is bounded below by $\eps$; up to a reparametrisation, this is equivalent to having the spectral gap of the graph bounded below by $\eps$.)
There has been a considerable line of work building upon this general result of \citeauthor{AR:cayley-expanders}.
(\textcite{P:cayley-expanders} proves a similar result.)
Their proof was simplified and extended, independently, by \textcite{LS:cayley-expanders} and \textcite{LR:cayley-expanders}; both were able to replace $\log_2 \abs G$ by $\log_2 D(G)$, where $D(G)$ is the sum of the dimensions of the irreducible representations of the group $G$;
for Abelian groups $D(G) = \abs G$.
%; the latter used some Chernoff-type bounds on operator valued random variables.
A `derandomised' argument for Alon--Roichman is given by \textcite{CMR:cayley-expanders}.
Both \cite{CMR:cayley-expanders,LR:cayley-expanders} use some Chernoff-type bounds on operator valued random variables.
%These also refined the bounds on $C_\eps$ to allow any $C_\eps > 4e/\eps^2$.

\textcite{CM:cayley-expanders} improve these further by using matrix martingales and proving a Hoeffding-type bound on operator valued random variables.
They also improved the quantification for $C_\eps$, showing that one may take $C_\eps \cq 1 + c_\eps$ with $c_\eps \to 0$ as $\eps \to 0$; this means that, \whp, the graph is an $\eps$-expander whenever $k \ge (1 + c_\eps) \log_2 D(G)$ and $c_\eps \to 0$ as $\eps \to 0$.
They also generalise Alon--Roichman to random coset graphs.
The proofs use tail bounds on the (random)~eigenvalues.

It is well-known that $D(G) \ge \sqrt{\abs G}$. Thus, all these results require at least $k \ge \tfrac12 \log_2 \abs G$.
Our result, on the other hand, applies to $k \ge c \log \abs G$ for \emph{any} constant $c > 0$, provided the underlying group is suitable---eg, this is the case if $G$ is Abelian and $d(G) \ll \log \abs G$; another example is given by $d \times d$ unit upper triangular matrix groups with entries in $\mbz_m$ if $m \gg 1$.

\textcite[Theorem~1.1]{H:cutoff-cayley-<} showed, for all diverging (sequences of) primes $p$, that the order of the relaxation time of the RW on the cyclic group $\mbz_p$ is $p^{2/k}$ when $1 \ll k \le \log p / \log \log p$.

\smallskip

In \cite[\cref{res-p3:intro:gap}]{HOt:rcg:abe:geom}, we restrict to Abelian groups under the assumption $k - 2 d(G) \asymp k$ and determine, via an altogether different method, the order of the relaxation time whenever $1 \ll k \lesssim \log \abs G$: it is $\abs G^{2/k}$ \whp.
%Thus, $k \asymp \log \abs G$ and $k - 2 d(G) \asymp k$ gives relaxation time order $1$, which is equivalent to being an expander by the Cheeger inequalities.
%Further, we show for `most $G$', in a precise sense, that $k - d(G) \asymp k$ suffices.
This extends \cref{res-p2:intro:gap} in the Abelian set-up to allow $1 \ll k \ll \log \abs G$.

\subsubsection{Cutoff for `Generic' Markov Chains and the Entropic Method}
\label{sec-p2:intro:previous-work:generic-ent}

We now put our results into a broader context.
A recurrent theme in the study of mixing times is that `generic' instances often exhibit the cutoff phenomenon.
In this set-up, a family of transition matrices chosen from a certain family of distributions is shown to give rise to a sequence of Markov chains which exhibits cutoff \whp.
A few notable examples include
	random birth and death chains \cite{DW:random-matrices,S:cutoff-birth-death},
	the simple or non-backtracking random walk on various models of sparse random graphs, including
		random regular graphs \cite{LS:cutoff-random-regular},
		random graphs with given degrees \cite{BhLP:comparing-mixing,BhS:cutoff-nonbacktracking,BLPS:giant-mixing,BCS:cutoff-entropic},
		the giant component of the Erd\H{o}s--R\'enyi random graph \cite{BLPS:giant-mixing} (where the authors consider mixing from a `typical' starting point)
	and
		a large family of sparse Markov chains \cite{BCS:cutoff-entropic},
	as well as random walks on a certain generalisation of Ramanujan graphs \cite{BL:cutoff-entropic-covered} and random lifts \cite{BL:cutoff-entropic-covered,Ck:cutoff-lifts}.

A recurring idea in the aforementioned papers is that the cutoff time can be described in terms of entropy.
	One can look at some auxiliary random process which up to the cutoff time can be coupled with, or otherwise related to, the original Markov chain---often in the above examples this is the RW on the corresponding Benjamini--Schramm local limit.
	The cutoff time is then shown to be (up to smaller order terms) the time at which the entropy of the auxiliary process equals the entropy of the invariant distribution of the original Markov chain.
It is a relatively new technique, and has been used recently in \cite{BLPS:giant-mixing,BCS:cutoff-entropic,BL:cutoff-entropic-covered,Ck:cutoff-lifts}.
For `most' regimes of $k$, this is the case for us too; further, for the non-Abelian groups considered in \cite{HOt:rcg:matrix} we use a similar idea.
As our auxiliary random process, we use a SRW, respectively DRW, in the undirected, respectively directed, case.

%To the best of our knowledge, in all previous application of the entropic method the Benjamini--Schramm limit had been a tree---eg a Poisson Galton--Watson tree in \cite{BLPS:giant-mixing} and a deterministic period tree in \cite{Ck:cutoff-lifts}.
%Ours is the first application in a set-up where the graphs are not close in the local topology to a tree.

With the exception of the very recent \cite{HSS:cutoff-random-matching}, to the best of our knowledge, in all previous instances where the entropic method was used the graphs were tree-like.
This is not the case for us: in the Abelian set-up, $G_k$ has cycles of length $4$ (potentially up to the direction of edges).
Admittedly, this has less of an impact on the walk since each vertex is of diverging degree.

\subsubsection{Subsequent Work}

The release of this multi-paper project in early 2021 spurred significant interest and progress on several related problems.
\textcite{S:cutoff-nonneg-curv} established a sufficient condition for cutoff involving an entropic concentration criterion,
and verified it for RW on any undirected Abelian Cayley graph which is an expander.
The collection of Cayley graphs was extended beyond expanders by \textcite{HHPS:info-groups}.

Regarding \emph{random} Cayley graphs, \textcite{HH:rcg-nilpotent} built on the ideas initiated here, hinging on multiple aspects of the proofs of \cref{res-p2:intro:tv,res-p2:intro:comp:nil-abe}.
They established cutoff for SRW on $G_k$ for nilpotent $G$, but required bounded step and \textit{rank} $d(\gab)$.
The generality permitted by
allowing divergent rank $d(G) = d(\gab)$ for Abelian $G$ is one of the major improvements of \cref{res-p2:intro:tv} over previous work.
Also, several of our examples in \cref{res-p2:intro:comp:cor},
such as	high-dimensional Heisenberg groups $\HH_{m,d}$ ($d \gg 1$),
%for which $d(\HH_{m,d}^\ab) \asymp d \gg 1$.
are ruled out.
Under further restriction, they showed that
the mixing times
%the relaxation time
of the RW on $\ggr$ and the projection to $\gab$ are
%exactly equal, and the mixing times
asymptotically equivalent.

Even more recently, \textcite{PS:cutoff-nonneg} introduced a new criterion for cutoff for Markov chains with non-negative \textit{curvature}; RW on an undirected Abelian Cayley graph has this property.
%Cutoff can be deduced when $k \ge (1 + \eps) \log_2 n$, for which the graph is an expander whp \cite{AR:cayley-expanders}, as well as some other cases.
Again, cutoff can be deduced when the graph is an expander, but also under some weaker, quantitative conditions.
No estimate on the mixing time itself is given in \cite{S:cutoff-nonneg-curv,PS:cutoff-nonneg}, though, so connection to the motivating Aldous--Diaconis conjecture is lost.
Again, directed graphs are excluded.

%ALTERNATIVE (longer) DESCRIPTION OF JUSTIN'S RESULT
%
%Even more recently, \textcite{PS:cutoff-nonneg} introduced a new criterion for cutoff for Markov chains with non-negative \textit{curvature}; RWs on undirected Abelian Cayley graphs have this property.
%Using our lower bound on the mixing time here and our spectral gap estimate from \cite{HOt:rcg:abe:geom}, cutoff can be deduced when $1 \ll \log k \ll \sqrt{\log n}$ and $k - 2 d(G) \asymp k$, for example---not $k - d(G) \asymp k$, though.
%
%If $k - d(G) \asymp k$ and $d(G) \ll \log \abs G$, then \cref{res-p2:intro:tv} shows that the mixing time depends only on $k$ and $\abs G$.
%No estimate on the actual mixing time is given in \cite{PS:cutoff-nonneg}.
%So, the connection to the motivating Aldous--Diaconis conjecture
%is lost.
%Further, directed Cayley graphs cannot be handled.

\subsection{Additional Remarks}
\label{sec-p2:intro:rmks}

\subsubsection{Precise Definition of Cayley Graphs}
\label{sec-p2:intro:cayley-def}

Let $G$ be a finite group and $Z$ a multisubset of $G$.
We focus on mixing properties of the \textit{Cayley graph} of $G$ with \textit{generators} $Z$.
The
	\textit{undirected}, respectively \textit{directed},
\textit{Cayley graph of $G$ generated by $Z$}, denoted
	$G^-(Z)$, respectively $G^+(Z)$,
is the multigraph with vertex set $G$ and edge multiset
\[
	\sbb{ \bra{ g,g \cdot z } \mid g \in G, \, z \in Z },
\Quad{respectively}
	\sbb{ \rbr{ g,g \cdot z } \mid g \in G, \, z \in Z }.
\]
%(We do not assume that the Cayley graph is connected; that is, $Z$ may not generate $G$.)
If the walk is at $g \in G$, then a step in $G^+(Z)$, respectively $G^-(Z)$, involves choosing a generator $z \in Z$ uniformly at random and moving to $g z$, respectively one of $g z$ or $g z^{-1}$ each with probability~$\tfrac12$.

We focus attention on the \emph{random} Cayley graph defined by choosing $Z_1, ..., Z_k \sim^\iid \Unif(G)$; when this is the case, denote $G^+_k \cq G^+(Z)$ and $G^-_k \cq G^-(Z)$.
Whilst we do not assume that the Cayley graph is connected (ie, $Z$ may not generate $G$), in the Abelian set-up the random Cayley graph $G_k$ is connected \whp whenever $k - d(G) \gg 1$; see \cite[\cref{res-p5:gen:dichotomy}]{HOt:rcg:abe:extra}.
In the nilpotent set-up, this is the case whenever $k - d(G/[G,G]) \gg 1$; see \cite[\cref{rmk-p1:intro:typdist:nil-gen}]{HOt:rcg:matrix}.

%When we write $G_k$ or $G^\pm_k$, this means ``either $G^-_k$ or $G^+_k$'', corresponding to the undirected, respectively directed, graphs with generators chosen independently and uniformly at random.
%(In fact, all our results also apply when $Z_1, ..., Z_k$ are chosen \emph{without replacement} and conditioned on jointly generating the group, ie having $\langle Z_1, ..., Z_k \rangle = G$; see \S\ref{sec-p3:intro:rmks:typ-simp}.)

The graph depends on the choice of $Z$.
Sometimes, it is convenient to emphasise this;
we use a subscript, writing $\pr[G(z)]{\cdot}$ if the graph is generated by the group $G$ and multiset~$z$.
Analogously, $\pr[G_k]{\cdot}$ stands for the \emph{random} law $\pr[G(Z)]{\cdot}$ where $Z = [Z_1, ..., Z_k]$ with $Z_1, ..., Z_k \sim^\iid \Unif(G)$.

\subsubsection{Typical and Simple Cayley Graphs}
\label{sec-p2:intro:rmks:typ-simp}

The directed Cayley graph $G^+(z)$ is simple if and only if no generator is picked twice, ie $z_i \ne z_j$ for all $i \ne j$.
The undirected Cayley graph $G^-(z)$ is simple if in addition no generator is the inverse of any other, ie $z_i \ne z_j^{-1}$ for all $i,j \in [k]$.
In particular, this means that no generator is of order 2, as any $s \in G$ of order 2 satisfies $s = s^{-1}$---this gives a multiedge between $g$ and $g s$ for each $g \in G$.

The RW on $G^-(z)$ is equivalent to an adjusted RW on $G^+(z)$ where,
when a generator $s \in z$ is chosen,
instead of applying a generator $s$,
either $s$ or $s^{-1}$ is applied, each with probability $\tfrac12$.
Abusing terminology, we relax the definition of simple Cayley graphs to allow order 2 generators.
%, ie remove the condition that $z_i \ne z_i^{-1}$ for all $i$.

Given a group $G$ and an integer $k$,
we are drawing the generators $Z_1, ..., Z_k$ independently and uniformly at random.
It is not difficult to see that the probability of drawing a given multiset depends only on the number of repetitions in that multiset.
Thus, conditional on being simple, $G_k$ is uniformly distributed on all simple degree-$k$ Cayley graphs.
Since $k \ll \sqrt{\abs G}$, the probability of simplicity tends to 1 \asinf{\abs G}.
So, when we say that our results hold ``\whp (over $Z$)'', we could equivalently say that the result holds ``for almost all degree-$k$ simple Cayley graphs of $G$''.

Our asymptotic evaluation does not depend on the particular choice of $Z$, so the statistics in question depend very weakly on the particular choice of generators for almost all choices.
In many cases, the statistics depend only on $G$ via $\abs G$ and $d(G)$.
This is a strong sense of `universality'.

\subsubsection{Overview of Random Cayley Graphs Project}
\label{sec-p2:intro:rmks:advert}

This paper is one part of an extensive project on random Cayley graphs.
%At the time of writing,
There are
	three main articles \cite{HOt:rcg:abe:cutoff,HOt:rcg:matrix,HOt:rcg:abe:geom}
	(including the current one \cite{HOt:rcg:abe:cutoff}),
	a technical report \cite{HOt:rcg:abe:extra}
and
	a supplementary document \cite{HOt:rcg:supp}
	containing deferred technical proofs.
\textit{Each main article is readable independently.}
%; knowledge of one is not required to understand another.}

The main objective of the project is to establish cutoff for the random walk and determining whether this can be written in a way that, up to subleading order terms, depends only on $k$ and $\abs G$; we also study universal mixing bounds, valid for all, or large classes of, groups.
Separately, we study the distance of a uniformly chosen element from the identity, ie typical distance, and the diameter; the main objective is to show that these distances concentrate and to determine whether the value at which these distances concentrate depends only on $k$ and $\abs G$.

\begin{itemize}[noitemsep, topsep = \smallskipamount, label = \bcdot]
	\item [\cite{HOt:rcg:abe:cutoff}]
	Cutoff phenomenon (and Aldous--Diaconis conjecture) for general Abelian groups; also, for nilpotent groups, expander graphs and comparison of mixing times with Abelian groups.
	
	\item [\cite{HOt:rcg:abe:geom}]
	Typical distance, diameter and spectral gap for general Abelian groups.
	
	\item [\cite{HOt:rcg:matrix}]
	Cutoff phenomenon and typical distance for upper triangular matrix groups.
	
	\item [\cite{HOt:rcg:abe:extra}]
	Additional results on cutoff and typical distance for general Abelian groups.
	
%	\item [\cite{HOt:rcg:supp}]
%	Deferred technical results mainly regarding random walk on $\mbz$ and the volume of lattice balls.
\end{itemize}

\subsubsection{Acknowledgements}
\label{sec-p2:intro:ackno}
%\addcontentsline{toc}{section}{Acknowledgements}

This whole random Cayley graphs project has benefited greatly from advice, discussions and suggestions from many of our peers and colleagues.
We thank a few of them specifically here.

\begin{itemize}[itemsep = 0pt, topsep = \smallskipamount, label = $\bcdot$]
	\item 
	Allan Sly for suggesting the underlying entropy idea for cutoff in Approach \#1 (\S\ref{sec-p2:cutoff1}).
	
	\item 
	Justin Salez for reading this paper in detail and giving many helpful and insightful comments as well as stimulating discussions ranging across the entire random Cayley graphs project.
	
	\item 
	Evita Nestoridi and Persi Diaconis for general discussions, consultation and advice.
\end{itemize}

%\section{Abelian Groups: Cutoff; Approach \#1} \label{sec-p2:cutoff1}
\section{TV Cutoff: Approach \#1} \label{sec-p2:cutoff1}

In this section, we prove the first part of the upper bound on mixing for arbitrary Abelian groups.
The main result of the section is \cref{res-p2:cutoff1:res}. %; see also \cref{hyp-p2:cutoff1,rmk-p2:cutoff1:hyp}.
The outline of the section is as follows.

\begin{itemize}[noitemsep, topsep = 0pt, label = \bcdot]
	\item 
	\S\ref{sec-p2:cutoff1:ent:method} introduces the \textit{entropic method}.
	
	\item 
	\S\ref{sec-p2:cutoff1:ent:def} defines \textit{entropic times} and states a CLT.
	
	\item 
	\S\ref{sec-p2:cutoff1:ent:eval} sketches arguments to evaluate these entropic times.
	
	\item 
	\S\ref{sec-p2:cutoff1:res} states precisely the main theorem of the section.
	
	\item 
	\S\ref{sec-p2:cutoff1:outline} outlines the argument.
	
	\item 
	\S\ref{sec-p2:cutoff1:lower} is devoted to the lower bound.
	
	\item 
	\S\ref{sec-p2:cutoff1:upper} is devoted to the upper bound.
\end{itemize}

\subsection{Entropic Times: Methodology}
\label{sec-p2:cutoff1:ent:method}

We use an `entropic method', as mentioned in \S\ref{sec-p2:intro:previous-work}; cf \cite{BLPS:giant-mixing,BCS:cutoff-entropic,BL:cutoff-entropic-covered,Ck:cutoff-lifts}.
The method is fairly general; we now explain the specific application in a little more depth.

\smallskip

We define an auxiliary random process $(W(t))_{t\ge0}$, recording how many times each generator has been used:
	for $t \ge 0$, for each generator $i = 1,...,k$, write $W_i(t)$ for the number of times that it has been picked by time $t$.
By independence, $W(\cdot)$ forms a rate-1 DRW on $\mbz_+^k$.
For the undirected case, recall that we either apply a generator or its inverse; when we apply the inverse of generator $i$, increment $W_i \to W_i - 1$ (rather than $W_i \to W_i + 1$).
In this case, $W(\cdot)$ is a SRW on $\mbz^k$.
% (rather than a DRW on $\mbz_+^k$).

Since the underlying group is Abelian, the order in which the generators are applied is irrelevant and generator-inverse pairs cancel.
Hence, we can write
\[
	S(t)
=
	\sumt[k]{i=1} W_i(t) Z_i
=
	W(t) \bcdot Z.
\]

Recall that the uniform distribution is invariant, regardless of the group and generators.
For an Abelian group $G$, we propose as the mixing time the time at which the auxiliary process $W$ obtains entropy $\log \abs G$.
The reason for this is the following:
	take $t$ to be slightly larger than the above entropic time;
	using the equivalence $- \log \mu \ge \log \abs G$ if and only if $\mu \le 1/\abs G$,
	`typically' $W(t)$ takes values to which it assigns probability smaller than $1/\abs G$;
	informally, this means that $W(t)$ is `well spread out'.
%If we could immediately deduce that $S(t)$ typically takes values to which it assigns probability approximately $1/\abs G$, we would be basically done.
We can have two independent copies $S$ and $S'$ (using the same generators $Z$) with $S(t) = S'(t)$ but $W(t) \ne W'(t)$. The uniformity of the generators will show that, on average, this is unlikely.
We thus deduce that $S(t)$ is well spread out, ie well mixed.

Contrastingly, if the entropy is much smaller than $\log \abs G$, then $W(t)$ is not well spread out: it is highly likely to lie in a set of size $\oh{1/\abs G}$. The same must be true for $S(t)$; hence, $S(t)$~is~not~mixed.

\subsection{Entropic Times: Definition and Concentration}
\label{sec-p2:cutoff1:ent:def}

We now define precisely the notion of \textit{entropic times}.
%, and give concentration of an `entropy random variable'.
Write $\mu_t$, respectively $\nu_s$, for the law of $W(t)$, respectively $W_1(sk)$;
so $\mu_t = \nu_{t/k}^{\otimes k}$.
Define
\[
	Q_i(t) \cq - \log \nu_{t/k}\rbb{W_i(t)},
\Quad{and set}
	Q(t) \cq - \log \mu_t\rbb{W(t)} = \sumt[k]{1} Q_i(t).
\]
So $\ex{Q(t)}$ and $\ex{Q_1(t)}$ are the entropies of $W(t)$ and $W_1(t)$, respectively.
Observe that $t \mapsto \ex{Q(t)} : [0,\infty) \to [0,\infty)$ is a smooth, increasing bijection.

\begin{defn}[Entropic and Cutoff Times]
\label{def-p2:cutoff1:ent-times}
	For all $k,n \in \mbn$ and all $\alpha \in \mbr$,
	define
	\(
		\tpro_\alpha
%	\cq
%		\tpro_\alpha(k,n)
	\)
	so~that
	\[
		\ex{ Q_1\rbr{\tpro_\alpha} } = \rbb{ \log n + \alpha \sqrt{v k} }/k
	\Qand
		\spro_\alpha \cq \tpro_\alpha/k,
	\Qwhere
		v \cq \Varb{Q_1\rbr{\tpro_0}},
	\]
	assuming that $\log n + \alpha \sqrt{v k} \ge 0$.
	We call $\tpro_0$ the \textit{entropic time} and the $\bra{\tpro_\alpha}_{\alpha\in\mbr}$ \textit{cutoff times}.
\end{defn}

Comparing with notation in the introduction, $\tpro_0 = \tent_\infty(G)$; see \cref{def-p2:intro:ent-time}.
The definition there was for \emph{cutoff} only; the \emph{profile} is described by the full range $(\tpro_\alpha)_{\alpha \in \mbr}$, in the regime handled~here.

Direct calculation with the Poisson distribution and SRW on $\mbz$ gives the following relations.
%These calculations are sketched below in \S\ref{sec-p2:cutoff1:ent:eval}; rigorous arguments are given in \cite[\S\ref{sec-p0:se}]{HOt:rcg:supp}.

\begin{prop}[Entropic and Cutoff Times]
%	[{\cite[\cref{res-p0:se:t0a}]{HOt:rcg:supp}}]
\label{res-p2:cutoff1:ent-times}
Assume that $1 \ll k \ll \log n$.
Then,
for all $\alpha \in \mbr$,
\[
	\tpro_\alpha \eqsim \tpro_0 \eqsim k \cdot n^{2/k}/(2 \pi e)
\Qand
	\rbr{ \tpro_\alpha - \tpro_0 } / \tpro_0 \eqsim \alpha \sqrt{2/k} \ll 1.
%\label{eq-p2:cutoff1:ent:ent-times}
%\nt
\]
\end{prop}

The idea is to approximate the SRW and DRW laws by a normal distribution, then calculate the entropy of this.
A rigorous proof is long and tedious, requiring many careful approximations.
We sketch the principal ideas below in \S\ref{sec-p2:cutoff1:ent:eval}.
The precise details are deferred to
\cite[\cref{res-p0:se:t0a}]{HOt:rcg:supp}.%
%\cite[\S\ref{sec-p0:se}]{HOt:rcg:supp}.

Since $Q = \sumt[k]{1} Q_i$ is a sum of $k \gg 1$ iid random variables, $Q(\tpro_0)$ concentrates around $\ex{Q(\tpro_0)} = \log N$.
One can show that multiplying the time a factor $1 + \xi$ for any constant $\xi > 0$ increases the entropy by a significant amount; similarly, if $\xi < 0$, then the entropy decreases significantly.
Further, the change is by an additive term of larger order than the standard deviation $\sqrt{\Var{Q(\tpro_0)}}$.
Thus, $Q\rbr{ (1 + \xi) \tpro_0 }$ concentrates around this new value.
In particular, the following~hold:
\[
	\mu_{(1 + \xi) \tpro_0}\rbb{ W\rbb{ (1 + \xi) \tpro_0 } } &= \expb{ - Q\rbb{ (1 + \xi) \tpro_0 } }
\ll
	1/n
\quad
	\whp;
\\
	\mu_{(1 - \xi) \tpro_0}\rbb{ W\rbb{ (1 - \xi) \tpro_0 } } &= \expb{ - Q\rbb{ (1 - \xi) \tpro_0 } }
\gg
	1/n
\quad
	\whp.
\]

The following proposition quantifies this change in entropy and this concentration.
%; see \cite[\S\ref{sec-p0:se}]{HOt:rcg:supp}.

%Since the $W_i$, and hence the $Q_i$, are iid, $Q$ is a sum of $k$ iid random variables.
%Also, it turns out that
%\(
%	\Var{ Q(t) } \eqsim \Var{ Q(\tpro_0) } \gg 1
%\Qwhere
%	t \eqsim \tpro_0;
%\)
%see \cite[\cref{res-p0:se:var:t0}]{HOt:rcg:supp}.
%It then stands to reason that a CLT holds for $Q = \sumt[k]{1} Q_i$; this is indeed the case.
%The following propositions, which will be of great importance, is proved in \cite[\S\ref{sec-p0:se}]{HOt:rcg:supp}.

\begin{prop}[CLT]
%	; {\cite[\cref{res-p0:se:CLT}]{HOt:rcg:supp}}]
\label{res-p2:cutoff1:CLT}
	Assume that $1 \ll k \ll \log n$.
	For all $\alpha \in \mbr$,
	we have
	\[
		\pr{ Q(\tpro_\alpha) \le \log n \pm \omega } \to \Psi(\alpha)
	\Qfor
		\omega \cq \Varb{ Q(\tpro_0) }^{1/4} = (vk)^{1/4}.
%	\label{eq-p2:cutoff1:CLT}
%	\nt
	\]
	(There is no specific reason for choosing this $\omega$. We just need some $\omega$ with $1 \ll \omega \ll (vk)^{1/2}$.)
\end{prop}

This follows without too much difficulty from the local CLT.
Again, though, the details are technical---albeit less so than for the entropic times.
We defer the proof to \cite[\cref{res-p0:se:CLT}]{HOt:rcg:supp}

\subsection{Entropic Times: Sketch Evaluation}
\label{sec-p2:cutoff1:ent:eval}

In this subsection, we sketch details towards a proof of \cref{res-p2:cutoff1:ent-times}.
%, which finds the entropic time $\tpro_0$ and determines the window $\tpro_\alpha - \tpro_0$, all up to subleading order terms.
The full, rigorous details can be found in \cite[\cref{res-p0:se:t0a}]{HOt:rcg:supp},
where all of the approximations below are carefully justified.
%We separate the sketch into three regimes.

Recall that $\tpro_0$ is the time $t$ at which the entropy of $W_1(t)$, which is a rate-$1/k$ RW, is $(\log n)/k$.
We need to find the variance $\Var{Q_1(\tpro_0)}$, as this is used in the definition of $\tpro_\alpha$, given in \cref{def-p2:cutoff1:ent-times}.
%\[
%	\ex{ Q_1\rbr{\tpro_\alpha} } = \rbb{ \log n + \alpha \sqrt{v k} }/k
%\Qwhere
%	v \cq \Varb{Q_1\rbb{\tpro_0}}.
%\]
In the sketch below, we replace $\Var{Q_1(\tpro_0)}$ by an approximation.

\smallskip

For $s \ge 0$, denote $X_s \cq W_1(sk)$ for $s \ge 0$ and the entropy of $X_s$ as $H(s)$.
The target entropy $\log n/k \gg 1$, and so the entropic time $\spro_0 \gg 1$.
For $s \gg 1$, the RW $X_s$ has approximately the normal $N(\ex{X_s},s)$ distribution.
Translating the random variable has no affect on its entropy, and so we approximate the entropy $H(s)$ of $X_s$ by the entropy $\widebar H(s)$ of a $N(0,s)$ random variable.
Direct calculation with the normal distribution gives
\[
	\widebar H(s) = \tfrac12 \log(2 \pi e s),
\Quad{and hence}
	\widebar H{}'(s) = 1/(2s).
\]
Define $\widebar \spro_\alpha$ as the entropic times for the approximation:
\[
	\widebar H(\widebar \spro_\alpha) = \rbb{ \log n + \alpha \sqrt{v k} }/k
\Qwhere
	\widebar v \cq \Varb{\widebar Q_1\rbb{\widebar \spro_0 k}},
\]
where $\widebar Q_1(sk)$ is the analogue of $Q_1(sk)$, except with $W_1(sk)$ replaced by $N(0,s)$.
Hence,
\[
	\widebar H(\spro_0) = \log n
\Quad{implies that}
	\widebar \spro_0 = n^{2/k} / (2 \pi e) \gg 1.
\]
By direct calculation, specific to the normal distribution,
%for $s \gg 1$
one finds
\[
	\Varb{ \widebar Q_1(sk) } = \tfrac12.
\]
As mentioned above, for this sketch, to ease the calculation of $\tpro_\alpha$ in \cref{def-p2:cutoff1:ent-times}, we replace $\Var{Q_1(\tpro_0)}$ by its approximation $\tfrac12$, and assume the above normal distribution approximation.

%\smallskip

In order to find the window, assuming for the moment that $\alpha > 0$, we write
\[
%\textstyle
	\spro_\alpha - \spro_0
=
	\int_0^\alpha
	\frac{d\spro_a}{da} \, da.
\]
Again, we replace $\spro_\alpha$ with $\widebar \spro_\alpha$.
By definition, $\widebar \spro_\alpha$ satisfies
\[
	\widebar H(\widebar \spro_\alpha) = \log n/k + \alpha / \sqrt{2k},
\Quad{and hence}
%\textstyle
	\frac{d\widebar \spro_\alpha}{d\alpha} \widebar H{}'(\widebar \spro_\alpha) = 1/\sqrt{2k}.
\]
Using the expressions for $d\widebar \spro_a/da$ and $\widebar H{}'(s) = 1/(2s)$ above, we find that
\[
	\widebar \spro_\alpha - \widebar \spro_0
=
	(2k)^{-1/2} \intt[\alpha]{0} 2 \widebar \spro_a \, da
\approx
	(2k)^{-1/2} \intt[\alpha]{0} 2 \widebar \spro_0 \, da
=
	\alpha \widebar \spro_0 \sqrt{2/k},
\]
since $\widebar \spro_a$ only varies by subleading order terms over $a \in [0,\alpha]$.
The argument is analogous for $\alpha < 0$.

\smallskip

We have now shown the desired result for $\widebar \spro_\alpha$, ie when approximating $W_1(sk)$ by $N(\ex{X_s},s)$.
It turns out that this approximation is sufficiently good for the results to pass over to the original case, ie to apply to $\spro_0$ and $\tpro_0 = \spro_0 k$.
This is made rigorous in \cite[\S\ref{sec-p0:se}]{HOt:rcg:supp} via a local CLT.

\subsection{Precise Statement and Remarks}
\label{sec-p2:cutoff1:res}

In this subsection, we state precisely the main theorem of the section.
There are some simple conditions on $k$, in terms of $d(G)$ and $\abs G$, needed for the upper bound.

\begin{hyp}
\label{hyp-p2:cutoff1}
	The sequence $(k_N, G_N)_\Ninn$ satisfies \textit{\cref{hyp-p2:cutoff1}} if
%	each individually is a non-decreasing sequences of positive integers and
	the following hold:
	\begin{gather*}
		\LIM\Ninf \abs{G_N} = \infty,
%	\quad
%		\LIMSUP\Ninf k_N / \log \abs{G_N} < \infty,
	\quad
		\LIM\Ninf \rbb{ k_N - d(G_N) } = \infty
	\quad
		\text{and}
	\\
		\frac{k_N - d_N(G_N) - 1}{k_N} \ge 5 \frac{k_N}{\log \abs{G_N}} + 2 \frac{d_N(G_N) \log\log k_N}{\log \abs{G_N}}
	\text{ for all }
		\Ninn.
	\end{gather*}
\end{hyp}

\begin{rmkt}
\label{rmk-p2:cutoff1:hyp}
Write $n \cq \abs G$.
Any of the following conditions imply \cref{hyp-p2:cutoff1}:
%The following conditions are sufficient for \cref{hyp-p2:cutoff1} to hold:
\begin{alignat*}{2}
	1 \ll k &\lesssim \sqrt{\log n / \log\log\log n}
&
	\Qand
	k - d &\gg 1;
\\
	1 \ll k &\lesssim \sqrt{\log n}
&
	\Qand
	k - d &\gg \log\log k;
\\
	1 \ll k &\ll \log n/\log\log\log n
&
	\Qand
	k - d &\ge \delta k
	\Quad{for some suitable} \delta = \oh1;
\\
	d &\ll \log n/\log\log\log n
&
	\Qand
	k - d &\asymp k \ll \log n.
\tag*{\qedhere}
\end{alignat*}
\end{rmkt}

%In \cref{rmk-p2:cutoff1:hyp} below, we give some sufficient conditions of \cref{hyp-p2:cutoff1} to hold.
Throughout the proofs, we drop the subscript-$N$ from the notation, eg writing $k$ or $n$, considering sequences implicitly.
Recall that we abbreviate the TV distance from uniformity at time~$t$~as
\[
	d_{G_k, N}(t) = \tvb{ \pr[{G_N([Z_1, ..., Z_{k_N}])}]{ S(t) \in \cdot } - \pi_{G_N} }
\Qwhere
	Z_1, ..., Z_{k_N} \sim^\iid \Unif(G_N).
\]

We now state the main theorem of this section.
Recall that $\Psi$ is the standard Gaussian tail.

\begin{thm}%[Cutoff: Upper Bound \#1]
\label{res-p2:cutoff1:res}
	Let $(k_N)_\Ninn$ be a sequence of positive integers and $(G_N)_\Ninn$ a sequence of finite, Abelian groups;
	for each $\Ninn$, define $Z_{(N)} \cq [Z_1, ..., Z_{k_N}]$ by drawing $Z_1, ..., Z_{k_N} \sim^\iid \Unif(G_N)$.
	
	Suppose that the sequence $(k_N, G_N)_\Ninn$ satisfies \cref{hyp-p2:cutoff1}.
	For all $\alpha \in \mbr$ and all $N \in \mbn$,
	write $\tpro_{\alpha,N} \cq \tpro_\alpha(k_N, \abs{G_N})$.
	Let $\alpha \in \mbr$.
	Then,
	\[
		\tpro_{\alpha,N} / \tpro_{0,N}
	\to
		1
	\Qand
		d_{G_k, N}\rbr{ \tpro_{\alpha,N} }
	\to^\mbp
		\Psi(\alpha)
	\Quad{(in probability)}
		\asinf N.
	\]
	That is, \whp there is TV cutoff at $\tpro_0$ with profile given by $\bra{\tpro_\alpha}_{\alpha \in \mathbb R}$:
		for all $\eps \in (0,1)$,
		the difference in the mixing times $\tmix(\eps) - \tmix(\tfrac12)$ is given, up to smaller order terms, by $\tpro_{\Psi^{-1}(\eps)} - \tpro_0$.
	Moreover, the implicit lower bound on the TV distance holds deterministically, ie for all choices of generators.
\end{thm}

\begin{rmkt*}
Using \cref{res-p2:cutoff1:ent-times},
we can write the cutoff statement in the form
\[
	\rbb{ \tmix(\eps) - \tpro_0 } / w
\to^\mbp
	\Psi^{-1}(\eps)
\quad
	\whp
\Qforall
	\eps \in (0,1),
\]
where $\tpro_0 \eqsim k \abs G^{2/k} / (2 \pi e)$ is the mixing time and $w \eqsim \sqrt k \abs G^{2/k} / (\sqrt2 \pi e)$ the window.
\end{rmkt*}

\begin{rmkt*}
The CLT, \cref{res-p2:cutoff1:CLT}, will give the dominating term in the TV distance:
%\cref{eq-p2:cutoff1:res}:
\begin{itemize}[noitemsep, topsep = 0pt, label = \bcdot]
	\item 
	on the event $\bra{Q(\tpro_\alpha) \le \log n - \omega}$, we lower bound the TV distance by $1 - \oh1$;

	\item 
	on the event $\bra{Q(\tpro_\alpha) \ge \log n + \omega}$, we upper bound the expected TV distance by $\oh1$.
\end{itemize}
Combining this with the CLT, we deduce that $d_{G_k}(\tpro_\alpha) \to \Psi(\alpha)$ in probability.
\end{rmkt*}

\begin{rmkt*}
	Observe that \cref{hyp-p2:cutoff1} does not cover the regime $k \gtrsim \log \abs G$.
	Under fairly mild conditions on the group we can apply a variation on the argument given below to obtain a limit profile result for any $1 \ll \log k \ll \log \abs G$.
	The detailed analysis is carried out in
%	We do not carry out the analysis here; see
	\cite[\S\ref{sec-p5:profile}]{HOt:rcg:abe:extra}.
\end{rmkt*}

\subsection{Outline of Proof}
\label{sec-p2:cutoff1:outline}

\renewcommand{\mm}{\ensuremath{\gamma}}

We now give a high-level description of our approach, introducing notations and concepts along the way.
No results or calculations from this section will be used in the remainder of the document.
%rather, this section merely introduces ideas.
%Recall the definitions from the previous sections.
Further, we restrict attention to establishing \emph{cutoff} only, not the \emph{limit profile}: take $t = (1 \pm \eps) \tpro_0$. 

In all cases, we show that cutoff occurs around the entropic time.
As $Q(t)$ is a sum of many iid random variables, we expected it to concentrate around its mean.
Loosely speaking, we show that the shape of the cutoff, ie the profile of the convergence to equilibrium, is determined by the fluctuations of $Q(t)$ around its mean, which in turn, by the CLT (\cref{res-p2:cutoff1:CLT}), are determined by $\Var{Q(t)}$, for $t$ `close' to $\tpro_0$.
Note that $\Var{Q(t)} = k \Var{Q_1(t)}$ since the $Q_i$ are iid.

Throughout this section (\S\ref{sec-p2:cutoff1:outline}),
we write $0$ for the identity element of the Abelian group $G$,
as is standard.
We now outline the proof in more detail.
We often drop $t$-dependence from the notation.

\medskip

We start by discussing the lower bound.
If $Q$ is sufficiently small, then $W$, and hence also $S$, is restricted to a small set.
Indeed,
\(
	Q \le \log n - \omega
\Quad{if and only if}
	 \mu\rbr{W} \ge n^{-1} e^\omega,
\)
and thus if this is the case then
\(
	W \in \bra{ w \mid \mu(w) \ge n^{-1} e^\omega }.
\)
Since we generate $S$ via $W$, it is thus also the case that
\[
	S \in E \cq \bra{ g \in G \mid \pr{ S = g } \ge n^{-1} e^\omega }.
\]
But clearly $\abs{E} \le n e^{-\omega}$.
Choosing the time $t$ slightly smaller than the entropic time $\tpro_0$ and $\omega \gg 1$ suitably, the event
\(
	\bra{ Q(t) \le \log n - \omega }
\)
will hold \whp.
Thus, \whp, $S(t)$ is restricted to a set of size $\oh{n}$.
Hence, it cannot be mixed.
This heuristic applies for \emph{any} choice of generators.

Precisely,
we show
%We show in \S\ref{sec-p2:cutoff1:lower},
	for any $\omega$ with $1 \ll \omega \ll \log n$,
	all $t$
and
	all $Z = [Z_1, ..., Z_k]$,
that
\[
	d_{G_k}(t)
\ge
	\pr{ Q(t) \le \log n - \omega } - e^{-\omega}.
\]
%Observe that the probability on the right-hand side is independent of $Z$.
Thus, we are interested in the fluctuations of $Q(t)$ for $t$ close to $\tpro_0$.
Using the CLT application above, ie \cref{res-p2:cutoff1:CLT} with $\omega \cq \Var{Q(\tpro_0)}^{1/4}$, we deduce the lower bound in \cref{res-p2:cutoff1:res}.

\medskip

We now turn to discussing the upper bound.
%Our aim is to show that, for all $\alpha \in \mbr$, we have
%\[
%	d_{G_k}(\tpro_\alpha) \le \Psi(\alpha) + \oh1
%\quad
%	\text{\whp over $Z$}.
%\]
The lower bound was valid for any choice of generators $Z$.
Here, we exploit the independence and uniformity of the elements of $Z$.
%For clarity of presentation, we concentrate here on $G = \mbz_n$.
%; we consider general Abelian groups later.

Let $W'(t)$ be an independent copy of $W(t)$, and let $V(t) \cq W(t) - W'(t)$.
Observe that, in both the undirected and directed case, the law of $V(t)$ is that of the rate-2 SRW in $\mbz^k$, evaluated at time $t$.
%For now, we suppress the $t$ from the notation.
%It is standard that the TV distance $\tv{\zeta - \pi_G}$ can be upper bounded by half the $L_2$ distance:
The standard $L_2$ calculation (using Cauchy--Schwarz) says that
\[
	2 \, \tv{ \zeta - \pi_G }
\le
	\norm{ \zeta - \pi_G }_2
=
	\sqrt{ n \sumt{x \in G} \rbb{ \zeta(x) - \tfrac1n }^2 },
\]
recalling that $\pi_G(x) = 1/n$ for all $x \in G$.
A standard, elementary calculation shows that
\[
	\normb{ \pr[G_k]{ S(t) \in \cdot } - \pi_G }_2
=
	\sqrt{ n \, \pr{ V(t) \bcdot Z = 0 \mid Z } - 1 }.
\]
Unfortunately, writing $X = (X(s))_{s\ge0}$ for a rate-1 SRW on $\mbz$, a simple calculation shows that
\[
	\pr{ V(\tpro_0) \bcdot Z = 0 \mid Z }
\ge
	\pr{ V(\tpro_0) = (0,...,0) \in \mbz^k }
=
	\pr{ X(2\tpro_0/k) = 0 }^k
\gg
	1/n.
\]
(This calculation differs amongst the regimes of $k$.)
Moreover, the $L_2$-mixing time can then be shown to be larger than the TV-mixing time by at least a constant factor;
hence, this is insufficiently precise for showing cutoff in TV.
We drop the $t$-dependence from the notation from now on.

This motivates the following `modified $L_2$ calculation'.
First, let $\mcw \subseteq \mbz^k$, and write
\[
	\typ \cq \brb{ W, W' \in \mcw },
\quad
	\widebar{\mbp}\rbr{\cdot} \cq \prt{ \, \cdot \mid \typ }
\Qand
	\widebar{\mbe}\rbr{\cdot} \cq \ext{ \, \cdot \mid \typ };
\]
note that here we are (implicitly) averaging over $Z$.
The set $\mcw \subseteq \mbz^k$ will be chosen later so that
\[
	\widebar{\mbp}\rbb{ V = 0 } = \pr{ V = 0 \mid \typ } \ll 1/n
\Qand
	\pr{ W \notin \mcw } = \oh1;
\]
we call this \textit{typicality}.
We now perform the same $L_2$ calculation, but for $\widebar{\mbp}$ rather than $\mbp$:
\begin{gather*}
	d_{G_k}(t)
=
	\tvb{ \pr[G_k]{ S \in \cdot } - \pi_G }
\le
	\tvb{ \pr[G_k]{ S \in \cdot \mid W \in \mcw } - \pi_G }
+	\pr{ W \notin \mcw };
\\
	4 \, \ex{ \tvb{ \pr[G_k]{ S \in \cdot \mid W \in \mcw } - \pi_G }^2 }
\le
	\ex{ \abs G \, \widebar{\mbp}\rbb{ V \bcdot Z = 0 \mid Z } - 1 }
=
	\abs G \, \widebar{\mbp}\rbb{ V \bcdot Z = 0 } - 1;
\end{gather*}
%\[
%	\tvb{ \pr{ S \in \cdot \mid Z } - \pi_G }
%&\le
%	\tvb{ \pr{ S \in \cdot \mid Z, \, W \in \mcw } - \pi_G }
%+	\pr{ W \notin \mcw };
%\\
%	2 \, \ex{ \tvb{ \pr{ S \in \cdot \mid Z, W \in \mcw } - \pi_G } }
%&\le
%	\ex{ \sqrt{ n \, \widebar{\mbp}\rbb{ V \bcdot Z \equiv 0 \mid Z } - 1} }
%\le
%	\sqrt{ n \, \widebar{\mbp}\rbb{ V \bcdot Z \equiv 0 } - 1 },
%\]
see \cref{res-p2:cutoff1:mod-l2}.
By taking expectation over $Z$ and doing a modified $L_2$ calculation, we transformed the quenched estimation of the mixing time into an annealed calculation concerning the probability that a random word involving random generators is equal to the identity.
This is a key step.

%We think of $\mcw$ as a set of `typical values' for $W$.
To have $w \in \mcw$, we impose \textit{local} and \textit{global typicality requirements}.
The \emph{global} part says that
\[
	- \log \mu(w) \ge \log n + \omega
\Qforall
	w \in \mcw,
\]
where $\omega \cq (v k)^{1/4}$ as above;
the \emph{local} part will come later.
For a precise statement of the typicality requirements, see \cref{def-p2:cutoff1:typ}.
These have the property that, when $t = (1 + \eps) \tpro_0$,
\[
	\pr{ W(t) \notin \mcw }
\ll
	1;
%\Qwhen
%	t = \tpro_\alpha;
\]
see Proposition~\ref{res-p2:cutoff1:typ}.
%Consider $\alpha$ large so that $\Psi(\alpha) \approx 0$.
Then, since $-\log p \ge \log n + \omega$ if and only if $p \le n^{-1} e^{-\omega}$, we have
\[
%	\pr{ V = (0,...,0) \mid W,W' \in \mcw }
	\widebar{\mbp}\rbb{ V = (0,...,0) }
\asymp
	\pr{ W = W' \mid W' \in \mcw }
\le
	n^{-1} e^{-\omega}.
\]
%Since $\pr{W \in \mcw} \asymp 1$, this implies that
%\[
%	\widebar{\mbp}\rbb{ V = (0,...,0) }
%\lesssim
%	n^{-1} e^{-\omega}.
%\]

\smallskip

Of course, there are other scenarios in which we may have $V \bcdot Z \equiv 0$.
To deal with these,
since linear combinations of independent uniform random variables in an Abelian group are uniform on their support,
%conditional on $V = v$,
we have $v \bcdot Z \sim \Unif(\mfgcd_v G)$ where $\mfgcd_v \cq \gcd(v_1, ..., v_k, n)$; see \cref{res-p2:cutoff1:unif-gcd}.
%(For an Abelian group $G$ and $\mm \in \mbn$,
%define
%\(
%	\mm G
%\cq
%	\brb{ \mm g \mid g \in G };
%\)
%eg, $2 \mbz_{2(m+1)} = \bra{ 0, 2, 4, ..., 2(m-1), 2m }$.)
Then,
\[
	\abs G \, \widebar{\mbp}\rbb{ V \bcdot Z = 0, \: V \ne 0 }
=
	\abs G \, \widebar{\mbe}\rbb{ \one{V \ne 0} / \abs{ \mfgcd_V G} } .
\]
%with the last inequality being elementary ().
(Recall that $V$ and $Z$ are independent.)
We use the \emph{local} typicality conditions to ensure that $\max_i \abs{W_i} \le r_*$, for some explicit $r_*$ which diverges a little faster than $n^{1/k}$.
This allows us to consider only values $\mfgcd \in [2 r_*]$ for the gcd.
%We then need to bound the expectation of the gcd.
It is here where the two approaches (\S\ref{sec-p2:cutoff1} and \S\ref{sec-p2:cutoff2}) diverge.

First (\S\ref{sec-p2:cutoff1}), we use a rather direct approach.
Elementary group theory gives
\[
	\abs G \, \widebar{\mbe}\rbb{ \one{V \ne 0} / \abs{ \mfgcd_V G} }
\le
	\widebar{\mbe}\rbb{ \mfgcd_V^{d(G)} \one{V \ne 0} }
\le
	1
+	\sumt[2r_*]{\mm=2}
	\mm^{d(G)} \pr{ \mfgcd_V = \mm };
\]
see \cref{res-p2:cutoff1:G/gammaG}.
Since the law of SRW on $\mbz$ is unimodal,
for each non-zero coordinate,
the probability that $\mm$ divides it is at most $1/\mm$.
Thus, in general, the probability is at most $1/\mm$ plus the probability that the coordinate is 0, the latter of which is order $1/\sqrt{t/k}$.
This leads to
\[
	\widebar{\mbp}\rbr{ \mfgcd_V = \mm }
\lesssim
	\rbb{ 2 / n^{1/k} + 1/\mm }^k
\Qwhen
	t \ge \tpro_0;
\]
see \cref{res-p2:cutoff1:divis}.
Provided at least one of $d(G)$ or $k$ is not too close to $\log n$,
we are able to use this to control the expectation, showing
\(
	\widebar{\mbe}\rbr{ \mfgcd_V^{d(G)} \one{V \ne 0} }
=
	1 + \oh1
\)
when $t \ge \tpro_0$;
%showing that it is $1 + \oh1$;
see \cref{res-p2:cutoff1:gcd-ex}.

Combining these two analyses,
we deduce that
\[
	n \, \widebar{\mbp}\rbr{ V \bcdot Z = 0 }
\le
	n \, \widebar{\mbp}\rbr{ V \bcdot Z = 0, \: V \ne 0 }
+	n \, \widebar{\mbp}\rbr{ V = 0 }
=
	1 + \oh1.
\]
The modified $L_2$ calculation then says that the TV distance tends to $0$ in probability, as required.

The only real adjustment needed to handle the profile is the use of the estimate
\[
	\pr{ W(\tpro_\alpha) \notin \mcw }
\eqsim
	\Psi(\alpha).
\]
The remainder of the analysis is fairly robust to the specific value of $t$.
%The modified $L_2$ calculation then says that the TV distance is roughly $\Psi(\alpha)$ plus a term $o_\mbp(1)$, ie tending to 0 in probability.
%This establishes a matching limiting upper bound of $\Psi(\alpha)$ in probability.

\smallskip

The second approach (\S\ref{sec-p2:cutoff2}) analyses the term $\widebar{\mbp}\rbr{ \mfgcd_V = \mm }$ and uses it to kill $\abs{G/\mm G}$ directly in
\[
	\abs{G} \, \widebar{\mbe}\rbb{ \one{V \ne 0} / \abs{ \mfgcd_V G} }
=
	\sumt{\mm \in \mbn}
	\widebar{\mbp}\rbr{ \mfgcd_V = \mm } \abs{G/\mm G}.
\]
We outline the details of the adaptation in \S\ref{sec-p2:cutoff2:outline}, including where Approach \#1 breaks down.

\medskip

This concludes the outline.
We now move onto the formal proofs.

\subsection{Lower Bound on Total-Variation Mixing}
\label{sec-p2:cutoff1:lower}

In this subsection,
we prove the lower bound on mixing, which holds for all choices of generators.

\begin{Proof}[Proof of Lower Bound in \cref{res-p2:cutoff1:res}]
For this proof only, to emphasise that $Z$ is fixed, not being averaged over, we add a subscript-$Z$ to the probabilities involving $Z$: $\pr[Z]{ S(\tpro_\alpha) \in \cdot }$.
%Those involving only $W$ are unaffected, of course, by $Z$.
%For this proof,
%Assume that $Z$ is given, and suppress it from notation.

For all $\alpha \in \mbr$, by the CLT (\cref{res-p2:cutoff1:CLT}),
\[
	\prt{\mce_\alpha}
\eqsim
	\Psi(\alpha)
\Qwhere
	\mce_\alpha
\cq
	\brb{ \mu\rbb{W(\tpro_\alpha)} \ge n^{-1} e^\omega }
=
	\brb{ Q(\tpro_\alpha) \le \log n - \omega };
\]
recall that $\omega \gg 1$.
Fix $\alpha \in \mbr$.
Consider the set
\[
	E_\alpha
\cq
	\brb{ x \in G \midb \exists \, w \in \mbz^d \st \mu_{\tpro_\alpha}(w) \ge n^{-1} e^\omega \text{ and } x = w \bcdot Z }.
\]
Then, $\pr[Z]{ S(\tpro_\alpha) \in E_\alpha \mid \mce_\alpha } = 1$ since $W$ generates $S = W \bcdot Z$.
Every element $x \in E_\alpha$ can be realised as $x = w_x \bcdot Z$ for some $w_x \in \mbz^k$ with $\mu_{\tpro_\alpha}(w_x) \ge n^{-1} e^\omega$.
Hence, for all $x \in E_\alpha$, we have
\[
	\pr[Z]{ S(\tpro_\alpha) = x }
\ge
	\pr{ W(\tpro_\alpha) = w_x }
=
	\mu_{\tpro_\alpha}(w_x)
\ge
	n^{-1} e^\omega.
\]
Taking the sum over all $x \in E_\alpha$,
we deduce that
\[
	1
\ge
	\sumt{x \in E_\alpha}
	\pr[Z]{ S(\tpro_\alpha) = x }
\ge
	\abs{E_\alpha} \cdot n^{-1} e^\omega,
\Quad{and hence}
	\abs{E_\alpha}/n \le e^{-\omega} = \oh1.
\]
Finally, we deduce the lower bound from the definition of TV distance:
\[
	\tvb{ \pr[Z]{ S(\tpro_\alpha) \in \cdot } - \pi_G }
\ge
	\pr[Z]{ S(\tpro_\alpha) \in E_\alpha } - \pi_G(E_\alpha)
\ge
	\prt{ \mce_\alpha } - \tfrac1n \abs{E_\alpha}
\ge
	\Psi(\alpha) - \oh1.
\qedhere
\]
\end{Proof}

\begin{rmkt*}
	Given an arbitrary group $G$, projecting the walk from $G$ to the Abelianisation $\gab = G / [G,G]$, which is Abelian, cannot increase the TV distance.
	Thus, $\tpro_0(k, \abs\gab)$ is a lower bound on mixing for the projected walk on the $\gab$, and hence for the original walk on $G$ too.
\end{rmkt*}

%\begin{rmkt*}
%	Using a variant of this argument, in \cite[\S\ref{sec-p1:cutoff:lower}]{HOt:rcg:matrix} we prove an analogous lower bound for general groups:
%		where $\tpro_0(k,\abs G)$ was the lower bound above (for Abelian groups),
%		we establish a lower bound of $\tpro_0(k,\abs{G/[G,G]})$ for any group.
%	(If a group is Abelian, then $[G,G]$ is trivial.)
%	In many cases, this is a significant improvement over previous best-known bound of $\log_{k-1} \abs G$.
%\end{rmkt*}

\subsection{Upper Bound on Total-Variation Mixing}
\label{sec-p2:cutoff1:upper}
\renewcommand*{\mm}{\gamma}

It is often easier to control $L_2$ distances, rather than $L_1$ (ie, TV).
However, $L_2$ is sensitive to rare events, unlike TV.
We use a `modified $L_2$ calculation' to bound the TV:
	first, condition that $W$ is `typical';
	then, use a standard $L_2$ calculation on the conditioned law.
Let $W'$ be an independent copy of $W$.
Then, $S' \cq W' \bcdot Z$ has the same law as $S$ and is conditionally independent of $S$ given~$Z$.
%then $S' \cq W' \bcdot Z$ is an independent copy of $S$ given $Z$.

\begin{lem}%[Modified $L_2$ Calculation]
\label{res-p2:cutoff1:mod-l2}
	For all $t \ge 0$ and all $\mcw \subseteq \mbz^k$,
	the following inequalities hold:
%	\begin{subequations}
%		\label{eq-p2:cutoff1:mod-l2}
	\begin{gather*}
		d_{G_k}(t)
	=
		\tvb{ \pr[G_k]{ S(t) \in \cdot } - \pi_G }
	\le
		\tvb{ \pr[G_k]{ S(t) \in \cdot \mid W(t) \in \mcw } - \pi_G }
	+	\pr{ W(t) \notin \mcw };
%	\label{eq-p2:cutoff1:mod-l2:triangle}
%	\nt
	\\
		4 \, \ex{ \tvb{ \pr[G_k]{ S(t) \in \cdot \mid W(t) \in \mcw } - \pi_G }^2 }
	\le
		n \, \pr{ S(t) = S'(t) \mid W(t), W'(t) \in \mcw } - 1.
%	\label{eq-p2:cutoff1:mod-l2:tv-l2}
%	\nt
	\end{gather*}
%	\end{subequations}
	We emphasise that $d_{G_k}$ is a random variable, a function of $Z_1, ..., Z_k \sim^\iid \Unif(G)$.
\end{lem}

%\begin{Proof}
%	The first claim follows from the triangle inequality; the second from Cauchy--Schwarz.
%\end{Proof}

\begin{Proof}
The first claim follows immediately from the triangle inequality.
For the second, using Cauchy--Schwarz, we upper bound the TV distance of the conditioned law by its $L_2$ distance:
\[
&
	4 \, \tvb{ \pr[G_k]{ S \in \cdot \mid W \in \mcw } - \pi_G }^2
%&
\le
	n \sumt{x} \rbb{ \pr[G_k]{ S = x \mid W \in \mcw } - \tfrac1n }^2
\\&
\qquad
=
	n \sumt{x} \pr[G_k]{ S = x \mid W \in \mcw }^2 - 1
%\\&
%\qquad
=
	n \sumt{x} \pr[G_k]{ S = S' = x \mid W,W' \in \mcw } - 1,
\]
as $S = W \bcdot Z$ and $S' = W' \bcdot Z$.
The claim follows by taking expectations over $[Z_1, ..., Z_k]$.
\end{Proof}

We now make the specific choice of the `typical' set $\mcw$; we make a different choice for each $\alpha \in \mbr$.
%Recall that we write $\Psi$ for the standard Gaussian tail.
The collection $\bra{\mcw_\alpha}_{\alpha \in \mbr}$ of sets will satisfy
\(
	\pr{ W(\tpro_\alpha) \notin \mcw_\alpha } \eqsim \Psi(\alpha),
\)
using the CLT (\cref{res-p2:cutoff1:CLT}); see Proposition~\ref{res-p2:cutoff1:typ}.
Recall that $\Psi$ is the standard Gaussian tail.
We show that the modified $L_2$ distance (given typicality) is $\oh1$ at $\tpro_\alpha$; see \cref{res-p2:cutoff1:l2}.
Applying \cref{res-p2:cutoff1:mod-l2},
we find that
%We show that the expression in \cref{eq-p2:cutoff1:mod-l2:tv-l2} is $\oh1$.
%Then applying \cref{eq-p2:cutoff1:mod-l2:triangle} gives
\(
	d_{G_k}(\tpro_\alpha)
\le
	\Psi(\alpha) + \oh1
\)
\whp over $Z$.
This matches the lower bound from \S\ref{sec-p2:cutoff1:lower}.

By considering all $\alpha \in \mbr$, we find the shape of the cutoff. If we only desire the order of the window, then we need only consider the limit $\alpha \to \infty$; in this case, $\pr{ W(\tpro_\alpha) \notin \mcw_\alpha } \approx \Psi(\alpha) \approx 0$, which explains the use of the word `typical' in describing $\mcw_\alpha$.

The typicality conditions will be a combination of `local' (coordinate-wise) and `global' ones.

\begin{defn}%[Typicality]%: Local and Global]
\label{def-p2:cutoff1:typ}
	For all $\alpha \in \mbr$, define the \textit{local} and \textit{global typicality conditions}, respectively:
	\[
		\mcw_{\alpha,\loc}
	&\cq
		\brb{ w \in \mbz^k \midb \abs{ w_i - \ex{W_1(\tpro_\alpha)} } \le r_* \: \forall \, i = 1,...,k }
	\Qwhere
		r_* \cq \tfrac12 n^{1/k} \logk[2];
	\\
		\mcw_{\alpha,\glo}
	&\cq
		\brb{ w \in \mbz^k \midb \pr{ W(\tpro_\alpha) = w } \le n^{-1} e^{-\omega} }.
	\]
	Define $\mcw_\alpha \cq \mcw_{\alpha,\loc} \cap \mcw_{\alpha,\glo}$, and say that $w \in \mbz^k$ is ($\alpha$-)\textit{typical} if $w \in \mcw_\alpha$.
\end{defn}

The following propositions determine the probability that $W(\tpro_\alpha)$ lies in $\mcw_\alpha$, ie of typicality.

\begin{subtheorem}{thm}
	\label{res-p2:cutoff1:typ}

\begin{prop}
\label{res-p2:cutoff1:typ:loc}
	Let $X = (X_s)_{s\ge1}$ be a rate-$1$ RW---either a SRW or a DRW---on $\mbz$.
	Then,
	\[
		\pr[0]{ \abs{X_s - \ex{X_s}} > r }
	\le
		1/k^{3/2}
	\Qforall
		r \ge r_* = \tfrac12 n^{1/k} (\log k)^2
	\Quad{whenever}
		s \lesssim n^{2/k} \log k,
	\]
	where the subscript $0$ indicates that $X$ starts from $X_0 = 0$.
	In particular,
	for all $\alpha \in \mbr$,
	we have
	\[
		\pr{ W(\tpro_\alpha) \notin \mcw_{\alpha, \loc} }
	\le
		k^{-1/2}
	=
		\oh1.
	\]
\end{prop}

\begin{Proof}
The proof of the first part of this proposition follows from standard large deviation estimates on the RW on $\mbz$, and the fact that $\tpro_\alpha \asymp k n^{2/k}$ for all $\alpha \in \mbr$, as stated in \cref{res-p2:cutoff1:ent-times}.
The precise details are arduously technical. We defer them to \cite[\S\ref{sec-p0:rp}]{HOt:rcg:supp}.
The second part follows immediately from the first part and the union bound over the $k$ coordinates.
\end{Proof}

\begin{prop}%[Probability of Typicality]
\label{res-p2:cutoff1:typ:glo}
	For each $\alpha \in \mbr$,
	we have
	\[
		\pr{ W(\tpro_\alpha) \notin \mcw_{\alpha, \glo} }
	\to
		\Psi(\alpha).
	\]
\end{prop}

\end{subtheorem}

\begin{Proof}
This follows immediately from our CLT (\cref{res-p2:cutoff1:CLT}).
\end{Proof}

Herein, we fix $\alpha \in \mbr$ and frequently suppress the time $\tpro_\alpha$ from the notation, eg writing $W$ for $W(\tpro_\alpha)$ or $\mcw$ for $\mcw_\alpha$.
Let $V \cq W - W'$,
so
\(
	\bra{ W \bcdot Z = W' \bcdot Z }
=
	\bra{ V \bcdot Z = 0 }.
\)
Write
\[
	D
\cq
	D(\tpro_\alpha)
\cq
	n \, \pr{ V(\tpro_\alpha) \bcdot Z = 0 \mid \typ_\alpha } - 1
\Qwhere
	\typ 
\cq
	\typ_\alpha
\cq
	\brb{ W(\tpro_\alpha), W'(\tpro_\alpha) \in \mcw_\alpha) }.
\]
It remains to show that $D(\tpro_\alpha) = \oh1$ for all $\alpha \in \mbr$.
Recall \cref{hyp-p2:cutoff1}, the crux of which is that
\[
	\frac{k-d-1}{k} - 2 \frac{d \log\logk}{\log n} \ge 5 \frac{k}{\log n}
\Qand
	k - d \gg 1.
\]

For $r_1, ..., r_\ell \in \mbz \setminus \bra{0}$,
we use the convention
\(
	\gcd(r_1, ..., r_\ell,0) \cq \gcd(\abs{r_1}, ..., \abs{r_\ell}).
\)

\begin{prop}%[Modified $L_2$ Distance]
\label{res-p2:cutoff1:l2}
	Suppose that $(d, n, k)$ jointly satisfy \cref{hyp-p2:cutoff1}.
	(Recall that, implicitly, $(d, n, k)$ is a sequence of triples of integers.)
	Write
	\(
		\mfgcd
	\cq
		\gcd\rbr{V_1, ..., V_k, n}.
	\)
	Then,
	for all $\alpha \in \mbr$,
	we have
	\[
		0
	\le
		D(\tpro_\alpha)
	=
		\sumt{\mm \in \mbn}
		\prt{ \mfgcd = \mm \mid \typ } \cdot \abs{G}/\abs{\mm G}
	-	1
	=
		\oh1.
	\]
\end{prop}

Given this proposition, we can prove the upper bound in the main theorem, \cref{res-p2:cutoff1:res}.

\begin{Proof}[Proof of Upper Bound in \cref{res-p2:cutoff1:res} Given \cref{res-p2:cutoff1:l2}]
Fix $\alpha \in \mbr$ and consider the TV distance at time $\tpro_\alpha$.
%\cref{hyp-p2:cutoff1} is precisely the conditions required for \cref{res-p2:cutoff1:l2}.
Apply the modified $L_2$ calculation, ie \cref{res-p2:cutoff1:mod-l2,def-p2:cutoff1:typ}, at time $\tpro_\alpha$:
	by \cref{res-p2:cutoff1:l2}, the modified $L_2$ distance (given typicality) is $\oh1$ in expectation;
	by Markov's inequality, it is thus $\oh1$ \whp.
Proposition~\ref{res-p2:cutoff1:typ} says that typicality holds with probability $\Psi(\tpro_\alpha)$ asymptotically.
Combined, this all says that
\(
	d_{G_k}(\tpro_\alpha)
\le
	\Psi(\alpha) + \oh1
\)
\whp.
\end{Proof}

%\begin{Proof}[Proof of Upper Bound in \cref{res-p2:cutoff1:res} Given \cref{res-p2:cutoff1:l2}]
%	%
%\cref{hyp-p2:cutoff1} is precisely the conditions required for \cref{res-p2:cutoff1:l2}.
%Apply the modified $L_2$ calculation, \cref{res-p2:cutoff1:mod-l2,def-p2:cutoff1:typ}, and use \cref{res-p2:cutoff1:typ,res-p2:cutoff1:l2} to control the two resulting terms.
%Combined, these say that
%\(
%	d_{G_k}(\tpro_\alpha)
%\le
%	\Psi(\alpha) + \oh1
%\)
%\whp.
%	%
%\end{Proof}

It remains to prove \cref{res-p2:cutoff1:l2}, ie to bound the modified $L_2$ distance.
The remainder of the section is dedicated to this goal.
To do this, we are interested in the law of $V \bcdot Z$.

Obviously, when $V = 0$, we have $V \bcdot Z = 0$.
The following auxiliary lemma controls this.
%probability; its proof is deferred to the end of this subsection.
% \S\ref{sec-p2:cutoff1:aux-pf}.

\begin{lem}%[Empty]
\label{res-p2:cutoff1:V=0}
	We have
	\[
		n \, \pr{ V = 0 \mid \typ }
	\le
		e^{-\omega} / \prt{\typ}
	\lesssim
		e^{-\omega}
	=
		\oh1.
	\]
\end{lem}

\begin{Proof}%[Proof of \cref{res-p2:cutoff1:V=0}]
By direct calculation, since $W$ and $W'$ are independent copies,
\[
	\pr{ V = 0, \, \typ }
=
	\pr{ W = W', \, W \in \mcw }
=
	\sumt{w \in \mcw} \pr{W = w}^2.
\]
Recall global typicality: $\pr{W = w} \le n^{-1} e^{-\omega}$ for all $w \in \mcw$.
Thus
\[
	n \, \pr{ V = 0 \mid \typ }
\le
	n \sumt{w \in \mcw} \pr{W = w}^2 / \prt{\typ}
\le
	e^{-\omega} / \prt{\typ}.
\qedhere
\]
\end{Proof}

We now analyse $v \bcdot Z = \sum_i v_i Z_i$ for $v \ne 0$.
Sums of independent uniform random variables are uniform. Some simple technicalities take care of the fact that the $v_i$-s need not be $1$, or even equal.
%eg, if $X, Y \sim^\iid \Unif(\mbz_n)$, then $2X + 4Y \sim \Unif(2\mbz_n)$.

\begin{lem}[{\cite[\cref{res-p0:deferred:unif-gcd}]{HOt:rcg:supp}}]
\label{res-p2:cutoff1:unif-gcd}
	For all $v \in \mbz^k$,
	we have
	\[
		v \bcdot Z \sim \Unif\rbr{ \mm G }
	\Qwhere
		\mm \cq \gcd(v_1, ..., v_k, n).
%	\Qand
%		\abs{\mm G} \ge \mm^{-d(G)} \abs G.
	\]
\end{lem}

We now need to control $\abs{\mm G}$,
since \cref{res-p2:cutoff1:unif-gcd} implies~that
\[
	\pr{ V \bcdot Z = 0 \mid \typ }
=
	\sumt{\mm \in \mbn}
	\prt{ \mfgcd = \mm \mid \typ } / \abs{\mm G}
\Qwhere
	\mfgcd \cq \gcd\rbb{ V_1, ..., V_k, n }.
\]

\begin{lem}
\label{res-p2:cutoff1:G/gammaG}
	For all Abelian groups $G$ and all $\mm \in \mbn$,
	we have
	\[
%		\abs{G/\mm G}
		\abs G / \abs{\mm G}
	\le
		\mm^{d(G)}.
	\]
\end{lem}

\begin{Proof}
	Decompose $G$ as $\oplus_1^d \: \mbz_{m_j}$ with $d = d(G)$ and some $m_1, ..., m_d \in \mbn$.
	Then $\mm G$ can be decomposed as $\oplus_1^d \: \gcd(\mm, m_j) \mbz_{m_j}$.
	Hence,
	\(
		\abs{\mm G}
	=
		\prodt[d]{1} \rbr{ m_j/\gcd(\mm, m_j) }
	\ge
		\prodt[d]{1} \rbr{ m_j/\mm }
	=
		\abs G / \mm^d.
%	\qedhere
	\)
\end{Proof}

These lemmas combine to produce a simple, but key, corollary.
%Recall that
%\(
%	\mfgcd
%=
%	\gcd\rbr{V_1, ..., V_k, n}.
%\)

\begin{cor}
\label{res-p2:cutoff1:pr-gcd}
	We have
	\[
		n \, \pr{ V \bcdot Z = 0, \: V \ne 0 \mid \typ }
	\le
		\ex{ \mfgcd^{d(G)} \, \one{V \ne 0} \mid \typ }.
	\]
\end{cor}

\begin{Proof}
The conditioning does not affect $Z$,
so the claim is immediate from the previous lemmas.
\end{Proof}

We control this gcd coordinate-by-coordinate,
using a crude divisibility bound.
%The proof is deferred to the end of this subsection.
%Remember that $V = V(\tpro_\alpha)$.

\begin{lem}
\label{res-p2:cutoff1:divis}
	Recall that $V = V(\tpro_\alpha)$.
	For all $\mm \in \mbn$, we have
	\[
		\pr{ V_1 \in \mm \mbz \mid V_1 \ne 0 }
	\le
		1/\mm
	\Qand
		\prt{ \mfgcd = \mm \mid \typ }
	\lesssim
		\rbb{ 1/\mm + 2/n^{1/k} }^k.
	\]
\end{lem}

\begin{Proof}%[Proof of \cref{res-p2:cutoff1:divis}]
For both the SRW and DRW, the difference $V = W - W'$ is a rate-$2$ SRW on $\mbz^k$.
Hence, each coordinate is an independent rate-$2/k$ SRW on $\mbz$, which is symmetric about $0$.

It is easy to see that any non-increasing distribution on $\mbn$ can be written as a mixture of $\Unif(\bra{1, ..., Y})$ distributions, for different $Y \in \mbn$.
The map $m \mapsto \pr{\abs{V_1(2t/k)} = m} : \mbn \to [0,1]$ is non-increasing for any $t \ge 0$.
Hence, $\abs{V_1}$ conditional on $V_1 \ne 0$ has such a distribution.~%
Thus,
\[
	\abs{V_1} = \abs{V_1(\tpro_\alpha)} \sim \Unif\bra{1,...,Y}
\Quad{conditional on}
	V_1 \ne 0,
\]
where $Y$ has some distribution.
Hence, we have
\[
	\pr{ V_1 \in \mm \mbz \mid V_1 \ne 0 }
=
	\ex{ \floorb{Y/\mm} \bigl/ Y }
\le
	1/\mm.
\]

If the gcd $\mfgcd = \mm$, then $V_i \in \mm \mbz$ for all $i \in [k]$.
By independence of coordinates, we then obtain
\[
	\prt{ \mfgcd = \mm \mid \typ }
%&
\le
	\prt{ \mfgcd = \mm } / \prt{\typ}
%\\&
\lesssim
	\prt{ V_1 \in \mm \mbz }^k
%\\&
\le
	\rbb{ \prt{ V_1 = 0 } + \prt{ V_1 \in \mm \mbz \mid V_1 \ne 0 } }^k,
\]
noting that $\prt{\typ} \asymp 1$.
Using \cref{res-p2:cutoff1:ent-times} to argue that
\(
	\pr{ V_1 = 0 }
\le
	2/n^{1/k},
\)
%and the previous part of this lemma,
we deduce that
\[
	\prt{ \mfgcd = \mm \mid \typ }
\lesssim
	\rbb{ 2/n^{1/k} + 1/\mm }^k.
\qedhere
\]
\end{Proof}

From this, using the \cref{hyp-p2:cutoff1}, we can deduce that
\(
	\ex{ \mfgcd^{d(G)} \, \one{V \ne 0} \mid \typ } = 1 + \oh1.
\)

\begin{cor}
\label{res-p2:cutoff1:gcd-ex}
	Recall that time $t = \tpro_\alpha$.
	Given \cref{hyp-p2:cutoff1},
	we have
	\[
		\exb{ \mfgcd^{d(G)} \, \one{V \ne 0} \mid \typ } = 1 + \oh1.
	\]
\end{cor}

This proof is briefly deferred.
First, we deduce \cref{res-p2:cutoff1:l2} from the above results.

\begin{Proof}[Proof of \cref{res-p2:cutoff1:l2}]
By \cref{res-p2:cutoff1:V=0,res-p2:cutoff1:pr-gcd,res-p2:cutoff1:gcd-ex},
we have
\[
	n \, \pr{ V \bcdot Z = 0 \mid \typ }
&
\le
	n \, \pr{ V = 0 \mid \typ }
+	n \, \pr{ V \bcdot Z = 0, \, V \ne 0 \mid \typ }
\\&
\le
	n \, \pr{ V = 0 \mid \typ }
+	\ex{ \mfgcd^{d(G)} \, \one{V \ne 0} \mid \typ }
=
	1 + \oh1.
\qedhere
\]
\end{Proof}

We close the analysis of Approach \#1 with the briefly-deferred proof of \cref{res-p2:cutoff1:gcd-ex}.

\begin{Proof}[Proof of \cref{res-p2:cutoff1:gcd-ex}]
Let $d \cq d(G)$.
By local typicality,
\(
	\mfgcd \le 2 r_* = n^{1/k} \logk[2]
\)
if $V \ne 0$.
Thus,
%For $\mm = 1$, we upper bound $\pr{g = \mm \mid \typ} \le 1$.
\[
	\ex{ \mfgcd^d \, \one{V \ne 0} \mid \typ }
=
	\sumt{\mm \in \mbn}
	\mm^d \pr{ \mfgcd = \mm \mid \typ }
\le
	1
+	\sumt[\floor{n^{1/k} \logk[2]}]{\mm = 2}
	\mm^d \, \prt{ \mfgcd = \mm \mid \typ }.
\]
For $\mm \ge 2$, we use \cref{res-p2:cutoff1:divis}.
Let $\delta \in (0,1)$. % with $\delta \ge 1/k$.
For $2 \le \mm \le \delta n^{1/k}$, we use the bound
\[
	\prt{ \mfgcd = \mm \mid \typ }
\lesssim
	\rbb{ 1/\mm + 2/(\mm/\delta) }^k
=
	(1 + 2\delta)^k / \mm^k.
\]
For $\mm \ge \delta n^{1/k}$, we use the slightly crude bound
\(
	(a + b)^k
\le
	2^k (a^k + b^k)
\)
for $a,b \ge 0$ to deduce that
\[
	\prt{ \mfgcd = \mm \mid \typ }
\lesssim
	2^k \rbb{ 1/\mm^k + 2^k / n }
=
	2^k/\mm^k + 4^k / n.
\]

Dividing the appropriate sum over $\mm$ into two parts according to whether or not $\mm \le \delta n^{1/k}$ and using the above inequalities, elementary algebraic manipulations can be used to deduce that
\[
	\ex{ \mfgcd^d \, \one{V \ne 0} \mid \typ } - 1
\lesssim
	e^{2\delta k} 2^{d+1-k}
+	2^k \delta^{d+1-k} n^{(d+1-k)/k}
+	4^k n^{(d+1)/k} \logk[2(d+1)]/n.
%\label{eq-p2:cutoff1:l2:gcd-decomp}
%\nt
\]
This is $\oh1$, by the conditions of \cref{hyp-p2:cutoff1},
%given in \cref{eq-p2:cutoff1:hyp},
as we now explain.
Write $\eta \cq (k-d-1)/k \in (0,1)$.
\begin{itemize}[noitemsep, label = \bcdot]
	\item 
	We wish to choose $\delta$ as large as possible, but with the first term $\oh1$:
	set $\delta \cq \tfrac14 \eta$.
	
	\item 
	\cref{hyp-p2:cutoff1} implies that $\eta \ge 4 k / \log n$,
	which shows that the second term is $\oh1$.
	
	\item 
	The inequality in \cref{hyp-p2:cutoff1} is designed precisely so that the final term is $\oh1$.
%	since $\eta k \ge \tfrac12(k - d)$.
\qedhere
\end{itemize}
\end{Proof}

\begin{rmkt*}
%\label{rmk-p2:cutoff1:relax-k-d}
	%
We have always assumed that $k - d(G) \gg 1$.
Our analysis does apply if $M \cq k - d(G) \ge 2$ is fixed (ie not diverging) too.
Then, however, it is not necessarily the case that the group is generated \whp---eg if $G = \mbz_2^d$ then it is not.
Our analysis shows that the mixing time is of order $\tpro_0$ with probability bounded away from $0$, and approaching $1$ as $M$ grows.
\end{rmkt*}

\section{TV Cutoff: Approach \#2} \label{sec-p2:cutoff2}
%\section{Abelian Groups: Cutoff; Approach \#2} \label{sec-p2:cutoff2}

Recall that \cref{res-p2:intro:tv} is established via two distinct approaches.
In the previous section, we used one approach to deal with the case that $k$ is `not too large'.
Here, we use a new approach to deal with the case that $k$ is `not too small'.
The main result of the section is \cref{res-p2:cutoff2:res}.
%see also \cref{hyp-p2:cutoff2,rmk-p2:cutoff2:hyp}.
The outline of the section is roughly the same as that of the previous one.

\begin{itemize}[noitemsep, topsep = 0pt, label = \bcdot]
	\item 
	\S\ref{sec-p2:cutoff2:ent:method} discusses the new, refined entropic methodology.
	
	\item 
	\S\ref{sec-p2:cutoff2:ent:def} defines the new \textit{entropic times}.
	
	\item 
	\S\ref{sec-p2:cutoff2:ent:growth-conc} states bounds on the growth rate of the entropy and concentration.
	
	\item 
	\S\ref{sec-p2:cutoff2:res} states precisely the main theorem of the section.
	
	\item 
	\S\ref{sec-p2:cutoff2:outline} outlines the differences between this argument and the previous approach.
	
	\item 
	\S\ref{sec-p2:cutoff2:lower} is devoted to the lower bound.
	
	\item 
	\S\ref{sec-p2:cutoff2:upper} is devoted to the upper bound.
\end{itemize}
%	\item 
%	\S\ref{sec-p2:cutoff2:special} discusses some special cases of particular interest.

\subsection{Entropic Times: New Methodology and Definition}
\label{sec-p2:cutoff2:ent:method}

The underlying principles of the method used in this section (\S\ref{sec-p2:cutoff2}) are the same as those of the previous one (\S\ref{sec-p2:cutoff1}), just adjusted slightly to deal with the cases not previously covered.
%We adjust the method slightly to deal with the cases which were not covered in \S\ref{sec-p2:cutoff1}.% \cref{res-p2:cutoff1:res}.

We first discuss where the previous approach broke down and how we might fix it.
The primary issue was when $d(G)$ was very large.
Eg, consider $\mbz_2^d$. All elements are of order $2$, so instead of looking at $W$, a RW on $\mbz$, we could equally have taken $W$ mod $2$. The entropy of $W_1(t)$ mod $2$ is significantly smaller than that of $W_1(t)$ once $t/k \gtrsim 1$.
This suggests a longer mixing time.

Now,
\(
	V \bcdot Z \sim \Unif(\mm G)
\)
when $\gcd(V_1, ..., V_k, n) = \mm$.
%(This assumes that the group $G$ is Abelian.)
This motivates defining $\tent_\mm$ to be the time at which the entropy of $W_1$ mod $\mm$ is $\log \abs{G/\mm G}$, and proposing $\tent_* \cq \max_{\mm \in \mbn} \tent_\mm$ as the upper~bound.

\subsection{Entropic Times: Definition and Concentration}
\label{sec-p2:cutoff2:ent:def}

In this section, we refine the definition of \textit{entropic times}.
The concept is highly analogous to that of the previous section.
There is, thus, some overlap in both verbal description and notation.
We have been careful, though, to set it up as not to cause confusion:
	we always use indices such as $\mm \in \mbn$ in the new entropic times $\tent_0(\mm, N)$ or $\tent_\mm$ below,
	whilst previously we used $\alpha \in \mbr$ in $\tpro_\alpha$.

\smallskip

We now define precisely the (updated) notion of \textit{entropic times}.
Let $W = \rbr{W_i(t) \mid i \in [k], \ t \ge 0}$ be a RW on $\mbz^k$, counting the uses of the generators, as in the previous sections.
This can be either a SRW on $\mbz^k$ or DRW on $\mbz_+^k$.
As before, $S(t) = W(t) \bcdot Z$.
%(We sometimes write ``RW'' in place of ``SRW/DRW''.)
For $\mm \in \mbn \cup \bra{\infty}$, define $W_\mm$ via
\[
	W_{\mm,i}(t) \cq W_i(t) \mod \mm
\Qfor
	\mm \in \mbn
\Qand
	W_\infty \cq W.
\]
Then, $W_\mm$ is a RW on $\mbz_\mm^k$.
So,
\(
	W_{\mm,i} \cq \rbr{ W_{\mm,i}(t) }_{t\ge0}
\)
forms an iid sequence
%(over $i \in [k]$)
of rate-$1/k$ RWs on $\mbz_\mm$.
As before, ``mod $\infty$'' has no effect:
	$w = w$ mod $\infty$ for all $w \in \mbz^k$,
	and $G / \infty G = G$ as $\infty G = \bra{\id}$.

Write $\mu_{\mm,t}$, respectively $\nu_{\mm,s}$, for the law of $W_\mm(t)$, respectively $W_{\mm,1}(sk)$;
so,
\(
	\mu_{\mm,t} = \nu_{\mm,t/k}^{\otimes k}.
\)~%
Define
\[
	Q_\mm(t) \cq - \log \mu_{\mm,t}\rbb{W_\mm(t)}
\Qand
	Q_{\mm,i}(t) \cq - \log \nu_{\mm,t/k}\rbb{W_{\mm,i}(t)}.
\]
So, $Q_{\mm,i}$ forms an iid sequence over $i \in [k]$;
also,
\[
	Q_\mm(t) = \sumt[k]{i=1} Q_{\mm,i}(t),
\quad
	h_\mm(t) \cq \ex{ Q_\mm(t) }
\Qand
	H_\mm(s) \cq \ex{ Q_{\mm,1}(sk) }.
\]
So, $h_\mm(t)$ and $H_\mm(s)$ are the entropies of $W_\gamma(t)$ and $W_{\gamma,1}(sk)$, respectively.
Note that $h_\mm(t) = k H_\mm(t/k)$ and that $h_\mm : [0,\infty) \to [0, \log(\mm^k))$ is a strictly increasing bijection.

Some of these expressions, such as $h_\mm$, depend on $k$; we usually suppress this from the notation.

\begin{defn}%[Entropic Times: Original]
\label{def-p2:cutoff2:ent-orig}
For $N < \mm^k$, define the \textit{entropic time}
\[
	\tent_0(\mm,N)
\cq
	h_\mm^{-1}( \log N )
\Qand
	\sent_0(\mm,N)
\cq
	\tent_0(\mm,N)/k
=
	H_\mm^{-1}( \log N / k ).
\]
We are interested primarily in $N \cq \abs{G/\mm G}$ for an Abelian group $G$:
set
\[
	\tent_*
\cq
	\tent_*(k, G)
\cq
	\maxt{\mm \in \mbn \cup \bra{\infty}}
	\tent_0(\mm, \abs{G/\mm G}).
\]
%where the notation $a \wr b$ ($a,b \in \mbn$) means that $a$ divides $b$, ie $b \in a \mbz$.
	%
\end{defn}

This $\tent_*$ is the same as defined in the introduction; see \cref{def-p2:intro:ent-time}.
Comparing with notation in Approach \#1, $\tent_0(\infty, \abs G) = \tpro_0$; see \cref{def-p2:cutoff1:ent-times}.
We establish \emph{cutoff} here, not the \emph{profile} as~well.

Recall that $\infty G = \abs G G = \bra{\id}$ and $1 G = G$.
So, the maximum is achieved at some $\gamma \in [2, n]$.
Below, for brevity, we write ``$\gamma \ge 2$'' to mean ``$\gamma \in \mbn \cup \bra{\infty} \setminus \bra{1}$''.

Our next result determines the asymptotics of $\tent_*$.

%\begin{subtheorem}{thm}
%	\label{res-p2:cutoff2:ent:ent-times}
%
%\begin{prop}[{\cite[\cref{res-p0:re:app:s*:order:asymp}]{HOt:rcg:supp}}]
%\label{res-p2:cutoff2:ent:ent-times:asymp}
%	If $k - d(G) \asymp k$, then
%	\(
%		\tent_*(k,G) \asymp k \abs G^{2/k}.
%	\)
%\end{prop}
%
%\begin{prop}[{\cite[\cref{res-p0:re:app:s*:order:>>}]{HOt:rcg:supp}}]
%\label{res-p2:cutoff2:ent:ent-times:>>}
%	If $k - d(G) > 1$, then
%	\(
%		\tent_*(k,G) \lesssim k \abs G^{2/k} \log k.
%	\)
%\end{prop}
%
%\end{subtheorem}

\begin{subtheorem}{thm}
	\label{res-p2:cutoff2:ent:eval}

\begin{prop}%[{\cite[Proposition~\ref{res-p0:re:app:s*:order}]{HOt:rcg:supp}}]
\label{res-p2:cutoff2:ent:eval:order}
	Suppose that $1 \ll k \lesssim \log \abs G$.
	The following hold:
	\begin{alignat*}{3}
		\text{if}
	\quad
		k - d(G) &\asymp k,&
	\Quad{then}&&
		\tent_*(k,G) &\asymp k \abs G^{2/k};
	\\
		\text{if}
	\quad
		k - d(G) &\ge 1,&
	\Quad{then}&&
		k \abs G^{2/k} \lesssim \tent_*(k,G) &\lesssim k \abs G^{2/k} \log k.
	\end{alignat*}
\end{prop}

\begin{prop}%[{\cite[\cref{res-p0:re:app:s*:eqsim:d<<log}]{HOt:rcg:supp}}]
\label{res-p2:cutoff2:ent:eval:eqsim}
	Suppose that $d(G) \ll \log \abs G$ and $k - d(G) \asymp k \gg 1$.
	Then,
	\[
		\tent_*(k,G)
	\eqsim
		\tent_0(\infty, G)
	=
		\tpro_0.
	\]
%	(Note that $\tent_0(\infty, \abs G) = \tent_0(k, \abs G)$ in the notation of \cref{def-p2:cutoff1:ent-times}.)
\end{prop}

\end{subtheorem}

As with earlier results, the rigorous proofs are technical.
%We give heuristics behind the ideas of proof below.
They boil down to comparing the RW on $\mbz_\mm$ with one on $\mbz$.
Precise details are deferred to \cite[Propositions~\ref{res-p0:re:app:s*:order} and~\ref{res-p0:re:app:s*:eqsim:d<<log}]{HOt:rcg:supp}.%

In \S\ref{sec-p2:cutoff2:lower}, we show that $\tent_0(\mm, \abs{G/\mm G})$ is a lower bound on mixing for all $\mm$, \forallZ.
Throughout this section, we work under the assumption that $k \lesssim \log \abs G$.
(Recall from \S\ref{sec-p2:intro:previous-work} that cutoff had already been established for all Abelian groups when $k \gg \log \abs G$.)
As a result of this, taking $\mm \cq \abs G$, we see that the mixing time is at least order $k$. Indeed, $|G|G = \bra{\id}$, so the target entropy per coordinate is $\log \abs G / k \gtrsim 1$, so each coordinate needs to be run for time $t/k \gtrsim 1$.
Hence, there exists a $\varsigma > 0$ so that the mixing time is at least $2 \varsigma k$.
This holds for all $Z$, not just \whp over $Z$.

Unfortunately, the $\mm$ with $\tent_0(\mm, \abs{G/\mm G}) \le \varsigma k$ cause some technical difficulties. For this reason, we take the maximum with $\varsigma k$ in the definition of the `adjusted' entropic time $\tent_\mm$ below.
Crucially,
\[
%	\maxt{\mm \in \mbn} \tent_\mm
%=
	\maxt{\mm \in \mbn}
	\rbb{ \tent_0(\mm, \abs{G/\mm G}) \vee \varsigma k }
=
	\maxt{\mm \in \mbn}
	\tent_0(\mm, \abs{G/\mm G})
\vee
	\varsigma k
=
	\tent_*.
\]
We emphasise that this last adjustment is purely technical.
On the other hand, the entropic times in \cref{def-p2:cutoff2:ent-orig} capture properties of the group to which those in \cref{def-p2:cutoff1:ent-times} are oblivious.

\begin{defn}%[Entropic Times: Adjusted]
\label{def-p2:cutoff2:ent-adj}
	For $\mm \ge 2$ and $s \ge 0$,
	define the \textit{(adjusted) entropic time} and \textit{relative entropy} via
	\[
		\sent_\mm \cq \sent_0(\mm, \abs{G/\mm G}) \vee \varsigma,
	\quad
		\tent_\mm \cq \sent_\mm k
	\Qand
		R_\mm(s) \cq \log \mm - H_\mm(s).
	\]
\end{defn}

The maximal entropy of a distribution on $\mbz_\mm$ is $\log \mm$, obtained uniquely by the uniform distribution $\Unif(\mbz_\mm)$.
Hence, $R_\mm(s) \to 0$ \asinf s, since the RW converges to $\Unif(\mbz_\mm)$.

\subsection{Entropic Times: Entropy Growth Rate and Concentration}
\label{sec-p2:cutoff2:ent:growth-conc}

We determine the rate of change of the entropy around the entropic time and establish concentration estimates on the `random entropy' $Q_\mm$ at a time shortly after the entropic time.
%This gives cutoff, but even more refined estimates---namely, a CLT---are required to obtain the profile.

\medskip

The first lemma controls the rate of change of the entropy near the entropic time.
%; see \cite[\S\ref{sec-p0:re}]{HOt:rcg:supp}.
%(The claim is not true if one defines $\tent_\mm$ \emph{without} the maximum with $\varsigma$.)

\begin{lem}%[{\cite[\cref{res-p0:re:app:ent-growth-rate}]{HOt:rcg:supp}}]
%	[Growth Rate of Entropy; {\cite[\cref{res-p0:re:app:ent-growth-rate}]{HOt:rcg:supp}}]
\label{res-p2:cutoff2:ent:entropy-growth}
	There exists a continuous function $\widebar c : (0,1) \to (0,1)$ so that,
	for
		all $\mm \ge 2$,
		all $\xi \in (-1,1) \setminus \bra{0}$
	and
		all $s \ge \varsigma$,
	we have
	\[
		\absb{ H_\mm\rbb{ s (1 + \xi) } - H_\mm\rbr{ s } }
	\ge
		2 \widebar c_{\abs \xi} \rbb{ R_\mm(s) \wedge 1 }.
	\]
\end{lem}

\begin{Proof}[Outline of Proof]
If $s \asymp 1$, then it is easy as all terms are order $1$.
If $s \ll \mm^2$, then the fact that the RW is on $\mbz_\mm$, not $\mbz$, is not significant. The entropy is thus approximately $\tfrac12 \log(2 \pi e s)$, ie that of $N(0, s)$, if also $s \gg 1$.
The case $s \gtrsim \mm^2$ follows from standard (modified) log-Sobolev arguments.

Making this proof rigorous is technical. Doing so is deferred to \cite[\cref{res-p0:re:app:ent-growth-rate}]{HOt:rcg:supp}.
\end{Proof}

Recall that $\tent_* = \max_{\mm \in \mbn} \tent_\mm$.
Abbreviate $d = d(G)$.
For $\mm \in \mbn$, write
\[
	\zeta_\mm
\cq
	\tfrac1k \rbb{ k - d(G) } \log \mm.
\]

\begin{prop}%[{\cite[\cref{res-p0:re:app:conc}]{HOt:rcg:supp}}]
%	[Concentration of $Q_\mm$; {\cite[\cref{res-p0:re:app:conc}]{HOt:rcg:supp}}]
\label{res-p2:cutoff2:ent:conc}
	Assume that $k > d$.
	There exists a continuous function $c : (0,1) \to (0,1)$ so that,
	for all $\mm \ge 2$ and all $\eps \in (0,1)$,
	the following hold:
	\[
		\pr{
			Q_\mm\rbb{ \tent_* (1 + \eps) }
		\le
			\log \abs{G/\mm G} + c_\eps \rbr{ \zeta_\mm \wedge 1 } k
		}
	&\le
		\expb{ - c_\eps \rbr{ \zeta_\mm \wedge 1 } k };
	\\
		\pr{
			Q_\mm\rbb{t(1 - \eps) }
		\ge
			\log \abs{G/\mm G} - c_\eps \rbr{ \zeta_\mm \wedge 1 } k
		}
	&=
		\oh1
	\Quad{for all}
		t \le \tent_\mm.
	\]
%	(The first inequality considers time $\tent_*$, while the second considers $t \le \tent_\mm$.)
\end{prop}

%\color{orange}
%\color{.!80!black}
%The technicalities of the proof are deferred to \cite[\cref{res-p0:re:app:conc}]{HOt:rcg:supp}.
%Heuristics are given below.%
%\color{black}

\begin{Proof}[Outline of Proof]
	%
%The proof of this proposition is given in \cite[\S\ref{sec-p0:re}]{HOt:rcg:supp}.
%We give a brief outline here.
%Consider $k - d \asymp k$, so that $\zeta_\mm \wedge 1 \asymp 1$.
Recall that $Q_\mm(t) = \sumt[k]{1} Q_{\mm,i}(t)$ is a sum of iid terms, each of mean $H_\mm(t/k)/k$.
By the entropy growth rate (\cref{res-p2:cutoff2:ent:entropy-growth}), for any $\xi \in (-1,1) \setminus \bra{0}$, the change in entropy between times $s$ and $(1 + \xi) s$ is order $R_\mm(s) \wedge 1$, with implicit constant depending on $\abs \xi$.
Taking $s \cq \sent_0(\mm, \abs{G/\mm G})$, recalling that $\abs{G/\mm G} \le \mm^{d(G)}$ by \cref{res-p2:cutoff1:G/gammaG},
gives
\[
	R_\mm(s)
=
	\log \mm - H_\mm(s)
=
	\log \mm - (\log \abs{G/\mm G})/k
\ge
	\tfrac1k \rbb{ k - d(G) } \log \mm
=
	\zeta_\mm.
\]
We are interested in the times $\sent_\mm$, not $\sent_0(\mm, \abs{G/\mm G})$; this is only a minor technical complication.

The quantitative concentration estimate requires first deterministically bounding $\ex{Q_{1,\mm}} - Q_{1,\mm}$ from above. A (one-sided) variant of Bernstein's inequality for a sum of iid, deterministically-bounded random variables is then applied.
The non-quantitative part is just an application of Chebyshev, once the variance $\Var{Q_{\mm,1}(sk)}$ has been uniformly bounded over $s \ge \varsigma$.

Again, making argument this rigorous is technical. It is deferred to \cite[\cref{res-p0:re:app:conc}]{HOt:rcg:supp}.
\end{Proof}

\subsection{Precise Statement and Remarks}
\label{sec-p2:cutoff2:res}

In this subsection, we state precisely the main theorem of the section.
There are some simple conditions on $k$, in terms of $d(G)$ and $\abs G$, needed for the upper bound.

\begin{hyp}
\label{hyp-p2:cutoff2}
	The sequence $(k_N, G_N)_\Ninn$ satisfies \textit{\cref{hyp-p2:cutoff2}} if
%	each individually is a non-decreasing sequence of positive integers and
	the following hold:
	\begin{gather*}
		\LIMSUP{\Ninf} k_N / \log \abs{G_N} < \infty,
	\quad
		\LIMINF{\Ninf} \bra{ k_N - d(G_N) } = \infty
	\Qand
		\LIMINF{\Ninf} k_N / \log(\abs{\mch_N} + 1) = \infty,
	\\
		\text{where}
	\quad
		\mch_N
	\cq
%		\brb{ \mm G_N \midb \mm \in [2, n_{*,N}] \text{ and } \abs{G_N} \in \mm \mbz }
		\brb{ \mm G_N \midb \mm \wr \abs{G_N} \text{ and } \mm \in [2, n_{*,N}] }
	\Qand
		n_{*,N} \cq \floor{ \abs{G_N}^{1/k_N} (\log k_N)^2 },
	\end{gather*}
	where the notation $a \wr b$
%	($a,b \in \mbn$)
	means that $a$ divides $b$, ie $b \in a \mbz$,
	for $a,b \in \mbn$.
\end{hyp}

\begin{rmkt}
\label{rmk-p2:cutoff2:hyp}
	If $k \gg \sqrt{\log n}$, then $k \gg \log(\abs \mch + 1)$, since $\abs \mch \le n_* \le n^{1/k} \logk[2]$.
%	We show in \S\ref{sec-p2:cutoff2:special} that if $\abs G$ lies in a certain density-$1$ subset of $\mbn$, then $k \gg \sqrt{\log\log \abs G}$ implies that $k \gg \log \abs \mch$.
\end{rmkt}

%In \cref{rmk-p2:cutoff2:hyp} below, we give a sufficient condition for \cref{hyp-p2:cutoff2} to hold.

Throughout the proofs, we suppress the subscript-$N$, eg writing $k$ or $n$, considering sequences implicitly.
Recall that we abbreviate the TV distance from uniformity at time~$t$~as
\[
	d_{G_k, N}(t) = \tvb{ \pr[{G_N([Z_1, ..., Z_{k_N}])}]{ S(t) \in \cdot } - \pi_{G_N} }
\Qwhere
	Z_1, ..., Z_{k_N} \sim^\iid \Unif(G_N).
\]

We now state the main theorem of this section.
Recall that
\(
	\tent_*
=
	\max_\mm \tent_0(\mm,\abs{G/\mm G})
=
	\max_\mm \tent_\mm.
\)

\begin{thm}%[Cutoff: Upper Bound \#2]
\label{res-p2:cutoff2:res}
	Let $(k_N)_\Ninn$ be a sequence of positive integers and $(G_N)_\Ninn$ a sequence of finite, Abelian groups;
	for each $\Ninn$, define $Z_{(N)} \cq [Z_1, ..., Z_{k_N}]$ by drawing $Z_1, ..., Z_{k_N} \sim^\iid \Unif(G_N)$.
	
	Suppose that the sequence $(k_N, G_N)_\Ninn$ satisfies \cref{hyp-p2:cutoff2}.
	Let $c \in (-1,1) \setminus \bra{0}$.
	Then,
	\[
		d_{G_k, N}\rbb{ (1 + c) \tent_*(k_N, G_N) }
	\to^\mbp
		\one{c < 0}
	\Quad{(in probability)}
		\asinf N.
	\]
	That is, \whp, there is TV cutoff at $\tent_* = \max_\mm \tent_0(\mm, \abs{G/\mm G})$.
	Moreover, the implicit lower bound on the TV distance holds deterministically, ie for all choices of generators.
\end{thm}

\subsection{Outline of Proof}
\label{sec-p2:cutoff2:outline}

The general outline is analogous to that before; see \S\ref{sec-p2:cutoff1:outline}.
%The general outline of this approach is the same as that of the previous; see \S\ref{sec-p2:cutoff1:outline} for an outline of the previous approach.
That approach failed once either $d$ or $k$ became too large or $k - d$ became too small.
We outline here the ideas used to cover these cases.

\smallskip

For the lower bound, we project the walk from $G$ to $G/\mm G$.
This cannot increase the TV distance.
The same argument shows that $\tent_0(\mm, \abs{G/\mm G})$ is a lower bound for all $\mm$.
%The idea, then, is that where before we looked at a RW on $\mbz^k$ and waited until it has entropy $\log \abs G$, instead we look at a RW on $\mbz_\mm^k$ and wait until it has entropy $\log \abs{G/\mm G}$; see \cref{def-p2:cutoff2:ent-orig}.
%We then take a worst-case over $\mm \in \mbn$.

For the upper bound,
fundamentally, we still wish to bound the same expression:
\[
	D\rbr{t}
=
	\sumt{\mm \in \mbn}
	\prt{ \mfgcd = \mm \mid \typ } \cdot \abs{G}/\abs{\mm G}
-	1;
\]
see \cref{res-p2:cutoff1:l2,res-p2:cutoff2:l2}.
If $\mfgcd = \mm$ then $W \equiv W' \mod \mm$.
But $W_\mm \cq W \mod \mm$ and $W'_\mm \cq W' \mod \mm$ are simply RWs on $\mbz_\mm^k$.
So, the same argument as in \cref{res-p2:cutoff1:V=0} gives
\[
	\pr{ \mfgcd = \mm \mid \typ }
\le
	\pr{ W_\mm = W'_\mm \mid \typ }
\ll
	1/\abs{G / \mm G}
=
	\abs{\mm G}/\abs G
\Qwhen
	t \ge \tent_0(\mm, \abs{G/\mm G}).
\]
Thus,
\(
	D(t)
=
	1 + \oh1
\)
when $t \ge \tent_0(\mm, \abs{G/\mm G})$ for all $\mm \ge 2$.
Hence, the proposed upper bound of
\[
	\tent_*
=
	\maxt{\mm \ge 2} \tent_0(\mm, \abs{G/\mm G}).
\]
The adjusted entropic times $\tent_\mm$ are only introduced to alleviate a technical problem.

\renewcommand*{\mm}{\gamma}

\subsection{Lower Bound on Total-Variation Mixing}
\label{sec-p2:cutoff2:lower}

The idea is to quotient out by $\mm G$, and show that the walk on this quotient is not mixed at time $(1 - \eps) \tent_0(\mm, \abs{G/\mm G})$, as in  \S\ref{sec-p2:cutoff1:lower}. Hence, the original walk is not mixed on $G$ either.
This is achieved via the entropy growth rate and variance bounds detailed in \cref{res-p2:cutoff2:ent:conc}.

\begin{Proof}[Proof of Lower Bound in \cref{res-p2:cutoff2:res}]
As in \S\ref{sec-p2:cutoff1:lower}, for this proof only, to emphasise that $Z$ is fixed, not being averaged over, we add a subscript-$Z$ to the probabilities involving $Z$: $\pr[Z]{ S(t) \in \cdot }$.
%Again, those involving only $W$ are unaffected, of course, by $Z$.

%For this proof, assume that $Z$ is given, and suppress it.

Fix $\eps \in (0,1)$ and set $t \cq (1 - \eps) \tent_0(\mm, \abs{G/\mm G})$.
Recall that
\(
	\zeta_\mm
=
	\tfrac1k \rbr{ k - d(G) } \log \mm.
\)
%Write
%\(
%	\zeta_\mm
%\cq
%	R_\mm\rbr{ \sent_0(\mm, \abs{G/\mm G}) }.
%\)
Then,
\[
	\prt{ \mce }
=
	1 - \oh1
\Qwhere
	\mce
\cq
	\brb{ \mu_{\mm,t}\rbb{ W_\mm\rbr{ t } } \ge \delta_{\mm}^{-1} / \abs{G/\mm G} }
\Qand
	\delta_{\mm}
\cq
	\expb{ - c_\eps (\zeta_\mm \wedge 1) k },
\]
by the entropy concentration (\cref{res-p2:cutoff2:ent:conc}).
Also, $\abs{G/\mm G} \le \mm^{d(G)}$ by \cref{res-p2:cutoff1:G/gammaG},
and so
\[
	R_\mm\rbb{ \sent_0(\mm, \abs{G/\mm G}) }
=
	\log \mm - \log \abs{G/\mm G}/k
\ge
	\tfrac1k \rbb{ k - d(G) } \log \mm
=
	\zeta_\mm.
\]
Also, $k - d(G) \gg 1$, by assumption.
Thus, $\delta_{\mm} = \oh1$ uniformly in $\mm$.
Consider the set
\[
	E
\cq
	\brb{ x \in G/\mm G \midb \exists \, w \in \mbz_\mm^k \st \mu_{\mm,t}(w) \ge \delta_{\mm}^{-1} / \abs{G/\mm G} \text{ and } x = (w \bcdot Z) \mm G }.
\]
Define $S_\mm$ to be the projection of $S = W_\infty \bcdot Z$ to $G/\mm G$.
Then, $\pr[Z]{ S_\mm(t) \in E \mid \mce } = 1$ since $W_\infty$ generates $S$.
Every element $x \in E$ can be realised as $x = w_x \bcdot Z$ for some $w_x \in \mbz_\mm^k$ with $\mu_{\mm,t}(w_x) \ge \delta_{\mm}^{-1} / \abs{G/\mm G}$, by definition of $E$.
Hence, for all $x \in E$, we have
\[
	\pr[Z]{ S_\mm(t) = x }
\ge
	\pr{ W_\mm(t) = w_x }
=
	\mu_{\mm,t}(w_x)
\ge
	\delta_{\mm}^{-1} / \abs{G/\mm G},
\]
recalling that $S_\mm$ lives in the quotient $G/\mm G$.
Taking the sum over all $x \in E$,
%Summing over $x \in E$,
we deduce that
\[
	1
\ge
	\sumt{x \in E}
	\pr[Z]{ S_\mm(t) = x }
\ge
	\abs E \cdot \delta_{\mm}^{-1} / \abs{G/\mm G},
\Quad{and hence}
	\abs E / \abs{G/\mm G} \le \delta_{\mm} = \oh1.
\]
Projecting onto $G/\mm G$ cannot increase the TV distance, so
\[
	\tvb{ \pr[G_k]{ S(t) \in \cdot } - \pi_G }
\ge
	\pr[Z]{ S_\mm(t) \in E } - \pi_{G/\mm G}(E)
\ge
	\prt{ \mce } - \abs E / \abs{G/\mm G}
=
	1 - \oh1.
%\qedhere
\]

Finally, recall that
\(
	\max_{\mm \in \mbn}
	\tent_\mm
=
	\max_{\mm \in \mbn}
	\tent_0(\mm, \abs{G/\mm G})
=
	\tent_*.
\)
This completes the proof.
\end{Proof}

\subsection{Upper Bound on Total-Variation Mixing}
\label{sec-p2:cutoff2:upper}

We use the same `modified $L_2$ calculation' as in \S\ref{sec-p2:cutoff1:upper}, conditioning on `typicality'; see \cref{res-p2:cutoff1:mod-l2}.

\newcommand{\ceps}{c}
Abbreviate
	$d = d(G)$
and recall that
	$\zeta_\mm = \tfrac1k (k - d) \log \mm$;
set
	$\hat \zeta_\mm \cq \zeta_\mm \wedge 1$.
Fix $\eps > 0$.
%We suppress the dependence on $\eps$ in the following.
The following depend on $\eps$; we suppress this.
Set $t \cq \tent_* (1 + \eps)$.
Recall the constant $\ceps = c_\eps > 0$~from~\cref{res-p2:cutoff2:ent:conc}.

\begin{defn}%[Typicality]
\label{def-p2:cutoff2:typ}
Set $t = \tent_* (1+\eps)$.
Define \textit{global} typicality sets for $\mm \in \mbn \cup \bra{\infty}$ by
\begin{alignat*}{3}
	\mcw_{\mm, \glo}
&
\cq
	\brb{ w \in \mbz_\mm^k \midb \prt{ W_\mm(t) = w } \le \delta_{\mm} / \abs{G/\mm G} }
&&\Quad{where}&
	\delta_{\mm} &\cq e^{-\ceps \hat \zeta_\mm k},
\intertext{%
	using the convention
	\(
		\hat \zeta_\infty
	=
		\zeta_\infty \wedge 1
%	=
%		\rbr{ \tfrac1k (k - d) \log \infty } \wedge 1
	=
		1,
	\)
	so
	\(
		\delta_\infty
	\cq
		e^{- \ceps k }.
	\)
	Define the \textit{local} typicality set by%
}
	\mcw_\loc
&\cq
	\brb{ w \in \mbz_\infty^k \midb \abs{w_i - \ext{W_{\infty,i}(t)}} \le r_* \: \forall \, i \in [k] }
&&\Quad{where}&
	r_* &\cq \tfrac12 \abs G^{1/k} \logk[2].
\intertext{%
	When $W'$ is an independent copy of $W$,
	define \textit{typicality} by%
}
	\typ
&\cq
\mathrlap{%
	\brb{ W_\infty(t), W'_\infty(t) \in \mcw_\loc }
\cap
	\rbb{ \cap_{\mm \in \MM} \brb{
		W_\mm(t), W'_\mm(t) \in \mcw_{\mm, \glo}
	} },
}
\end{alignat*}
where $\MM$ is a specific subset of $[2,\abs G]$ to be defined below in \cref{def-p2:cutoff2:Gamma-def}.
\end{defn}

We are going to use a union bound over $\mm \in \MM$, so desire control on
$\sumt{\mm \in \MM} \delta_\mm$.

\begin{lem}
\label{res-p2:cutoff2:sum-deltagamma}
	For all $\MM \subseteq \mbn \setminus \bra{1}$,
	we have
	\(
		\sumt{\mm \in \MM} \delta_\mm
	\le
		\delta_\infty \abs \MM + 2^{-c(k-d)+1}
	=
		\delta_\infty \abs \MM + \oh1.
	\)
\end{lem}

\begin{Proof}
Since $\min \MM \ge 2$ and $k - d \gg 1$,
we have
\(
	\sum_{\mm\in\MM}
	\mm^{-c(k-d)}
\le
	\sum_{\mm\ge2}
	\mm^{-c(k-d)}
\le
	2^{-c(k-d)+1}.
\)~%
So,
\[
	\sumt{\mm \in \MM}
	\delta_\mm
%=
%	\sumt{\mm \in \MM}
%	e^{- \ceps (\zeta_\mm \wedge 1) k}
\le
	\sumt{\mm \in \MM}
	\rbr{ e^{-\ceps k} + e^{-\ceps \zeta_\mm k} }
=
	e^{-\ceps k} \abs \MM
+	\sumt{\mm \in \MM}
	\mm^{-\ceps (k-d)}
\le
	\delta_\infty \abs \MM + 2^{-c(k-d)+1}.
\qedhere
\]
\end{Proof}

\begin{prop}%[Typicality]
\label{res-p2:cutoff2:typ}
	For all $\eps > 0$ and any subset $\MM \subseteq \mbn \setminus \bra{1}$,
	we have
	\[
		\prt{\typ}
	\ge
		1 - \delta_\infty \abs \MM - \oh1.
	\]
\end{prop}

\begin{Proof}
Suppress the time-dependence from the notation:
eg, write $W_\mm$ for $W_\mm(t)$ and $Q_\mm$ for $Q_\mm(t)$.

We consider global typicality first.
Observe that
\[
	Q_\mm = - \log \mu_\mm\rbr{W_\mm}
\ge
	\log \abs{G/\mm G} + \ceps \hat \zeta_\mm k
\quad
	\text{if and only if}
\quad
	\mu_\mm\rbr{W_\mm} \le e^{-\ceps \hat \zeta_\mm k} / \abs{G/\mm G}.
\]
Hence, recalling that $\delta_{\mm} = \exp{-\ceps \hat \zeta_\mm k}$ with $\hat \zeta_\mm = \zeta_\mm \wedge 1$, by \cref{res-p2:cutoff2:ent:conc}, we have
\[
	\pr{ W_\mm \notin \mcw_{\mm, \glo} }
%=
%	\pr{ \mu_\mm\rbr{W_\mm} \le \delta_\mm / \abs{G/\mm G} }
\le
	\delta_\mm,
\Quad{and hence}
	\pr{ \cap_{\mm \in \MM} \brb{ W_\mm \in \mcw_{\mm, \glo} } }
\ge
	1 - \sumt{\mm \in \MM} \delta_\mm,
\]
by the union bound.
Recall that $\zeta_\mm = \tfrac1k (k - d) \log \mm$.
Applying \cref{res-p2:cutoff2:sum-deltagamma},
we deduce that
\[
	\pr{ \cap_{\mm \in \MM} \brb{ W_\mm \in \mcw_{\mm, \glo} } }
\ge
	1 - \delta_\infty \abs \MM - \oh1.
%\Qwhere
%	\delta_\infty
%=
%	e^{- \ceps k}.
\]
%This covers the global part of typicality.

We turn to local typicality.
\cref{res-p2:cutoff2:ent:eval:order} gives $t/k \le \abs{G}^{2/k} \log k$.
So, \cref{res-p2:cutoff1:typ:loc}~gives%
\[
	\pr{ \cap_i \brb{ \abs{ W_{\infty,i} - \ext{W_{\infty,i}} } \le r_* } }
=
	1 - \oh1,
\Quad{and hence}
	\pr{ W_\infty \in \mcw_\loc }
=
	1 - \oh1.
\]
%This covers local typicality.

The claim follows by combining local and global typicality and applying the union bound.
\end{Proof}

We now choose the set $\MM$, to make sense of typicality.
Recall that $a \wr b$ means that $a$ divides $b$.
%	that ``$H \le G$'' means that $H$ is a subgroup of $G$
%and
%	that we write $\alpha \wr \beta$, for $\alpha,\beta \in \mbn$, if $\alpha$ divides $\beta$.
%Abbreviate $n_* \cq (n-1) \wedge \floor{2r_*}$.

\begin{defn}%[$\MM$ and $\mch$ Definition]
\label{def-p2:cutoff2:Gamma-def}
Define
\(
	\Delta
\cq
	\bra{ \mm \in [2, n_*] \mid \mm \wr n },
\)
where $n_* = \floor{2r_*}$.
Recall that
\[
	\mch
\cq
	\brb{ \mm G \mid \mm \in \Delta, \: \mm G \ne G }
=
	\brb{ H \mid H = \mm G \ne G \text{ for some } \mm \in \Delta }.
%	\brb{ H \mid H = \mm G \ne G \text{ for some } \mm \wr n \text{ with } 2 \le \mm \le n_* }.
\]
Given $H \in \mch$, write
\(
	\MM_H
\cq
	\bra{ \mm \in \Delta \mid H = \mm G }
\)
and denote by $\mm_H$ the minimal $\mm \wr n$ with $H = \mm G$, ie
\(
	\mm_H
\cq
	\inf \MM_H.
\)
Finally, define
\(
	\MM
\cq
	\bra{ \mm_H \mid H \in \mch } \cup \bra{ n };
\)
so,
\(
	\MM
\subseteq
	\Delta \cup \bra{ n }
\subseteq
	[2, n_*] \cup \bra{ n }.
\)
\end{defn}

The following lemma, whose proof is deferred to the end of this subsection, is also needed.

\begin{lem}%[Divisibility]
\label{res-p2:cutoff2:divis}
	For all $H \in \mch$ and all $\mm \in \MM_H$, we have $\mm_H \wr \mm$.
\end{lem}

Recall that $\tent_* = \max_\mm \tent_\mm$.
In analogy with \S\ref{sec-p2:cutoff1:upper} and \cref{res-p2:cutoff1:l2},
write
\[
	D
\cq
	D(t)
\cq
	n \, \pr{ V_\infty\rbr{ t } \bcdot Z = 0 \mid \typ } - 1.
\]

\begin{prop}%[Modified $L_2$ Distance]
\label{res-p2:cutoff2:l2}
%	Suppose that $(d,n,k)$ jointly satisfy \cref{hyp-p2:cutoff2}.
%	(Recall that, implicitly, $(d,n,k)$ is a sequence of triples of integers.)
	Write
	\(
		\mfgcd
	\cq
		\gcd(V_{\infty,1}, ..., V_{\infty,k}, n).
	\)
	Recall that $\eps > 0$ and $t = \tent_* (1 + \eps)$.
%	Then, for all $\eps \in (0,1)$, we have
	We have
	\[
		0
	\le
		D\rbb{ t(1 + \eps) }
	=
		\sumt{\mm \in \mbn}
		\prt{ \mfgcd = \mm \mid \typ } \cdot \abs{G}/\abs{\mm G}
	-	1
	\le
		\rbb{ \delta_\infty \abs \MM + \oh1 } / \prt{\typ}.
	\]
	The conditions of \cref{hyp-p2:cutoff2} imply immediately that this last term is $\oh1$, as $\delta_\infty = e^{- \ceps k}$.
\end{prop}

It is straightforward to deduce the upper bound on mixing from \cref{res-p2:cutoff2:typ,res-p2:cutoff2:l2}.

\begin{Proof}[Proof of Upper Bound in \cref{res-p2:cutoff2:res}]
We use a modified $L_2$ calculation at time $(1 + \eps) \max_\mm \tent_\mm$.

\begin{itemize}[itemsep = 0pt, topsep = \medskipamount, label = \bcdot]
	\item 
	Condition that $W$ satisfies typicality; see \cref{def-p2:cutoff2:typ,res-p2:cutoff2:typ}.
	
	\item 
	Perform the standard TV--$L_2$ bound on $S = W \bcdot Z$ conditional that $W$ is typical; cf \cref{res-p2:cutoff1:mod-l2}.
%	Perform the standard TV--$L_2$ upper bound on the law of $S$ conditioned that $W$ is typical.
	
	\item 
	Upper bound the expected $L_2$ distance by
	\(
		\rbr{ \delta_\infty \abs \MM + \oh1 } / \prt{\typ};
	\)
	see \cref{res-p2:cutoff2:l2}.
	
	\item 
	This gives an upper bound on the expected TV distance of
	\(
		\rbr{ \delta_\infty \abs \MM + \oh1 } / \prt{\typ} + \prt{\typ^c}.
	\)
	
	\item 
	Clearly, $\abs \MM \le \abs \mch + 1$.
	So, $k \gg \log(\abs \mch + 1)$ implies
	\(
		\delta_\infty \abs \MM
	\le
		\delta_\infty (\abs \mch + 1)
	\ll
		1,
	\)
	as
	\(
		\delta_\infty
	=
		e^{-\ceps k}.
%	\ll
%		1.
	\)
	Thus, $\prt{\typ} = 1 - \oh1$ by \cref{res-p2:cutoff2:typ}.
	Hence, the expected TV distance is $\oh1$.
	
	\item 
	This means that the TV distance is $\oh1$ \whp, by Markov's inequality.
\qedhere
\end{itemize}
%These calculations are all performed at time $t = (1 + \eps) \max_\mm \tent_\mm$.
%This completes the proof.
	%
\end{Proof}

We now prove \cref{res-p2:cutoff2:l2}.
All terms are evaluated at $t = \tent_* (1 + \eps)$, and this is suppressed.
%To ease exposition, while all terms are evaluated at time
%$t = (1 + \eps) \max_{\mm \in \mbn} \tent_\mm$,
%we suppress this from the notation.

\begin{Proof}[Proof of \cref{res-p2:cutoff2:l2}]
	%
%For this proof, we we use a subscript $\infty$ to emphasise when a random variable lies in $\mbz$, as opposed to some $\mbz_\mm$.
%
Write $V_\infty \cq W_\infty - W'_\infty$ and $\mfgcd \cq \gcd(V_{\infty,1}, ..., V_{\infty,k}, n)$.
If $\mfgcd = \mm$, which must have $\mm \wr n$ as the gcd is with $n$, then $V_\infty \bcdot Z \sim \Unif(\mm G)$ by \cref{res-p2:cutoff1:unif-gcd}.~%
Then,
%by the standard $L_2$ calculation,
\[
	D
=
	n \, \pr{ V_\infty \bcdot Z = 0 \mid \typ } - 1
=
	\abs G
	\sumt{\mm \wr n}
	\prt{ \mfgcd = \mm \mid \typ } / \abs{\mm G}
-	1.
\]

We consider various cases.
%For $\mm$ such that $\mm G = G$, we have $\abs{\mm G} = \abs G$ and upper bound
First, combining together all $\mm$ such that $\mm G = G$, we upper bound
\[
	\abs G \, \pr{ \mfgcd \in \bra{ \mm \mid \mm G = G } } / \abs{\mm G} \le 1.
\]
If $V_{\infty} = 0$ in $\mbz^k$, then $\mfgcd = \mm = n$, which gives $\mm G = \bra{\id}$;
using the definition of typicality,
\[
	\abs G \, \pr{ V_\infty = 0 \mid \typ } / \abs{\mm G}
=
	\abs G \, \ex{ \pr{ W_\infty = W'_\infty \mid W'_\infty, \, \typ } \mid \typ }
\le
	\delta_\infty / \prt{\typ};
\]
cf \cref{res-p2:cutoff1:V=0}.
If $V_{\infty} \ne 0$, then, given (local) typicality, $\mfgcd \le n_* = \floor{2 r_*}$.

It remains to study
\(
	\mm \in \Delta.
% \cq \bra{ \mm \in [2,n_*] \mid \mm \wr n }.
\)
By
\cref{res-p2:cutoff2:divis}, for any $H \in \mch$, we have
\[
	\brb{ V_\mm = 0 \text{ for some } \mm \in \MM_H }
\subseteq
	\brb{ V_{\mm_H} = 0 }.
\]
(Recall that $V_\mm \in \mbz_\mm^k$ for each $\mm$.)
This collapses the sum over all $\mm \in \MM_H$ into the single term $\mm_H$:
\[
&	\sumt{\mm \in \MM_H}
	\prt{ \mfgcd = \mm \mid \typ } / \abs{\mm G}
=
	\pr{ \cup_{\mm \in \MM_H} \bra{ \mfgcd = \mm } \mid \typ } / \abs H
\\&\qquad
\le
	\pr{ V_\mm = 0 \text{ for some } \mm \in \MM_H } / \abs H
\le
	\pr{ V_{\mm_H} = 0 \mid \typ } / \abs H
\le
	(\delta_{\mm_H} / \abs G) / \prt{\typ},
\]
using typicality for the final inequality.
%with the final inequality using typicality, as above.
We decompose $\sumt{\mm \in \Delta}$ into $\sumt{H \in \mch} \sumt{\mm \in \MM_H}$:
\[
	\abs G
	\sumt{\mm \in \Delta}
	\prt{ \mfgcd = \mm \mid \typ } / \abs{\mm G}
=
	\abs G
	\sumt{H \in \mch}
	\sumt{\mm \in \MM_H}
	\prt{ \mfgcd = \mm \mid \typ } / \abs{\mm G}
\le
	\sumt{H \in \mch}
	\delta_{\mm_H} / \prt{\typ}
\]

Combining all parts and using \cref{res-p2:cutoff2:sum-deltagamma},
we deduce the proposition:
\[
	0
\le
	n \, \pr{ V \bcdot Z = 0 \mid \typ } - 1
=
	\abs G
	\sumt{\mm \wr n}
	\prt{ \mfgcd = \mm \mid \typ } / \abs{\mm G}
-	1
\le
	\rbb{ \delta_\infty \abs \MM + \oh1} / \prt{\typ}.
\qedhere
\]
\end{Proof}

%It remains to give the deferred proof of the divisibility lemma (\cref{res-p2:cutoff2:divis}).
It remains to give the deferred proof \cref{res-p2:cutoff2:divis}.
Recall that $a \wr b$ means that $a$ divides $b$.

\begin{Proof}[Proof of \cref{res-p2:cutoff2:divis}]
Decompose $G = \oplus_1^r \: \mbz_{m_j}$ arbitrarily.
Fix $\beta \in \MM_H$.
Then, $H = \oplus_1^r \: h_j \mbz_{m_j}$ where $h_j \cq \gcd(\beta, m_j)$ for all $j$,
since $\alpha G = \beta G$ if and only if $\gcd(\alpha, m_j) = \gcd(\beta, m_j)$ for all $j$.
Set $\mm_* \cq \lcm(h_1, ..., h_r)$.
We show that $\mm_* G = H$ and $\mm_* \wr \alpha$ for all $\alpha \in \MM_H$, proving the lemma.

Fix $j \in [r]$.
Now, $h_j \wr \mm_* = \lcm(h_1, ..., h_r)$ and $h_j \wr m_j$ by assumption.
Hence, $h_j \wr \gcd(\mm_*, m_j)$.
Conversely, if $x \wr z$ and $y \wr z$, then $\lcm(x,y) \wr z$, and so $\mm_* = \lcm(h_1, ..., h_r) \wr \beta$ since $h_j \wr \beta$.
Hence, $\gcd(\mm_*, m_j) \wr \gcd(\beta, m_j) = h_j$.
Thus, $h_j = \gcd(\mm_*, m_j)$.
Hence, $\mm_* G = H$.
Now consider any $\alpha$ with $\alpha G = H$; so, $h_j = \gcd(\alpha, m_j)$ for all $j$.
Hence, $h_j \wr \alpha$ for all $j$, and so $\lcm(h_1, ..., h_r) \wr \alpha$, ie $\mm_* \wr \alpha$.
\end{Proof}

\section{TV Cutoff: Combining Approaches \#1 and \#2}
\label{sec-p2:cutoff3}

The only regime which we have not yet covered is
\[
	\sqrt{ \log \abs G / \log \log \log \abs G }
\lesssim
	k
\lesssim
	\sqrt{ \log \abs G }
\Quad{with}
	1
\ll
	k - d(G)
\ll
	k;
\]
see \cref{rmk-p2:cutoff1:hyp,rmk-p2:cutoff2:hyp} for the regimes covered by Approaches \#1 and \#2, respectively.
We combine the approaches here, using the refined notion of entropic times (\S\ref{sec-p2:cutoff2:ent:def}), to handle the rest.

\subsection{Precise Statements and Remarks}
\label{sec-p2:cutoff3:res}

There are some simple conditions on $k$, in terms of $d(G)$ and $\abs G$, needed for the upper bound.

\begin{hyp}
\label{hyp-p2:cutoff3}
	The sequence $(k_N, G_N)_\Ninn$ satisfies \textit{\cref{hyp-p2:cutoff3}} if
	the following hold:
	\begin{gather*}
		\LIMINF{\Ninf} k_N / \sqrt{ \log \abs{G_N} / \log \log \log \abs{G_N} } > 0,
	\quad
		\LIMSUP{\Ninf} k_N / \sqrt{ \log \abs{G_N} } < \infty,
	\\
		\LIMINF{\Ninf} \rbb{ k_N - d(G_N) } = \infty
	\Qand
		\LIMSUP{\Ninf} \rbb{ k_N - d(G_N) } / k_N = 0.
	\end{gather*}
\end{hyp}

\begin{rmkt}
\label{rmk-p2:cutoff3:hyp}
	In short, the conditions of \cref{hyp-p2:cutoff3} say that
	\[
		\sqrt{ \log \abs G / \log \log \log \abs G }
	\lesssim
		k
	\lesssim
		\sqrt{ \log \abs G}
	\Qand
		1
	\ll
		k - d(G)
	\ll
		k.
	\qedhere
	\]
%	The regime of smaller $k$ is covered by Approach \#1 and of larger $k$ by Approach \#2
\end{rmkt}

%In \cref{rmk-p2:cutoff3:hyp} below, we give a sufficient condition for \cref{hyp-p2:cutoff3} to hold.
Throughout the proofs, we drop the subscript-$N$ from the notation, eg writing $k$ or $n$, considering sequences implicitly.
Recall that we abbreviate the TV distance from uniformity at time~$t$~as
\[
	d_{G_k, N}(t) = \tvb{ \pr[{G_N([Z_1, ..., Z_{k_N}])}]{ S(t) \in \cdot } - \pi_{G_N} }
\Qwhere
	Z_1, ..., Z_{k_N} \sim^\iid \Unif(G_N).
\]

We now state the main theorem of this section.
Recall that
\(
	\tent_* = \max_{\mm \in \mbn} \tent_0(\mm,\abs{G/\mm G}).
\)

\begin{thm}%[Cutoff: Upper Bound \#3]
\label{res-p2:cutoff3:res}
	Let $(k_N)_\Ninn$ be a sequence of positive integers and $(G_N)_\Ninn$ a sequence of finite, Abelian groups;
	for each $\Ninn$, define $Z_{(N)} \cq [Z_1, ..., Z_{k_N}]$ by drawing $Z_1, ..., Z_{k_N} \sim^\iid \Unif(G_N)$.
	
	Suppose that the sequence $(k_N, G_N)_\Ninn$ satisfies \cref{hyp-p2:cutoff3}.
	Let $c \in (-1,1) \setminus \bra{0}$.
	Then
	\[
		d_{G_k, N}\rbb{ (1 + c) \tent_*(k_N, G_N) }
	\to^\mbp
		\one{c < 0}
	\Quad{(in probability)}
		\asinf N.
	\]
	That is, \whp, there is TV cutoff at $\max_\mm \tent_0(\mm, \abs{G/\mm G})$.
	Moreover, the implicit lower bound on the TV distance holds deterministically, ie for all choices of generators.
\end{thm}

\begin{rmkt*}
The TV lower bound from \S\ref{sec-p2:cutoff2:lower} is valid whenever $k - d(G) \gg 1$. Thus, it suffices to consider only the upper bound here.
The asymptotic evaluation of $\tent_*$ depends on the regime of~$k$.
\end{rmkt*}

\subsection{Outline of Proof}

Fundamentally, we still wish to bound the same expression that we did in previously:
\[
	\sumt{\mm \wr \abs G}
	\prt{ \mfgcd = \mm \mid \typ } \cdot \abs{G}/\abs{\mm G}
-	1;
\]
see \cref{res-p2:cutoff1:l2,res-p2:cutoff2:l2}.
In \S\ref{sec-p2:cutoff1:upper}, we used
\(
	\abs{G/\mm G}
\le
	\mm^{d(G)}.
\)
In \S\ref{sec-p2:cutoff2:upper}, we used typicality to get
\[
	\pr{ \mfgcd = \mm \mid \typ }
\le
	\pr{W_\mm = W'_\mm \mid \typ}
\ll
	1/\abs{G/\mm G}.
\]
The idea here, for this interim regime of $k$ near $\sqrt{\log \abs G}$, is to improve the bound $\abs{G/\mm G} \le \mm^{d(G)}$ for all but $e^{\oh k}$ of the $\mm$;
for the remaining $\mm$, we use $\abs{G/\mm G} \le \mm^{d(G)}$ and the second approach.

\subsection{Upper Bound on Total-Variation Mixing}
\label{sec-p2:cutoff3:upper}

Let $G$ be an Abelian group; set $n \cq \abs G$.
One can find a decomposition $\oplus_1^d \: \mbz_{m_j}$ of $G$ such that $d = d(G)$, the minimal size of a generating set, and $m_i \wr m_j$ for all $i \le j$.
Fix such a decomposition.

\medskip

Let $\eps > 0$ and let $t \cq (1 + \eps) \tent_*(k,G)$.
We frequently suppress the $t$ and $\eps$ dependence in the notation.
Let $c \cq c_\eps > 0$ be the constant from \cref{res-p2:cutoff2:ent:conc}.
Recall some notation from \S\ref{sec-p2:cutoff3:upper}:
\[
	\zeta_\mm
=
	\tfrac1k (k - d) \log \mm,
\quad
	\hat \zeta_\mm
=
	\zeta_\mm \wedge 1,
\quad
	\delta_\mm
=
	e^{-c \hat \zeta_\mm k}
\Qand
	r_*
=
	\tfrac12 \abs G^{1/k} \logk[2].
\]
Since $k - d \gg 1$ and $k \lesssim \sqrt{\log n}$, we have $\hat \zeta_n = 1$; set $\hat \zeta_\infty \cq 1$.
Recall that $W$ is a RW on $\mbz$ and we define $W_\mm$ by $W \mod \mm$; set $W_\infty \cq W$.
We now define typicality for this section precisely.
%We repeat the definition of typicality for convenience.

\begin{defn}[cf \cref{def-p2:cutoff2:typ}]
\label{def-p2:cutoff3:typ}
Define typical sets for $\mm \in \mbn \cup \{\infty\}$ by the following:
\begin{alignat*}{3}
	\mcw_{\mm, \glo}
&\cq
	\brb{ w \in \mbz_\mm^k \midb \prt{ W_\mm(t) = w } \le \delta_\mm / \abs{G/\mm G} }
&&\Qwhere&
	\delta_\mm &= e^{-c \hat \zeta_\mm k};
\\
	\mcw_\loc
&\cq
	\brb{ w \in \mbz^k \midb \abs{w_i - \ext{W_i(t)}} \le r_* \: \forall \, i \in [k] }
&&\Qwhere&
	r_* &= \tfrac12 n^{1/k} \logk[2].
\end{alignat*}
Choose $L$ to be the maximal integer in $[1,d]$ with $m_L \le M$ where
\[
	M
\cq
	\expb{ \sqrt{ \log n \log \log n } };
\Quad{set}
	\MM
\cq
	\brb{ r m \midb r \in [k^{1/2}], \: m \wr m_L, \: rm \wr n } \setminus \bra{1}.
\]
When $W'$ is an independent copy of $W$,
define \textit{typicality} by
\[
	\typ
\cq
	\brb{ W(t), W'(t) \in \mcw_\loc }
\cap
	\rbb{ \cap_{\mm \in \MM} \brb{
		W_\mm(t), W'_\mm(t) \in \mcw_{\mm, \glo}
	} }.
\]
\end{defn}

\begin{lem}
\label{res-p2:cutoff3:|Gamma|}
	We have
	\(
		\log \abs \MM
	\ll
		k.
	\)
	In particular, $\delta_\infty \abs \MM = \oh1$.
\end{lem}

\begin{Proof}
We have
\(
	\abs \MM
\le
	k^{1/2} \, \div m_L
\)
where $\div m$ is the number of divisors of $m \in \mbn$.
It is a standard number-theoretic result that
\(
	\log \div m
\lesssim
	\log m / \log \log m
\)
uniformly in $m \in \mbn$;
see, eg, \cite[\S 18.1]{HW:number-theory}.
By the definition of $m_L$ and the assumption that $k \gtrsim \sqrt{\log n / \log \log \log n}$, we obtain
\[
	\log \div m_L
\lesssim
	\log M / \log \log M
\lesssim
	\sqrt{ \log n \log \log n } / \log \log n
\ll
	k.
\]
%Also $\log(k^{1/2}) \asymp \log k \ll k$.
Thus,
\(
	\log \abs \MM
\ll
	k.
\)
Finally, recall that $\log(1/\delta_\infty) = c k \asymp k$.
\end{Proof}

%We use a union bound over $\mm \in \MM$, which we then bound via \cref{res-p2:cutoff3:|Gamma|}.
The following result is an immediate consequence of \cref{res-p2:cutoff2:sum-deltagamma,res-p2:cutoff2:typ,res-p2:cutoff3:|Gamma|}.
% combined with the bound \cref{res-p2:cutoff3:|Gamma|} on $\abs \MM$.

\begin{lem}[cf \cref{res-p2:cutoff2:sum-deltagamma,res-p2:cutoff2:typ}]
\label{res-p2:cutoff3:delta-typ}
	We have
	\(
		\sumt{\mm \in \MM} \delta_\mm
	=
		\oh1
	\)
	and
	\(
		\pr{\typ}
	=
		1 - \oh1.
	\)
\end{lem}

%\begin{lem}[\cref{res-p2:cutoff2:sum-deltagamma}]
%\label{res-p2:cutoff3:sum-deltagamma}
%	We have
%	\(
%		\sumt{\mm \in \MM} \delta_\mm
%%	=
%%		\delta_\infty \abs \MM + \oh1
%	=
%		\oh1.
%	\)
%\end{lem}
%
%\begin{prop}[\cref{res-p2:cutoff2:typ}]
%\label{res-p2:cutoff3:typ}
%	For all $\eps > 0$,
%	we have
%	\(
%		\prt{\typ}
%	=
%		1 - \oh1.
%	\)
%\end{prop}

Thus, by applying the modified $L_2$ calculation,
as before,
%exactly as for \cref{res-p2:cutoff1:l2,res-p2:cutoff2:l2},
it suffices to prove the following result.

\begin{prop}
\label{res-p2:cutoff3:l2}
	Let $\eps > 0$ be fixed and set $t \cq (1 + \eps) \tent_*(k,G)$.
	Then,
	\[
		\abs G \, \pr{ S = S' \mid \typ } - 1
	=
		\sumt{\mm \in \mbn}
		\abs{G/\mm G} \, \pr{ \mfgcd = \mm \mid \typ }
	-	1
	=
		\oh1.
	\]
\end{prop}

In order to prove this, we first show that $L \eqsim d \eqsim k$.

\begin{lem}
	We have
	\(
		0
	\le
		d - L
	\le
		\sqrt{ \log n / \log \log n }
	\ll
		k.
	\)
	In particular, $L \eqsim d \eqsim k$.
\end{lem}

\begin{Proof}
By definition, $L \in [1, d]$, so $L \le d$.
If $L < d$, then $m_L \le M \le m_{L+1}$.
Now, $n = m_1 \cdots m_d$, so then
\(
	M^{d-L}
\le
	m_{L+1}^{d-L}
\le
	m_{L+1} \cdots m_d
\le
	n.
\)
Recall that $k \gtrsim \sqrt{\log n / \log \log \log n}$.
Rearranging,
\[
	d-L
\le
	\log n / \log M
=
	\sqrt{ \log n / \log \log n }
\ll
	k
\eqsim
	d.
\qedhere
\]
\end{Proof}

We prove \cref{res-p2:cutoff3:l2} by separating the sum over $\mm$ into two parts according to $\MM$.

\begin{Proof}[Proof of \cref{res-p2:cutoff3:l2}]
Observe that
\(
	\abs{G/\mm G} \pr{ \mfgcd = \mm \mid \typ }
\le
	1
\)
when $\mm = 1$.
Also, $\mfgcd \wr n$.
Thus,
\[
	\sumt{\mm \in \mbn}
	\abs{G/\mm G} \, \pr{ \mfgcd = \mm \mid \typ }
-	1
\le
	\sumt{\mm \in \MM'}
	\abs{G/\mm G} \, \pr{ \mfgcd = \mm \mid \typ }
+	\sumt{\mm \in \MM}
	\abs{G/\mm G} \, \pr{ \mfgcd = \mm \mid \typ }
\]
where
\(
	\MM'
\cq
	\bra{ \mm \in [2,n] \mid \mm \wr n } \setminus \MM.
\)
We analyse these sums with Approach \#1 and \#2, respectively:
namely, we show below that both sums are $\oh1$, when $t \cq (1 + \eps) \tent_*(k,G)$ with $\eps > 0$ a constant.
%From this, the proposition follows.

\begin{itemize}
\item [\bfseries\#1]
Suppose that $\mm \in \MM'$, so $\mm \notin \MM \cup \bra{1}$.
%We have
%\(
%	\abs{G/\mm G}
%=
%	\prodt[d]{1} \gcd(\mm, m_j).
%\)
For each $j \in [L]$, we may write
\[
	\mm = r_j  \cdot \gcd(\mm, m_j)
\Qand
	m_j = r'_j \cdot \gcd(\mm, m_j)
\Qwhere
	\gcd(r_j,r'_j) = 1.
\]
%In particular,
By definition of $\MM$,
if
\(
	\mm = \tilde r \cdot m
\)
for some $m \wr m_j$, then $\tilde r > k^{1/2}$,
as $\mm \notin \MM$.
Hence,
\(
	\gcd(\mm, m_j)
=
	\mm / r_j
\le
	\mm / k^{1/2}
\)
for $j \in [L]$.
Applying this to the first $L$ terms of the product gives
\[
	\abs{G/\mm G}
=
	\prodt[d]{1} \gcd(\mm, m_j)
\le
	\mm^d / k^{L/2}.
\]
Let $\delta \in (0,1)$.
Exactly the same analysis as in the proof of \cref{res-p2:cutoff1:gcd-ex} then leads us to
\[
	\sumt{\mm \in \MM'}
	\abs{G/\mm G} \, \pr{ \mfgcd = \mm \mid \typ }
\le
	e^{2 \delta k} 2^{d+1-k}
+	2^k \delta^{d+1-k} n^{(d+1-k)/k}
+	4^k \logk[2(d+1)] / k^{L/2}.
\]
Setting $\delta \cq \tfrac14 (k - d - 1)/k$ makes the first two terms $\oh1$, as in \cref{res-p2:cutoff1:gcd-ex}.
%	\footnote{%
%		Set $\eta \cq (k - d - 1)/k$ and $\delta \cq \tfrac14 \eta$.
%		Then $k - d \gg 1$ implies that the first is $\oh1$.
%		For the second, is it not difficult to see that $\eta \ge 4 k / \log n$ suffices, ie $k - d - 1 \ge 4 k^2 / \log n$, but $4 k^2 / \log n \lesssim 1$ so $k - d \gg 1$ is sufficient}
For the third term,
\(
	4^k \logk[2(d+1)] / k^{L/2}
\ll
	1
\)
as $L \eqsim k \eqsim d$.
Hence, the sum over $\mm \in \MM'$ is $\oh1$.

\item [\bfseries\#2]
The typicality conditions set out in \cref{def-p2:cutoff3:typ} imply that
\[
	\pr{ \mfgcd = \mm \mid \typ }
\le
	\pr{ W_\mm = W'_\mm \mid \typ }
\le
	\delta_\mm / \abs{G/\mm G};
\]
cf \cref{res-p2:cutoff1:V=0}.
Hence,
combining this with \cref{res-p2:cutoff3:delta-typ},
the sum over $\mm \in \MM$ is $\oh1$.
\qedhere
%\[
%	\sumt{\mm \in \MM}
%	\abs{G/\mm G} \, \pr{ \mfgcd = \mm \mid \typ }
%\le
%	\sumt{\mm \in \MM}
%	\delta_\mm
%=
%	\oh1.
%\qedhere
%\]
	%
\end{itemize}
\end{Proof}

\section{Separation Cutoff}
\label{sec-p2:sep}

Recall that separation distance is defined by
\[
	s(t)
\cq
	\maxt{x,y} \bra{ 1 - P_t(x,y) / \pi(y) }
\Qfor
	t \ge 0,
\]
%where $P_\cdot$ is the heat kernel (ie transition probabilities)
where
	$P_t(x,y)$ is the time-$t$ transition probability from $x$ to $y$
and
	$\pi$ the invariant distribution.
We write $s_{G_k,N}$ when considering sequences $(k_N, G_N)_\Ninn$, analogously to $d_{G_k,N}$ for total variation.

\subsection{Precise Statement and Remarks}

%We now state the main theorem;
As for the previous theorems, conditions are imposed on $(k,G)$.

\begin{hyp}
\label{hyp-p2:sep}
	The sequence $(k_N, G_N)_\Ninn$ satisfies \textit{\cref{hyp-p2:sep}} if the following hold:
	\[
		\LIMINF{\Ninf} \,
		\frac{ k_N - d(G_N) }{ \max\brb{ (\log \abs{G_N} / k_N)^2, \: (\log \abs{G_N})^{1/2} } }
	=
		\infty
	\Qand
		\LIMSUP{\Ninf} \,
		\frac{ \log k_N }{ \log \abs{G_N} }
	=
		0.
	\]
\end{hyp}

\begin{rmkt}
	It is easy to check that \cref{hyp-p2:sep} is satisfied when
	$\log k \ll \log \abs G$ and
	\[
		\text{\emph{either}}
	\quad
		k \gtrsim (\log \abs G)^{3/4}
	\text{ and }
		k - d(G) \gg (\log \abs G)^{1/2}
	\Quad{\emph{or}}
		k \gg (\log \abs G)^{2/3}
	\text{ and }
		k - d(G) \asymp k.
%	\qedhere
	\]
%	\[
%		k \gtrsim (\log \abs G)^{3/4},
%	\quad
%		k - d(G) \gg (\log \abs G)^{1/2}
%	\Qand
%		\log k \ll \log \abs G
%	\quad
%		\text{(simultaneously)}.
%	\qedhere
%	\]
	In particular, the latter condition holds whenever $k - \log_2 \abs G \gtrsim \log \abs G$, eg $k \gg \log \abs G$.
\end{rmkt}

\begin{thm}%[Separation Cutoff]
\label{res-p2:sep:res}
	Let $(k_N)_\Ninn$ be a sequence of positive integers and $(G_N)_\Ninn$ a sequence of finite, Abelian groups;
	for each $\Ninn$, define $Z_{(N)} \cq [Z_1, ..., Z_{k_N}]$ by drawing $Z_1, ..., Z_{k_N} \sim^\iid \Unif(G_N)$.
	
	Suppose that the sequence $(k_N, G_N)_\Ninn$ satisfies \cref{hyp-p2:sep}.
	Let $c \in (-1,1) \setminus \bra{0}$.
	Then,
	\[
		s_{G_k,N}\rbb{ (1 + c) \tent_*(k_N, G_N) }
	\to^\mbp
		\one{ c < 0 }
	\Quad{(in probability)}
		\asinf N.
	\]
	That is, \whp there is separation cutoff at $\tent_*(k,G)$.
	Moreover, the implicit lower bound on the separation distance holds deterministically, ie for all choices of generators.
\end{thm}
	
%\begin{rmkt}
%	We show in \cref{res-p2:sep:hyp-special} that
%	\cref{hyp-p2:sep} is satisfied in the following special cases:
%	\[
%		k \asymp \log \abs G \Quad{with} k - d(G) \gg \sqrt{d(G)};
%	\qquad
%		k \gg \log \abs G.
%	\qedhere
%	\]
%\end{rmkt}

%\begin{rmkt}
%%	Write $d \cq d(G)$.
%	While we only state and prove the result for $k \gtrsim \log \abs G$ with $k - d(G) \asymp k$, the argument can be extended to larger regimes in a couple of ways:
%	\begin{itemize}[itemsep = 0pt, topsep = \smallskipamount, label = \bcdot]
%		\item 
%			$k \ll \log \abs G$ with $k - d(G) \asymp k$,
%		provided
%			$\log \abs G / k$ diverges sufficiently slowly;
%		
%		\item 
%			$k \asymp \log \abs G$ with $1 \ll k - d(G) \ll k$
%		provided
%			$k - d(G)$ diverges sufficiently rapidly.
%	\end{itemize}
%	These regimes require a little more care; we do not explore the details here.
%\end{rmkt}

\begin{rmkt*}
	Total-variation distance is a lower bound on separation distance;
	see, eg, \cite[Lemma~6.16]{LPW:markov-mixing}.
	Hence, the lower bound on separation distance is immediate from that on TV distance.
\end{rmkt*}

The proof uses the previously established TV mixing time upper bound as a building block.
%In particular, the TV mixing time gives a lower bound on the separation mixing time.
%It thus suffices to prove only the upper bound.

%\subsection{Lower Bound on Separation Mixing---to be removed}
%
%\sott{%
%	Should we just remove this whole section, and add a remark after the theorem:
%	``TV is a lower bound on separation; see, eg, \cite[Lemma~6.16]{LPW:markov-mixing}. So, the lower bound follows immediately from the TV results''?
%	We did similarly in \S\ref{sec-p2:cutoff3:res}
%}
%
%%We start with the lower bound, which is a straightforward corollary of already established results.
%
%%\begin{Proof}[Proof of \cref{res-p2:sep:res}: Lower Bound]
%%	%
%Since TV is a lower bound on separation (see, eg, \cite[Lemma~6.16]{LPW:markov-mixing}), the lower bound follows from the TV result.
%References for the TV result are as follows.
%	See \cref{res-p2:cutoff2:res}, specifically \S\ref{sec-p2:cutoff2:lower} for the lower bound on mixing, for the regime $k \asymp \log \abs G$.
%	For $k \gg \log \abs G$, TV cutoff had already been established at time $\tent_*(k,G)$; see \S\ref{sec-p2:intro:previous-work:ad-conj} and \cite[\cref{res-p0:re:app:s*:eqsim:k>>log}]{HOt:rcg:supp}.
%	\sotf{should we adjust this, so as not to be quoting our \cite{HOt:rcg:supp}?}
%%	%
%%\end{Proof}

\subsection{Upper Bound on Separation Mixing}

\begin{Proof}[Preliminaries]
\qedtriangle
	%
%For $y,z \in G$ and $t \ge 0$, write $P_t(y,z) \cq \pr[y]{ S(t) = z }$ for the time-$t$ transition probability from $y$ to $z$.
Write $n \cq \abs G$.
We want to show, for fixed $\xi > 0$, that
\[
	\mint{x \in G}
	P^\pm_t(0,x)
\ge
	\tfrac1n \rbb{ 1 - \oh1 }
\Quad{for some}
	t \le (1 + 2\xi) \tent^\pm_*(k, G).
\]

Abbreviate $d \cq d(G)$.
Let $\chi = \oh1$, to be specified later.
Throughout the course of the proof, we impose conditions on $\chi$; at the end, we show that these are equivalent to \cref{hyp-p2:sep}.

Set $k' \cq k - \chi(k - d)$; then, $k' \eqsim k$ and $k' - d = (1 - \chi) (k - d) \eqsim k - d \gg 1$.
Let $A \cq \sbr{ Z_1, ..., Z_{k'} }$ be the first $k'$ generators and $B \cq \sbr{ Z_{k' + 1}, ..., Z_k }$ be the remaining $k - k' = \chi(k - d)$.
Since $G$ is Abelian,
%we may write
\(
	P_t = P_{t,A} P_{t,B}
\)
where in $P_{t,A}$, respectively $P_{t,B}$, we pick each generator of $A$, respectively $B$, at rate $1/k$ independently.
In words, first apply the generators from $A$, then those from~$B$.
\end{Proof}

Let $\xi > 0$ be a constant; let $t' \cq (1 + \xi) \tent_*\rbr{ k', G }$.
Since there is cutoff on $G(A)$ \whp at $\tent_*\rbr{ k', G }$, we can choose $\delta = \oh1$ so that $t'$ is larger than the $\delta^2$-TV mixing time for the rate-1 RW on $G(A)$ for a typical choice of $A$.
In the regime $k \gg \log n$, simply having $\delta = \oh1$ will be sufficient.
In the regime $k \lesssim \log n$, we quantify this to be $\delta = e^{-2c(k-d)}$.
\cref{hyp-p2:sep} implies that $k \gg \sqrt{\log n}$; combined with $k \lesssim \log n$, this means that \cref{hyp-p2:cutoff2}, used in Approach \#2 (\S\ref{sec-p2:cutoff2}), is satisfied.
We also compare $\tent_*(k', G)$ and $\tent_*(k,G)$, the \whp-cutoff times for $G(A)$ and $G(Z)$, respectively.

We need two auxiliary lemmas, which we state now;
%, both of which have two parts.
%The following two auxiliary lemmas, both with two parts, have
their proofs are deferred to \S\ref{sec-p2:sep:aux}.

\begin{lem}
\label{res-p2:sep:tv-quant}
	Assume \cref{hyp-p2:sep}.
	Then, there exists $\delta \ll 1$ such that the $\delta^2$-mixing time of the RW on $G(A)$ is at most
%	$t'$
	$t' = (1 + \xi) \tent_*(k', G)$
	\whp.
	Further, when additionally $k \lesssim \log n$, there exists a constant $c > 0$ such that we may take $\delta \cq e^{-4c(k-d)}$.
\end{lem}

%\begin{subtheorem}{thm}
%\label{res-p2:sep:tv-quant}
%
%\begin{lem}
%\label{res-p2:sep:tv-quant:<}
%	Assume \cref{hyp-p2:sep}.
%	When $k \lesssim \log n$,
%	there exists a constant $c > 0$ such that, writing $\delta \cq e^{-4c(k-d)} \ll 1$,
%	the $\delta^2$-mixing time of the RW on $G(A)$ is at most $t'$~\whp.
%\end{lem}
%
%\begin{lem}
%\label{res-p2:sep:tv-quant:>}
%	When $k \gg \log n$,
%	there is some $\delta \ll 1$ such that
%	the $\delta^2$-mixing time of the RW on $G(A)$ is at most $t'$~\whp.
%\end{lem}
%
%\end{subtheorem}
%\begin{Proof}
%	We defer the proof of this auxiliary lemma to \S\ref{sec-p2:sep:aux}.
%\end{Proof}

\begin{lem}
\label{res-p2:sep:t*A-t*Z}
	We have
	\(
		\tent_*(k', G)
	\eqsim
		\tent_*(k, G)
	\)
	if and only if
	\(
		\chi (k - d) k^{-2} \log n
	\ll
		1.
	\)
\end{lem}

%\begin{subtheorem}{thm}
%\label{res-p2:sep:t*A-t*Z}
%
%\begin{lem}
%\label{res-p2:sep:t*A-t*Z:<}
%	When $k \lesssim \log n$,
%	we have
%	\(
%		\tent_*(k', G)
%	\eqsim
%		\tent_*(k, G)
%	\)
%	if and only if
%	\(
%		\chi (k - d) k^{-2} \log n
%	\ll
%		1.
%	\)
%\end{lem}
%
%\begin{lem}
%\label{res-p2:sep:t*A-t*Z:>}
%	When $k \gg \log n$,
%	we have
%	\(
%		\tent_*(k', G)
%	\eqsim
%		\tent_*(k, G).
%	\)
%%	for all $\alpha \in (0, \infty)$.
%\end{lem}
%
%\end{subtheorem}

%\begin{Proof}
%	We defer the proof of this auxiliary lemma to \S\ref{sec-p2:sep:aux}.
%\end{Proof}

Note that $k \gg \log n$ implies $k - d \eqsim k$, and so $(k - d) k^{-2} \log n \eqsim \log n / k \ll 1$ already; so, any $\chi \ll 1$ suffices.
Assume that $\chi (k - d) k^{-2} \log n \ll 1$ so that $\tent_*(k', G) \eqsim \tent_*(k, G)$.
To relate this to the rate-1 RW on $G(Z)$, rescale time by $k / \abs A = 1/(1 - \chi(k-d)/k)$: set $t \cq t' / (1 - \chi(k-d)/k)$.~%
Thus,
\[
	t \eqsim (1 + \xi) \tent_*(k, G)
\Quad{as}
	\chi \ll 1 \le k/(k-d);
\Quad{in particular,}
	t \le (1 + 2 \xi) \tent_*(k, G).
\]
By monotonicity of the separation distance with respect to time,
it thus suffices to show that
\[
	\mint{x \in G} P_t(0,x)
\ge
	\tfrac1n \rbb{ 1 - \oh1 }.
\quad
	\whp.
\]

Recall that the generators $Z = A \cup B$ are separated into $A$, the first $k'$, and $B$, the remainder.

\begin{lem}
\label{res-p2:sep:union-bound}
	Assume \cref{hyp-p2:sep}.
	Let $t \ge t' = (1 + \xi) \tent_*(k', G)$.
	Assume that
	\[
		1/(k-d) \ll \chi \ll 1
	\Qwhen
		k \lesssim \log n
	\Qand
		1/t' \ll \chi \ll 1
	\Qwhen
		k \gg \log n.
	\]
	Let $\delta$ be as in Lemma~\ref{res-p2:sep:tv-quant}.
	For $y,z \in G$, define
	\[
		Q_B(y,z)
	\cq
		\abs{B_\pm}^{-1} \sumt{b \in B_\pm} \one{ y + b^{-1} = z };
	\]
	ie, $Q_B$ is the transition matrix for the time-reversed RW on the Cayley graph $G(B)$.
	Suppose also that
%	there is some $\eta \ll 1$ such that
	for
		all (deterministic) sets $D \subseteq G$ with $\abs{G \setminus D} \le \delta \abs G$
	and
		all $x \in G$
	uniformly,
	\[
		\pr{ Q_B(x,D) \le 1 - \eta }
	=
		\oh{1/\abs G}
	\Quad{for some}
		\eta = \oh1,
	\]
	where $B^+ \cq B$ and $B^- \cq B \cup B^{-1}$ (as multisets).
%	for some $\eta = \oh1$.
	Then,
	\[
		\mint{x \in G} P_t(0,x) \ge \tfrac1n \rbb{ 1 - \oh1 }
	\quad
		\whp.
	\]
\end{lem}

\begin{Proof}
We condition on a typical realisation of $A$: write
\(
	\mca
\cq
	\bra{ a \mid \tmix\rbr{\delta^2; G(a)} \le t' }
\)
and condition on $A = a$ for a fixed $a \in \mca$.
%Let $\delta$ be as in Lemma~\ref{res-p2:sep:tv-quant}:
%\begin{itemize}[noitemsep, topsep=\smallskipamount, label = \bcdot]
%	\item 
%	$\delta = e^{-2c(k-d)}$, for some constant $c > 0$, when $k \lesssim \log n$;
%	
%	\item 
%	$\delta$ is allowed to vanish arbitrarily slowly when $k \gg \log n$.
%\end{itemize}
Then,
\(
	\pr{ A \in \mca }
=
	1 - \oh1
\)
by Lemma~\ref{res-p2:sep:tv-quant}.
%Indicate this by writing $\pr[A=a]{\cdot}$.
Given $A = a \in \mca$,
the set
\[
	D
\cq
	\brb{ z \in G \midb P_{t,a}(0,z) \ge \tfrac1n \rbr{1 - \delta} }
\Quad{satisfies}
	\abs D
\ge
	n (1 - \delta).
\]
Indeed, using the distinguishing-statistic representation of total-variation distance,
\[
	\delta^2
\ge
	\pi(D^c) - P_{t,a}(0,D^c)
\ge
	\tfrac1n \abs{D^c} - \tfrac1n (1 - \delta) \abs{D^c}
=
	\tfrac1n \delta \abs{D^c},
\Quad{so}
	\tfrac1n \abs{D^c}
\le
	\delta.
\]
For the undirected case (ie the RW on $G^-_k$),
by reversibility,
conditional on $A = a \in \mca$,
we have
\[
	P^-_t(0,x)
\ge
	P^-_{t,B}(x,D) \cdot \tfrac1n \rbr{ 1 - \delta }.
\]
While the RW on $G^+_k$ is not reversible, Cayley graphs have the special property that a step `backwards' with a generator $z$ corresponds to a step `forwards' with $z^{-1}$.
Thus,
\[
	P^+_t(0,x)
\ge
	Q^+_{t,B}(x,D) \cdot \tfrac1n \rbr{ 1 - \delta }
\]
where $Q^+_{\cdot,B}$ is the heat kernel for the RW on $G^+(B^{-1})$ where $B^{-1} \cq [ z^{-1} \mid z \in B ]$, rather than on $G^+(B)$.
For the RW on $G^-_k$, replacing the generators with their inverses has no effect on the graph (or RW); set $Q^-_{\cdot,B} \cq P^-_{\cdot,B}$.
We want to show that
\(
	Q_{t,B}(x,D) = 1 - \oh1
\)
uniformly in $x \in G$ \whp.
% (over the randomness in $B$).

Now, $Q_{t,B}$ corresponds to a RW on $G^\pm(B^{-1})$ run for time $t$.
By considering just the final step of this RW, we now argue that the hypothesis of the lemma is sufficient.
Indeed,
first note that
\[
	\mint{x}
	Q_{t,B}(x,D)
\ge
	\rbb{ 1 - e^{-t \abs B / k} }
\cdot
	\mint{x}
	Q_B(x,D),
\]
where $e^{-t \abs B / k}$ is the probability that none of the generators in $B$ are applied by time $t$.
%We assume that $\chi$ decays slowly enough so that $\chi (k - d) \gg 1$.
%Then, $t \abs B / k = \chi t (k - d) / k$.
%
%If $k \lesssim \log n$, then the condition $(k - d)^2 \gg \log n$ from \cref{hyp-p2:sep} and that on $\chi$ from the lemma implies that
%\cref{hyp-p2:sep}
%\(
%	t \abs B / k = \chi t (k - d) / k \gg 1.
%\)
%To see this, first note that $t' \gg 1 $ as $\log k' \ll \log n$.
\begin{itemize}[noitemsep, topsep = \smallskipamount, label = \bcdot]
	\item 
	If $k \gg \log n$, then $(k - d)/k \eqsim 1$,
	so $t \abs B / k \eqsim \chi t \gg 1$ whenever $1/ t' \ll \chi \ll 1$.
	
	\item 
	If $k \lesssim \log n$, then $t \ge t' \gtrsim k$,
	so $t \abs B / k \gtrsim \chi (k - d) \gg 1$ whenever $1/(k-d) \ll \chi \ll 1$.
\end{itemize}
Thus,
the assumptions in the lemma allow us to perform a union bound over $x \in G$:
\[
	\pr{ \mint{x} Q_{t,B}(x,D) \le 1 - 2 \eta \midb A = a }
=
	\oh1,
\]
where the randomness is over the generators $B$, provided $\eta$ decays sufficiently slowly.

For each $a \in \mca$. we have the desired lower bound on $\mint{x} P_t(0,x)$ conditional on $A = a$.
Finally, we average over $A$ and use $\pr{A \in \mca} = 1 - \oh1$ to complete the proof:
\[
	\pr{ \mint{x} P_t(0,x) \ge \tfrac1n \rbb{ 1 - \oh1 } }
\ge
	\pr{ \mint{x} P_t(0,x) \ge \tfrac1n \rbb{ 1 - \oh1 } \midb A \in \mca }
	\pr{ A \in \mca }
=
	1 - \oh1.
\qedhere
\]
%\(
%	\mint{x} P_t(0,x) \ge \tfrac1n \rbr{ 1 - \oh1 }
%\)
%\whp.
	%
\end{Proof}

%Next, we find conditions under which the supposition of the lemma is satisfiable.
We need to check that the supposition of the previous lemma is satisfiable.

\begin{lem}
\label{res-p2:sep:o(1/n)-bound}
%	Recall the definition of $\delta$ and of $Q_B$.
	Suppose that $\chi (k - d)^2 \gg \log n$ if $k \lesssim \log n$ and that $\chi \ll 1$ sufficiently slowly if $k \gg \log n$.
	Then, we can choose $\eta \ll 1$ vanishing sufficiently slowly so that
	for
		all (deterministic) sets $D \subseteq G$ with $\abs{G \setminus D} \le \delta \abs G$
	and
		all $x \in G$
	uniformly,
	\[
		\pr{ Q_B(x,D) \le 1 - \eta } = \oh{1/\abs G}.
%	\Qwhere
%		Q_B(y,z) \cq \abs{B_\pm}^{-1} \sumt{b \in B_\pm} \one{ y + b^{-1} = z }
	\]
%	for $y,z \in G$ where $B_+ \cq B$ and $B_- \cq B \cup B^{-1}$ (as multisets).
\end{lem}

\begin{Proof}
Fix an arbitrary $x \in G$.
We desire a proportion at least $1 - \eta$ of the generators in $B$ to connect $x$ to $D$.
The generators are chosen independently, and each connect with probability $\abs D / \abs G \ge 1 - \delta$.
Since there are $\chi (k-d)$ generators, it thus suffices to choose $\eta \ll 1$ so that
\[
	\pr{ \Bin(\chi (k-d), 1 - \delta) \le \chi (k-d) (1 - \eta) }
=
	\oh{1/\abs G}.
\]
Let $L \cq \chi (k - d)$; then, $L \gg 1$.
Direct calculation, using standard inequalities, gives
\[
	\pr{ \Bin(L, 1 - \delta) \le L (1 - \eta) }
=
	\pr{ \Bin(L, \delta) \ge \eta L }
\le
	\binomt{ L }{ \eta L } \delta^{\eta L}
\le
	(\delta e/\eta)^{\eta L}
=
	(\delta e/\eta)^{\eta \chi (k - d)}.
\]
%We require this to be $\oh{1/n}$.

\begin{itemize}[itemsep = 0pt, topsep = \smallskipamount, label = \bcdot]
	\item 
	Consider first $k \gg \log n$; necessarily, $k - d \eqsim k$.
	Here, we do not quantify $\delta$: we simply assume $\delta = \oh1$.
	Requiring $\eta$ and $\chi$ to vanish sufficiently slowly (compared with $\delta$) gives
	\(
		(\delta e / \eta)^{\eta \chi}
	=
		\oh1.
	\)
	Raising this to the power $k - d \eqsim k \gg \log n$ gives super-polynomial decay.
	
	\item 
	Consider now $k \lesssim \log n$.
	Here, $\delta = e^{-2c(k-d)}$.
	Choosing $\eta \ge e \sqrt \delta = e^{-c(k-d)+1}$ gives
	\[
		(\delta e/\eta)^{\eta \chi (k - d)}
	=
		\delta^{\eta \chi (k-d)/2}
	=
		\exp{ - c \eta \chi (k - d)^2 }.
	\]
	We can choose $\eta \ll 1$ so that
	\(
		\eta \chi (k - d)^2 \gg \log n,
	\)
	giving super-polynomial decay.
\end{itemize}

These bounds is independent of $x$, and hence holds for all $x \in G$ uniformly, as required.
\end{Proof}

On top of \cref{hyp-p2:sep},
we need $\chi$ to satisfy the following simultaneously when $k \lesssim \log n$:
\begin{itemize}[noitemsep, topsep = \smallskipamount, label = \bcdot]
	\item 
	we need $\chi (k-d) k^{-2} \log n \ll 1$ to hold for \cref{res-p2:sep:t*A-t*Z};
	
	\item 
	we need
	$1/(k-d) \ll 1 \ll \chi$ and $\chi (k-d)^2 \gg \log n$ to hold for \cref{res-p2:sep:union-bound,res-p2:sep:o(1/n)-bound}.
\end{itemize}
No requirements beyond ``sufficiently slowly'', in terms of $\delta$, are imposed on $\chi$ when $k \gg \log n$.
The next lemma states that having such a $\chi$ is equivalent to \cref{hyp-p2:sep}; its proof is deferred~to~\S\ref{sec-p2:sep:aux}.

%The following lemma shows that the conditions $\chi (k-d) k^{-2} \log n \ll 1$ and $\chi (k-d)^2 \gg \log n$ are satisfied under \cref{hyp-p2:sep}.
%Its proof is also deferred to \S\ref{sec-p2:sep:aux}.

\begin{lem}
\label{res-p2:sep:ass-hyp}
	Assume \cref{hyp-p2:sep} and that $k \lesssim \log n$.
	Then, we can choose $\chi \in (0,1)$ satisfying
	\[
		\chi (k - d) k^{-2} \log n \ll 1,
	\quad
		\chi (k - d)^2 \gg \log n
	\Qand
		1/(k-d) \ll \chi \ll 1.
	\]
	In fact, \cref{hyp-p2:sep} are equivalent to being able to find such a $\chi$.
\end{lem}

%\begin{Proof}
%	We defer the proof of this auxiliary lemma to \S\ref{sec-p2:sep:aux}.
%\end{Proof}

These lemmas combine easily to establish the upper bound in \cref{res-p2:sep:res}, as we expose now.

\begin{Proof}[Proof of Upper Bound in \cref{res-p2:sep:res}]
Assume that $1 \ll \log k \ll \log n$.
\cref{res-p2:sep:ass-hyp} guarantees that a $\chi$ satisfying the conditions of Lemmas~\ref{res-p2:sep:t*A-t*Z}, \ref{res-p2:sep:union-bound} and~\ref{res-p2:sep:o(1/n)-bound} simultaneously can be found under \cref{hyp-p2:sep} when $k \lesssim \log n$;
when $k \gg \log n$, simply take $\chi \ll 1$ sufficiently slowly.

The conclusion of these lemmas is that the separation distance is $\oh1$ \whp at time $t'$, and that $t' = (1 + \xi) \tent_*(k', G) \le (1 + 2\xi) \tent_*(k, G)$.
This completes the proof of the upper bound.
\end{Proof}

%\cref{res-p2:sep:union-bound,res-p2:sep:o(1/n)-bound} combine to imply that the separation distance is $\oh1$ \whp at time $t' = (1 + \xi) \tent_*(k', G)$, under \cref{hyp-p2:sep}.
%$\chi$ can be chosen to satisfy the additional conditions in these lemmas \cref{hyp-p2:sep} by \cref{res-p2:sep:ass-hyp} when $k \lesssim \log n$, and by simply taking it to vanish sufficiently slowly when $k \gg \log n$.
%
%Suppose $k \lesssim \log n$.
%By \cref{res-p2:sep:ass-hyp}, we can choose $\chi$ with
%\[
%	\chi (k - d) k^{-2} \log n \ll 1,
%\quad
%	\chi (k - d)^2 \gg \log n
%\Qand
%	1/(k-d) \ll \chi \ll 1
%\]
%under \cref{hyp-p2:sep}.
%Then, the separation distance is $\oh1$ \whp at time $t' = (1 + \xi) \tent_*(k', G)$ by \cref{res-p2:sep:union-bound,res-p2:sep:o(1/n)-bound}.
%Finally, $t' \eqsim \tent_*(k, G)$ by \cref{res-p2:sep:t*A-t*Z:<}, which completes the proof of $k \lesssim \log n$.
%
%Suppose now that $k \gg \log n$.
%We do not need to invoke \cref{hyp-p2:sep} in order to apply \cref{res-p2:sep:union-bound,res-p2:sep:o(1/n)-bound} here%
%	---although, \cref{hyp-p2:sep} always holds when $k \gg \log n$ if also $\log k \ll \log n$.
%As above, these imply that the separation distance is $\oh1$ \whp at $t'$.
%Finally, $t' \eqsim \tent_*(k, G)$.

\subsection{Auxiliary Lemmas}
\label{sec-p2:sep:aux}

It remains to give the deferred proofs of the auxiliary lemmas:
	Lemmas~\ref{res-p2:sep:tv-quant}, \ref{res-p2:sep:t*A-t*Z} and \ref{res-p2:sep:ass-hyp}.

\begin{Proof}[Proof of Lemma~\ref{res-p2:sep:tv-quant}]
Consider $k' \eqsim k \lesssim \log n$ first.
We use Approach \#2 (\S\ref{sec-p2:cutoff2:upper}), applied to $G(A)$; recall that $A$ has $k'$ iid generators, and that $t' = (1 + \xi) \tent_*(k', G)$.
Quantifying the \textit{Proof of Upper Bound in \cref{res-p2:cutoff2:res}},
%(Approach \#2, \S\ref{sec-p2:cutoff2:upper}),
using \cref{res-p2:cutoff2:sum-deltagamma,res-p2:cutoff2:l2},
gives a bound of
\[
	\rbb{ \delta_\infty (\abs \MM + 1) + 2^{-c(k-d)+1} }
/	\prt{\typ}
+	\prt{\typ^c}
\]
on the expected TV distance at time $t'$, under \cref{hyp-p2:cutoff2} with typicality given by \cref{def-p2:cutoff2:typ}.

The proof of \cref{res-p2:cutoff2:typ} shows that \emph{global} typicality fails with probability at most $\delta_\infty \abs \MM + e^{-c(k-d)+1}$.
The \emph{local} conditions as stated in \cref{def-p2:cutoff2:typ} do not fail with sufficiently low probability---only at most $\eta$, for some $\eta = \oh1$.
This is proved via standard large deviation estimates from \cite[\S\ref{sec-p0:rp}]{HOt:rcg:supp}.
Replacing $r_*$ with $r_* \log n$ gives failure probability $\eta^{\log n} \le e^{-k}$.

Increasing $r_*$ like this increases $\abs \MM$, but not enough to cause any issues. Indeed, we have
\[
	\log \abs \MM
\le
	\log n_*
\le
	\log\rbb{ n^{1/k} (\log k)^2 \log n }
\le
	\tfrac1k \log n + 3 \log \log n.
\]
But, \cref{hyp-p2:sep} implies that $k - d \gtrsim (\log n)^{1/2}$ and that
\[
	k
\ge
	k - d
\gg
	(\log n / k)^2,
\Quad{so}
	k
\gg
	(\log n)^{2/3};
\Quad{hence,}
	\log \abs \MM
\ll
	(\log n)^{1/3}
\ll
	k - d.
\]
Thus, the dominating term in the upper bound is $2^{-c(k-d)}$.
Adjusting $c$, this completes the case.

The case $k' \eqsim k \gg \log n$ is trivial, since here all Abelian groups have TV cutoff at $\tent_*(k', G)$.
\end{Proof}

\begin{Proof}[Proof of Lemma~\ref{res-p2:sep:t*A-t*Z}]
We have $k \eqsim k'$ and $k - d \eqsim k' - d$.
Observe that $n^{2/k} \eqsim n^{2/k'}$ if and only if
\[
	1
\gg
	\rbb{ \tfrac2{k'} - \tfrac2k } \log n
=
	\rbb{ \tfrac2{k - \chi(k-d)} - \tfrac2k } \log n,
\Quad{ie}
	\chi (k - d) k^{-2} \log n
\ll
	1.
\]
The claim  follows by \cref{res-p2:cutoff2:ent:eval:order} for $1 \ll k \lesssim \log n$.
On the other hand,
if $k \gg \log n$, then
\[
	\tent_*(k,G)
\eqsim
	T(k,n)
\cq
	\log n \mathrel{/} \log(k/\log n);
\]
see \S\ref{sec-p2:intro:previous-work:ad-conj}.
Hence, $T(k,n) \eqsim T(\alpha k,n)$ for all $\alpha \in (0,\infty)$.
Thus, $T(k,n) \eqsim T(k',n)$ as $k \eqsim k'$.
\end{Proof}
%\[
%	T(k,n)
%\cq
%	\tfrac{\rho}{\rho-1} \log n / \log k
%\Qwhere
%	\rho \cq \log k / \log \log n.
%\]
%Simple algebraic manipulations give
%\[
%	T(k,n)
%=
%	\frac1{\log k / \log\log n - 1} \frac{\log n}{\log \log n}
%=
%	\frac{\log n}{\log k  - \log \log n}
%=
%	\frac{\log n}{\log(k/\log n)}.
%\]

\begin{Proof}[Proof of \cref{res-p2:sep:ass-hyp}]
Rearranging the conditions, they are equivalent to having
\[
	\sqrt{ \log n / \chi }
\ll
	k - d
\ll
	k^2 / \rbb{\chi \log n}
\Quad{for some}
	\chi \in (0,1)
\Quad{with}
	1/\rbr{k-d} \ll \chi \ll 1.
\]
%By replacing $\chi$ with $\chi \omega$ for some $\omega \gg 1$ diverging sufficiently slowly, replacing the relation $\chi \ge 1/(k-d)$ with $\chi \gg 1/(k-d)$ gives an equivalent set of conditions.
We reparametrise these conditions.
Let $\eps \in (0,\infty)$ and set
\[
	\chi \cq \frac{\eps k^2}{(k - d) \log n};
\Quad{then}
	\sqrt{ \frac{\log n}{\chi} } = \frac{ \sqrt{k - d} \log n }{\sqrt \eps k}.
\]
The conditions on $\chi$ then, in terms of $\eps$, become
\[
	\frac{(\log n)^2}{(k-d) k^2} \ll \eps \ll 1
\Qand
	\frac{\log n}{k^2} \ll \eps \ll \frac{(k - d) \log n}{k^2}.
\]
We can find such an $\eps \in (0,\infty)$, implicitly a sequence, if and only if
\[
	\max\brbb{
		\frac{(\log n)^2}{(k-d) k^2},
	\:	\frac{\log n}{k^2}
	}
\ll
	\min\brbb{
		1,
	\:	\frac{(k - d) \log n}{k^2}
	}.
\]
Some case analysis shows that this condition is equivalent to the first condition of \cref{hyp-p2:sep}.
\end{Proof}

\section{Nilpotent Groups: Mixing Comparison and Expansion}
\label{sec-p2:comp}
\renewcommand{\mm}{\ensuremath{\gamma}}

We compare the mixing times for a nilpotent group $G$ with a `corresponding' Abelian group $\widebar G$:
	\(
		\tmix(G_k) / \tmix(\widebar G_k) \le 1 + \oh1.
	\)
We use this to show that $G_k$ is an expander \whp if $k - d(\widebar G) \gtrsim \log \abs G$.

To emphasise, the material in this section applies to both the un- and directed cases, simultaneously.
Throughout, $(G_{(\ell)})_{\ell\ge0}$ is the lower central series of $G$ and $L \cq \min\bra{ \ell \ge 0 \mid G_{(\ell)} = \bra{\id} }$.

\subsection{Precise Statements}

We prove \cref{res-p2:intro:comp:nil-abe}, which we recall here for the reader's convenience as \cref{res-p2:comp:res}.

\begin{thm}
\label{res-p2:comp:res}
	Let $G$ be a nilpotent group.
	Set
	$\widebar G \cq \oplus_1^L \: (G_{(\ell-1)}/G_{(\ell)})$.
%	where $(G_{(\ell)})_{\ell\ge0}$ is the lower central series of $G$ and $L \cq \min\bra{ \ell \ge 0 \mid G_{(\ell)} = \bra{\id} }$.
	Suppose that $1 \ll \log k \ll \log \abs G$ and $k - d(\widebar G) \gg 1$.
	Let $\eps > 0$ and let $t \ge (1 + \eps) \tent_*(k,\widebar G)$.
	Then,
	\(
		d_{G_k}(t) = \oh1
	\)~%
	\whp.
\end{thm}

\begin{rmkt*}
	An equivalent bound, depending only on $k$ and $\abs G = \abs{\widebar G}$, valid for all groups has already been established when $k \gg \log \abs G$; recall \S\ref{sec-p2:intro:previous-work:ad-conj}.
%	\cref{rmk-p2:intro:cutoff:k>>logn}.
	Thus, we need only consider~$1 \ll k \lesssim \log \abs G$.
\end{rmkt*}

We use this mixing time bound to show that $G_k$ for nilpotent $G$ is an expander
%, ie has spectral gap order 1,
\whp when $k - d(\widebar G) \gtrsim \log \abs G$.
The isoperimetric constant was defined in \cref{def-p2:intro:isoperimetric} for $d$-regular graphs:
\[
	\Phi_*
\cq
	\MIN{1 \le \abs S \le \frac12 \abs V}
	\Phi(S)
\Qwhere
	\Phi(S)
\cq
	\tfrac1{d \abs S} \absb{ \sbb{ \bra{ a, b } \in E \midb a \in S, \, b \in S^c } }.
\]
Specifically, we prove \cref{res-p2:intro:gap}, which we recall here for the reader's~convenience as \cref{res-p2:gap:res}.

\begin{thm}
\label{res-p2:gap:res}
	Let $G$ be a nilpotent group.
	Set
	$\widebar G \cq \oplus_1^L \: (G_{(\ell-1)}/G_{(\ell)})$.
%	where $(G_{(\ell)})_{\ell\ge0}$ is the lower central series of $G$ and $L \cq \min\bra{ \ell \ge 0 \mid G_{(\ell)} = \bra{\id} }$.
	Suppose that $k - d(\widebar G) \gtrsim \log \abs G$.
	Then, $\Phi_*(G_k) \asymp 1$~\whp.
\end{thm}

%We prove \cref{res-p2:comp:res} in \S\ref{sec-p2:comp:outline}--\S\ref{sec-p2:comp:eval} and \cref{res-p2:gap:res} in \S\ref{sec-p2:gap}.

\subsection{Outline of Proof}
\label{sec-p2:comp:outline}

%Let $L$ be the minimal integer such that $G_L$ is the trivial group.
Consider the series of quotients $(Q_\ell \cq G_{(\ell-1)}/G_{(\ell)})_{\ell=1}^L$.
For each $\ell \in [L]$, choose a set $R_\ell \subseteq G_{(\ell-1)}$ of coset representatives for $Q_\ell = G_{(\ell-1)}/G_{(\ell)}$.
To sample $Z_i \sim \Unif(G)$, it suffices to sample $Z_{i,\ell} \sim \Unif(R_\ell)$ for each $\ell$ independently, and then take the product: $Z_i \cq Z_{i,1} \cdots Z_{i,L}$; see \cref{res-p2:comp:unif-rep}.
Then $Z_{i,\ell} G_{(\ell)} \sim \Unif(Q_\ell)$ independently for each $i$ and $\ell$; see \cref{res-p2:comp:sample-Z}.

Suppose that $M$ steps are taken; let $\sigma : [M] \to [k]$ indicate which generator is used in each step.
%(Consider just the directed case for now.)
Set $S \cq \prodt[M]{m=1} Z_{\sigma(m)}$.
For each $\ell \in [L]$,
let $S_\ell \cq \prodt[M]{m=1} Z_{\sigma(m),\ell}$; this is the projection of $S$ to $Q_\ell$.
Then, each $S_\ell G_{(\ell)}$ is a RW on $Q_\ell$, which is an Abelian group, but all using the same choice $\sigma$.

These are RWs on Abelian groups, so the ordering in $\sigma$ does not matter.
For each $i \in [k]$, let $W_i$ be the number of times in $\sigma$ that generator $Z_i$ has been applied minus the number of times that $Z_i^{-1}$ has been applied.
Let $\sigma'$ be an independent copy of $\sigma$ and define $S'$ and $W'$ via $\sigma'$ and $Z$; for each $\ell \in [L]$, define $S'_\ell \cq \prodt[M]{m=1} Z_{\sigma(m),\ell}$.
Then, $S$ and $S'$ are iid conditional on $Z$.

To compare the RW on the nilpotent group with one on an Abelian group, we show~that
\begin{gather*}
	n \, \pr{ S = S' \mid (W, W') }
\le
	n \, \prodt[L]{1}
	\pr{ S_\ell G_{(\ell)} = S'_\ell G_{(\ell)} \mid (W,W') }
=
	\absb{ \widebar G / \mfgcd \widebar G },
\\
	\text{where}
\quad
	\mfgcd \cq \gcd(W_1 - W'_1, ..., W_k - W'_k, n);
\end{gather*}
see \cref{res-p2:comp:prod-decomp,res-p2:comp:l2-gcd}.
By analysing $\abs{ \widebar G / \mfgcd \widebar G }$, we showed in \S\ref{sec-p2:cutoff1}--\S\ref{sec-p2:cutoff3} that the RW on $\widebar G_k$ is mixed \whp shortly after $\tent_*(k, \widebar G)$; see specifically \cref{res-p2:cutoff1:unif-gcd}.
From this and the inequality above, we are able to deduce that the RW on $G_k$ is mixed \whp shortly after the same time.

\subsection{Reduction to Abelian-Type Calculations}
\label{sec-p2:comp:prelim}

%Let $L$ be the minimal integer such that $G_L$ is the trivial group.
Consider the series of quotients $(Q_\ell \cq G_{(\ell-1)}/G_{(\ell)})_{\ell=1}^L$.
For each $\ell \in [L]$, choose a set $R_\ell \subseteq G_{(\ell-1)}$ of coset representatives for $Q_\ell = G_{(\ell-1)}/G_{(\ell)}$, ie a set $R_\ell$ with $\abs{R_\ell} = \abs{Q_\ell}$ and $\bra{r G_{(\ell)}}_{r \in R_\ell} = G_{(\ell-1)}/G_{(\ell)} = Q_\ell$.
We sample the uniform generators via uniform random variables on each of the quotients.
In this way, projecting to one of the quotients, we get a RW on this quotient.

\begin{lem}%[{\cite[\cref{res-p0:nil:unif-rep}]{HOt:rcg:supp}}]
\label{res-p2:comp:unif-rep}
	For each $\ell \in [L]$,
	let $Y_\ell \sim \Unif(R_\ell)$ independently.
	Then, $Y \cq Y_1 \cdots Y_L \sim \Unif(G)$.
\end{lem}

\begin{Proof}
Let $r_0 \in G$ and consider the event $\bra{Y = r_0}$.
If $r_0 = Y_1 \cdots Y_L$, then $r_1 \cq Y_1^{-1} r_0 = Y_2 \cdots Y_L$.
Clearly, the right-hand side is in $G_{(1)}$, and so the left-hand side is too.
Hence, $r_0 \equiv Y_1$ mod~$G_{(1)}$.
% ie $\pi_1(r_0) = Y_1$.
But, $Y_1 \sim \Unif(R_1)$, so the probability of this is $1/\abs{R_1} = 1/\abs{G_{(0)}/G_{(1)}}$.
Similarly, $r_2 \cq Y_2^{-1} r_1 = Y_3 \cdots Y_L \in G_{(2)}$, so $r_1 \equiv Y_2$ mod $G_{(2)}$, the probability of which is $1/\abs{R_2} = 1/\abs{G_{(1)}/G_{(2)}}$.
Iterating,
%this argument, recalling that the $Y_\ell$ are independent, we deduce that
\[
	\pr{ Y = r_0 }
=
	\prodt[L]{1} 1/\abs{G_{(\ell-1)}/G_{(\ell)}}
=
	\prodt[L]{1} \abs{G_{(\ell)}} / \abs{G_{(\ell-1)}}
=
	\abs{G_{(L)}} / \abs{G_{(0)}}
=
	1/\abs G,
\]
since the $Y_\ell$ are independent.
Since $r_0 \in G$ was arbitrary, we deduce that $Y \sim \Unif(G)$.
\end{Proof}

This gives the following corollary.

\begin{cor}%[{\cite[\cref{res-p0:nil:sample-Z}]{HOt:rcg:supp}}]
\label{res-p2:comp:sample-Z}
	For all $(i, \ell) \in [k] \times [L]$,
	sample $Z_{i,\ell} \sim \Unif(R_\ell)$ independently and set $Z_i \cq Z_{i,1} \cdots Z_{i,L}$.
	Then, $Z_1, ..., Z_L \sim^\iid \Unif(G)$.
	Further, $Z_{i,\ell} G_{(\ell)} \sim \Unif(Q_\ell)$ independently for all~$(i, \ell)$.
%	$(i, \ell) \in [k] \times [L]$.
\end{cor}

\begin{Proof}
All the independence claims are immediate.
The first claim is immediate from \cref{res-p2:comp:unif-rep}.

For the second claim, we have $Z_{i,\ell} \sim \Unif(R_\ell)$ and $\abs{R_\ell} = \abs{Q_\ell}$.
Now, $x G_{(\ell)} = y G_{(\ell)}$ if and only if $y^{-1} x G_{(\ell)} = G_{(\ell)}$.
If $X \sim \Unif(R_\ell)$ and $H \in Q_\ell$, say $H = y G_{(\ell)}$ with $y \in R_\ell$, then $y^{-1} X \sim \Unif(R_\ell)$ independently of $y$.
So, $\pr{ X G_{(\ell)} = y G_{(\ell)} } = 1/\abs{R_\ell}$.
Hence, $X G_{(\ell)} \sim \Unif(Q_\ell)$.
\end{Proof}

Assume that $Z$ is drawn in this way for the remainder of the section.
The next main result (\cref{res-p2:comp:prod-decomp}) is the key element of the proof of \cref{res-p2:comp:res}.
It reduces the problem to a collection of Abelian calculations, the like of which were handled previously.

%As is standard, we write $0$ for the identity of an Abelian group.
First, we need an auxiliary lemma, showing that $v \bcdot Z = 0$ is the `worst case'.

\begin{lem}
\label{res-p2:comp:abe-worst-0}
	Let $H$ be an Abelian group.
%	, with identity $0$.
	Let $Z_1, ..., Z_k \sim^\iid \Unif(H)$ and $v \in \mbz^k$.
%	Write $v \bcdot Z \cq \sumt[k]{1} v_i Z_i$.
	Then,
	\[
		\maxt{h \in H}
		\pr{ v \bcdot Z = h }
%		\pr{ \sumt[k]{1} v_i Z_i = h }
	=
		\pr{ v \bcdot Z = \id(H) }.
%		\pr{ \sumt[k]{1} v_i Z_i = 0 }.
	\]
\end{lem}

\begin{Proof}
	%
%In essence, we show that every solution to $v \bcdot z = h$ gives rise to a to solution $v \bcdot z = 0$.
Write $0 \cq \id(H)$.
Let $h \in H$.
Write
\(
	A(h)
\cq
	\bra{ z \in H^k \mid v \bcdot z = h }.
\)
If $w \in A(h)$, then
\(
	B
\cq
	\bra{ z - w \mid z \in A(h) }
\subseteq
	A(0);
\)
also, clearly, $\abs B = \abs{A(h)}$, so $\abs{A(h)} \le \abs{A(0)}$.
Hence,
\[
	\pr{ v \bcdot Z = h }
=
	\abs{A(h)} / \abs{H}^k
\le
	\abs{A(0)} / \abs{H}^k
=
	\pr{ v \bcdot Z = 0 }.
\qedhere
\]
\end{Proof}

The following theorem decomposes the probability $\pr{S = S'}$ into a product of probabilities $\pr{ S_\ell G_{(\ell)} = S'_\ell G_{(\ell)} }$.
The latter correspond to probabilities of RWs on \emph{Abelian groups} $Q_{\ell}$.
%The proof uses crucially the nilpotency of the group.

\begin{prop}
\label{res-p2:comp:prod-decomp}
	Let $M, M' \in \mbn$.
	Let $\sigma : [M] \to [k]$ and $\sigma' : [M'] \to [k]$.
	Let $\eta \in \bra{\pm1}^M$ and $\eta' \in \bra{\pm1}^{M'}$.
	Recall that $(G_{(\ell)})_{\ell\ge0}$ is the lower central series and
	that $G_{(L)} = \bra{\id}$.
%	$L = \min\bra{ \ell \ge 0 \mid G_{(\ell)} = \bra{\id} }$.
	For $\ell \in [L]$,
	set
	\[
		S_\ell \cq \prodt[M]{m=1} Z_{\sigma(m), \ell}^{\eta_m},
	\quad
		S'_\ell \cq \prodt[M]{m=1} Z_{\sigma'(m), \ell}^{\eta'_m},
	\quad
		S \cq \prodt[M]{m=1} Z_{\sigma(m)}^{\eta_m}
	\Quad{and}
		S' \cq \prodt[M]{m=1} Z_{\sigma'(m)}^{\eta'_m}.
	\]
	For $i \in [k]$,
	write
	\(
		v_i
	\cq
		\sumt{m \in [M'] : \sigma'(m) = i} \eta'_m
	-	\sumt{m \in [M] : \sigma(m) = i} \eta_m.
	\)
	Then,
	\[
		\pr{ S = S' }
	\le
		\prodt[L]{\ell=1} \pr{ S_\ell G_{(\ell)} = S'_\ell G_{(\ell)} }
	=
		\prodt[L]{\ell=1} \pr{ \sumt[k]{i=1} v_i Z_{i,\ell} G_{(\ell)} = \id(Q_\ell) }.
	\]
\end{prop}

The randomness in the above set-up comes from the choice $(Z_k)_{i=1}^k$ of \emph{generators}, not from the RW aspect $(\sigma, \eta)$ or $(\sigma', \eta')$:
	it is valid for \emph{any} choices of $(\sigma, \sigma', \eta, \eta')$.
In particular, it applies to both the undirected and directed Cayley graphs, the latter requiring $\eta$ and $\eta'$ to be all-$1$ sequences.

\begin{Proof}
The claimed equality follows immediately from the fact that $Q_\ell$ is Abelian.

We now set up a little notation.
Write
\(
	A_{i,\ell} \cq Z_{i,1} \cdots Z_{i,\ell-1}
\)
and
\(
	B_{i,\ell} \cq Z_{i,\ell+1} \cdots Z_{i,L};
\)
then
\(
	Z_i = A_{i,\ell} Z_{i,\ell} B_{i,\ell}.
\)
(Here, $A_{i,1} \cq \id$ and $B_{i,L} \cq \id$.)
Note that $B_{j,\ell} \in G_{(\ell)}$ for all $j \in [k]$ and $\ell \in [L]$.

Let $\mce_\ell \cq \bra{S' S^{-1} \in G_{(\ell)}}$.
Then,
\[
	\pr{ S = S' }
=
	\prodt[L]{1} \pr{ \mce_\ell \mid \mce_{\ell-1} }.
\]
Now,
\(
	[g, h] \in G_{(\ell)}
\)
and
\(
	hg = gh [h^{-1}, g^{-1}] = gh [g, h]^{-1}
\)
for all  $g \in G$ and $h \in G_{(\ell-1)}$.
So,
%We can hence write $S' S^{-1}$ in the following way:
\[
	S' S^{-1}
=
	M_\ell N_\ell
\cdot
	\rbb{ \prodt[M']{m=1} B_{\sigma'(m), \ell}^{\eta'_m} C'_{\sigma'(m), \ell} }
\cdot
	\rbb{ \prodt[M]{m=1} B_{\sigma(M+1-m), \ell}^{-\eta_{M+1-m}} C'_{\sigma(M+1-m), \ell} }
\]
for some $C_{j,\ell}, C'_{j,\ell} \in G_{(\ell)}$ and $M_\ell$ and $N_\ell$ defined as follows:
\[
	M_\ell
&\cq
	\rbb{ \prodt[M']{m=1} A_{\sigma'(m), \ell}^{\eta'_m} }
\cdot
	\rbb{ \prodt[M]{m=1} A_{\sigma(M+1-m), \ell}^{-\eta_{M+1-m}} };
\\
	N_\ell
&\cq
	\rbb{ \prodt[M']{m=1} Z_{\sigma'(m), \ell}^{\eta'_m} }
\cdot
	\rbb{ \prodt[M]{m=1} Z_{\sigma(M+1-m), \ell}^{-\eta_{M+1-m}} }
\in
	G_{(\ell-1)}.
\]
We thus see that
\(
	\mce_{\ell-1} = \bra{ S' S^{-1} \in G_{(\ell-1)} }
\)
holds if and only if
\(
	\bra{ M_\ell \in G_{(\ell-1)} }
\)
holds.
Crucially, this implies that the indicator $\one{\mce_{\ell-1}}$ of this event's occurrence is independent of $N_\ell$.

We claim that the following is true:
\[
	\text{given that}
\quad
	S' S^{-1} \in G_{(\ell-1)},
\Quad{we have}
	S' S^{-1} \in G_{(\ell)}
\Quad{if and only if}
	M_\ell N_\ell \in G_{(\ell)}.
\]
To prove this, we first make three observations, recalling that $G_{(\ell-1)}/G_{(\ell)}$ is Abelian:
\begin{itemize}[noitemsep, topsep = \smallskipamount, label = \bcdot]
	\item 
	for all $\alpha \in G_{(\ell-1)}$,
	we have
	$\alpha G_{(\ell)} = G_{(\ell)}$ and $(\alpha \beta) G_{(\ell)} = (\alpha G_{(\ell)}) (\beta G_{(\ell)})$ for all $\beta \in G$;
	
	\item 
	$B_{j,\ell}, C_{j,\ell}, C'_{j,\ell} \in G_{(\ell)}$ for all $j \in [k]$ and $N_\ell \in G_{(\ell-1)}$;
	
	\item 
	$S' S^{-1} \in G_{(\ell-1)}$ if and only if $M_\ell \in G_{(\ell-1)}$, and so $M_\ell N_\ell \in G_{(\ell-1)}$.
\end{itemize}
Assume that $S' S^{-1} \in G_{(\ell-1)}$.
Applying these observations in the above formula above gives
\[
	S' S^{-1} G_{(\ell)}
&
=
	\rbr{ M_\ell N_\ell G_{(\ell)} }
\cdot
	\rbb{ \prodt[M']{m=1}
		\rbr{ B_{\sigma'(m), \ell}^{\eta'_m} G_{(\ell)} }
		\rbr{ C'_{\sigma'(m), \ell} G_{(\ell)} }
	}
\\&\qquad
\cdot
	\rbb{ \prodt[M]{m=1}
		\rbr{ B_{\sigma(M+1-m), \ell}^{-\eta_{M+1-m}} G_{(\ell)} }
		\rbr{ C_{\sigma(M+1-m), \ell} G_{(\ell)} }
	}
=
	M_\ell N_\ell G_{(\ell)}.
\]
Thus, $S' S^{-1} \in G_{(\ell-1)}$ if and only if $M_\ell N_\ell \in G_{(\ell-1)}$, as claimed.

Now, $M_\ell$ is independent of $N_\ell$, and so $N_\ell$ is independent also of $\one{\mce_{\ell-1}}$.
Thus,
\[
	\pr{ \mce_\ell \mid \mce_{\ell-1} }
=
	\pr{ M_\ell N_\ell \in G_{(\ell)} \mid \mce_{\ell-1} }
\le
	\maxt{x \in G_{(\ell-1)}} \pr{ x N_\ell \in G_{(\ell)} }.
\]
Now, $Q_\ell = G_{(\ell-1)}/G_{(\ell)}$ is Abelian and $N_\ell$ is a product of generators $Z_{j,\ell}$ and $Z_{j,\ell}^{-1}$ for different $j \in [k]$.
Hence, we are in the set-up of \cref{res-p2:comp:abe-worst-0}.
Applying that lemma,
%we deduce that
\[
	\pr{ \mce_\ell \mid \mce_{\ell-1} }
\le
	\pr{ N_\ell \in G_{(\ell)} }
=
	\pr{ S_\ell G_{(\ell)} = S'_\ell G_{(\ell)} },
\]
using the definition of $N_\ell$.
This proves the desired inequality.
\end{Proof}

\subsection{Evaluation of Abelian-Type Calculations}
\label{sec-p2:comp:eval}

The quotients $Q_\ell$ are Abelian, so the order in which the generators are applied does not matter.
%Given $\sigma : [M] \to [k]$,
Define $W_i \cq \sumt[M]{m=1} \one{\sigma(m) = i}$ for each $i$.
Then, $W = (W_i)_{i=1}^k$ is the RW on $\mbz^k$ run for $M$ steps.

Key in analysing the Abelian-type terms are gcds:
for all $w,w' \in \mbz^k$,
define
\[
	\mfgcd_{(w,w')}
\cq
	\gcd\rbb{ w_1 - w'_1, \: w_2 - w'_2, \: ..., \: w_k - w'_k, \: \abs G }.
%\Quad{abbreviate}
%	\mfgcd \cq \mfgcd_{(W,W')}.
\]
We use these to evaluate the right-hand side of \cref{res-p2:comp:prod-decomp}, culminating in \cref{res-p2:comp:l2-gcd}.
%, using \cite[\cref{res-p0:deferred:unif-gcd}]{HOt:rcg:supp}.

First, we prove an auxiliary lemma akin to \cref{res-p2:cutoff1:unif-gcd}.

\begin{lem}
\label{res-p2:comp:unif-gcd}
	Let $\ell \in [L]$.
%	Write $d_\ell \cq d(Q_\ell)$.
	For all $w,w' \in \mbz^k$,
	we have
	\[
		\sumt[k]{i=1} v_i Z_{i,\ell} G_{(\ell)}
	\sim
		\Unif\rbb{ \mfgcd_{(w,w')} Q_\ell }.
	\]
\end{lem}

\begin{Proof}
\cref{res-p2:comp:sample-Z} says that $Z_{i,\ell} G_{(\ell)} \sim \Unif(Q_\ell)$ independently.
The quotients $Q_\ell$ are Abelian.
%\cite[\cref{res-p0:deferred:unif-gcd}]{HOt:rcg:supp} in the supplementary material
\cref{res-p2:cutoff1:unif-gcd} says that linear combinations of independent random variables in an Abelian group are also uniform, but on the subgroup given by the gcd of the coefficients.
This proves the lemma.
\end{Proof}

This leads us to a bound on $\pr[(w,w')]{ S = S' }$ in terms of a product of $\abs{Q_\ell/\mm Q_\ell} = \abs{Q_\ell} / \abs{\mm Q_\ell}$ over $\ell \in [L]$, for some $\mm$ which is a suitable gcd.
The following lemma controls this product.
%We convert it to a statement about $\widebar G$.

\begin{lem}
\label{res-p2:comp:prod-Ql-sizes}
	For all $\mm \in \mbn$,
	we have
	\(
		\prodt[L]{\ell=1}
		\abs{\mm Q_\ell}
	=
		\abs{\mm \widebar G}.
	\)
\end{lem}

\begin{Proof}
	%
%This is straightforward.
For any Abelian groups $A$ and $B$ and any $\mm \in \mbn$,
we have
\(
	\mm(A \oplus B) = (\mm A) \oplus (\mm B)
\)
and
\(
	\abs{A \oplus B} = \abs A \abs B.
\)
Since $\widebar G$ was defined to be a direct sum of the $Q_\ell$, the claim now follows.
\end{Proof}

%\begin{lem*}
%%\label{res-p2:comp:prod-Ql-sizes}
%	For all nilpotent $G$ and all $\mm \in \mbn$,
%	\blueb{for the old version of $\widebar G$},
%	we have
%	\[
%		\prodt[L]{\ell=1} \abs{\mm Q_\ell}
%	\ge
%		\abs{\mm \widebar G}.
%	\]
%\end{lem*}
%
%\begin{Proof}
%	%
%By induction, it suffices to show that $\abs{\mm \mbz_{pq}} \ge \abs{\mm(\mbz_p \oplus \mbz_q)}$ for all $p,q \in \mbn$.
%This is easy:
%\[
%	\frac{\abs{\mm \mbz_{pq}}}{\abs{\mm (\mbz_p \oplus \mbz_q)}}
%=
%	\frac{ (pq) / \gcd(\mm, pq) }{ (p/\gcd(\mm, p)) (q/\gcd(\mm, q)) }
%=
%	\frac{\gcd(\mm, p) \gcd(\mm, q)}{\gcd(\mm, pq)}
%\ge
%	1.
%\qedhere
%\]
%	%
%\end{Proof}

%Let $S \cq S(t)$ be the RW, and write $W \cq W(t)$ for its auxiliary process, both evaluated at time $t$ (yet to be chosen).
Let $(S',W')$ be an independent copy of $(S,W)$.
Combining \cref{res-p2:comp:prod-decomp,res-p2:comp:unif-gcd,res-p2:comp:prod-Ql-sizes} gives the following corollary.
For $w,w' \in \mbz^k$,
write
%use the shorthand
\[
	\prt[(w,w')]{\cdot}
\cq
	\pr{\, \cdot \mid (W, W') = (w, w')}.
\]

\begin{cor}
\label{res-p2:comp:l2-gcd}
	For all $w, w' \in \mbz^k$,
%	writing $v \cq w - w'$,
	we have
	\[
		n \, \pr[(w,w')]{ S = S' }
	\le
		\prodt[L]{\ell=1} \abs{ Q_\ell } / \abs{ \mfgcd_{(w,w')} Q_\ell }
	=
		\abs{\widebar G} / \abs{ \mfgcd_{(w,w')} \widebar G }
	=
		\absb{ \widebar G / \mfgcd_{(w,w')} \widebar G }.
	\]
\end{cor}

\begin{Proof}
Note that $\abs{Q_\ell}$ divides $\abs G$, and so
\(
	\gcd\rbr{ v_1, ..., v_k, \abs{Q_\ell} }
\le
	\gcd\rbr{ v_1, ..., v_k, \abs G }
\)
for all $v \in \mbz^k$.
Also, for any Abelian subgroup $H$ of $G$, if $\alpha \wr \abs H$ and $\alpha \wr \beta$, then $\beta H \le \alpha H$.
Combined with \cref{res-p2:comp:prod-decomp,res-p2:comp:unif-gcd}, this proves the inequality.
To emphasise, \cref{res-p2:comp:prod-decomp} is valid for any choice of $(\sigma, \sigma', \eta, \eta')$, so, in particular, applying under this conditioning.
The first equality follows immediately from \cref{res-p2:comp:prod-Ql-sizes}.
The second equality follows from Lagrange's theorem.
\end{Proof}

The right-hand side of this corollary depends only on the Abelian group $\widebar G$.
We apply the theory developed in \S\ref{sec-p2:cutoff1}--\S\ref{sec-p2:cutoff3} to bound the mixing time for this Abelian group.
%By applying the results used for Abelian groups, we can prove \cref{res-p2:comp:res}; we explain this now.
%Here, as there, we use a modified $L_2$ calculation; see \cref{res-p2:cutoff1:mod-l2}.

%\begin{lem}[\cref{res-p2:cutoff1:mod-l2}]
%\label{res-p2:comp:mod-l2}
%	For all $t \ge 0$ and all $\mcw \subseteq \mbz^k$,
%	the following inequalities hold:
%	\begin{gather*}
%		d_{G_k}(t)
%	=
%		\tvb{ \pr[G_k]{ S(t) \in \cdot } - \pi_G }
%	\le
%		\tvb{ \pr[G_k]{ S(t) \in \cdot \mid W(t) \in \mcw } - \pi_G }
%	+	\pr{ W(t) \notin \mcw };
%	\\
%		4 \, \ex{ \tvb{ \pr[G_k]{ S(t) \in \cdot \mid W(t) \in \mcw } - \pi_G }^2 }
%	\le
%		n \, \pr{ S(t) = S'(t) \mid W(t), W'(t) \in \mcw } - 1.
%	\end{gather*}
%\end{lem}

\begin{Proof}[Proof of \cref{res-p2:comp:res}]
Let $\mcw \subseteq \mbz^k$ be arbitrary for the moment.
Set
\[
	D
\cq
	n \, \pr{ S = S' \mid \typ } - 1
\Qwhere
	\typ \cq \bra{ W, W' \in \mcw }.
\]
Abbreviate $\mfgcd \cq \mfgcd_{(W,W')}$.
Applying now \cref{res-p2:comp:l2-gcd}, we obtain
\[
	D
\le
	\sumt{\mm \in \mbn}
	\pr{ \mfgcd = \mm \mid \typ } \cdot \absb{\widebar G/\mm \widebar G}
-	1.
\]

This latter expression is purely a statistics of the Abelian group $\widebar G$.
We established the upper bound on mixing by looking at \emph{precisely} this quantity.
Bounding it was one of the main challenges.
There were three different arguments for bounding it, corresponding to different regimes of $k$.
We briefly outline these arguments now.
The choice of $\mcw$ varies from argument to argument.
%We used two distinct methods, as we outline now.
\begin{itemize}[itemsep = 0pt, topsep = \smallskipamount, label = \bcdot]
	\item 
	In \S\ref{sec-p2:cutoff1:upper},
	we upper bounded $\abs{\widebar G/\mm \widebar G} \le \mm^{d(\widebar G)}$;
	we then used unimodality to show that
	\(
		\pr{ \mm \wr W_i \mid W_i \ne 0 } \le 1/\mm,
	\)
	from which we deduced that
	\(
		\pr{\mfgcd = \mm \mid \typ} \le (1/\mm + \pr{W_1 = 0 \mid \typ})^k.
	\)
	
	\item 
	In \S\ref{sec-p2:cutoff2:upper},
	we analysed $(W,W')$ taken modulo $\mm$, for each $\mm$;
	we then used entropic considerations to bound
	\(
		\pr{\mfgcd = \mm \mid \typ} \ll \abs{\widebar G/\mm \widebar G}
	\)
	in a quantitative sense.
	
	\item 
	In \S\ref{sec-p2:cutoff3:upper},
	we combined these two approaches.
\end{itemize}
Instead of reconstructing these arguments, we reference the appropriate places in the previous sections.
For each approach, there are conditions on $(k,\widebar G)$; see \cref{hyp-p2:cutoff1,hyp-p2:cutoff2,hyp-p2:cutoff3}.
At least one of these is satisfied if $1 \ll k \lesssim \log \abs G - \log \abs{\widebar G}$ and $k - d(\widebar G) \gg 1$; see \cref{rmk-p2:cutoff1:hyp,rmk-p2:cutoff2:hyp,rmk-p2:cutoff3:hyp}.

We need to choose the set $\mcw$; see \cref{def-p2:cutoff1:typ,def-p2:cutoff2:typ} for the respective definitions,
replacing $G$ with $\widebar G$ in those definitions.
See \cref{res-p2:cutoff1:l2,res-p2:cutoff2:l2,res-p2:cutoff3:l2} specifically for the results bounding the above sum.
The conclusion of these results is that
\[
	D
\le
	\sumt{\mm \in \mbn}
	\pr{ \mfgcd = \mm \mid \typ } \cdot \absb{\widebar G/\mm \widebar G}
-	1
=
	\oh1.
\]

Combined with the modified $L_2$ calculation of \cref{res-p2:cutoff1:mod-l2}, this completes the proof.
\end{Proof}

\subsection{Cutoff for Nilpotent Groups with Small Commutators}
\label{sec-p2:comp:cor}

We now prove \cref{res-p2:intro:comp:cor:gen,res-p2:intro:comp:cor:special,res-p2:intro:comp:cor:heis},
which relate to nilpotent groups with large Abelianisation.

First, we compare $\tent_*(k, \widebar G = \gab \oplus \widebar{[G, G]})$ with $\tent_*(k, \gab)$.

\begin{prop}
\label{res-p2:comp:cor:t*-comp}
	Let $A$ and $B$ be finite, Abelian groups and $k$ be such that $1 \ll \log k \ll \log \abs A$.
	\begin{itemize}[noitemsep, topsep = \smallskipamount, label = \bcdot]
		\item 
		If $k \eqmathsbox{ent-com}{\lesssim} \log \abs{A \oplus B}$, then
		suppose that $k \gg d(B) \log \abs B$ and $k - d(A) \gg d(B)$.
		
		\item 
		If $k \eqmathsbox{ent-com}{\gg} \log \abs{A \oplus B}$, then
		suppose only that $\log \abs B \ll \log \abs A$.
	\end{itemize}
	(In either case, $\log \abs B \ll \log \abs A$.)
	Then,
	\(
		\tent_*(k, A \oplus B)
	\eqsim
		\tent_*(k, A).
	\)
\end{prop}

These conditions imply that $A \oplus B$ should be viewed as a `small perturbation' of $A$.
The proof is technical, relying on auxiliary results on entropy; it is deferred to \cite[\cref{res-p0:re:abe-com:eqsim}]{HOt:rcg:supp}.

\begin{Proof}[Proof of \cref{res-p2:intro:comp:cor:gen}]
The lower bound of $\tent_*(k, \gab)$ follows by projecting onto the Abelianisation, which is an Abelian group.
The argument is analogous to that of Approach \#2 in \S\ref{sec-p2:cutoff2:lower}.

The conditions of \cref{res-p2:intro:comp:cor:gen} match those of \cref{res-p2:comp:cor:t*-comp} when $A = \gab$ and $B = \widebar{[G, G]}$.
Hence, $\tent_*(k, \gab) \eqsim \tent_*(k, \widebar G = \gab \oplus \widebar{[G, G]})$.
The upper bound of $\tent_*(k, \widebar G)$ is immediate from \cref{res-p2:intro:comp:nil-abe} if $k \lesssim \log \abs G$.
It was already known for $k \gg \log \abs G$; recall \S\ref{sec-p2:intro:previous-work:ad-conj} or see \cref{res-p2:conc-rmks:roichman}.
\end{Proof}

%Indeed, $\widebar G = \gab \oplus \widebar{[G,G]}$, so certainly $\tent_*(k, \widebar G) \ge \tent_*(k, \gab)$.

To prove \cref{res-p2:intro:comp:cor:special,res-p2:intro:comp:cor:heis}, we need an asymptotic evaluation of the entropic times.
Recall that $\tent_\mm(k, \mm^r)$ is the time at which the entropy of RW on $\mbz_\mm^k$ reaches $\log(\mm^r) = \log \abs{\mbz_\mm^r}$.

\begin{prop}[{\cite[\cref{res-p0:re:p:t0}]{HOt:rcg:supp}}]
\label{res-p2:comp:cor:t*-eval}
	Let
	\(
		\zeta_\mm
	\cq
		\tfrac1k (k - r) \log \mm.
	\)
	Suppose that $1 \ll k \lesssim r \log \mm$.
%	The following hold:
	\begin{alignat*}{2}
		\text{If}
	\quad
		\zeta_\mm &\eqmathsbox{re:m:t0}{\ll} 1,
	&\Quad{then}
		\tent_\mm(k, \mm^r)/k &\eqsim \tfrac12 \log(1/\zeta_\mm) / \rbb{ 1 - \cos(2\pi/\mm) }.
	\\
		\text{If}
	\quad
		\zeta_\mm &\eqmathsbox{re:m:t0}{\gtrsim} 1,
	&\Quad{then}
		\tent_\mm(k, \mm^r)/k &\asymp \mm^2 e^{-2\zeta_\mm} = (\mm^r)^{2/k}.
	\intertext{Suppose further that $1 \ll k \ll r \log m$.}
		\text{If}
	\quad
		\zeta_\mm &\eqmathsbox{re:m:t0}{\gg} 1,
	&\Quad{then}
		\tent_\mm(k, \mm^r)/k &\eqsim \mm^2 e^{-2\zeta_\mm} / (2 \pi e) = (\mm^r)^{2/k} / (2 \pi e).
	\end{alignat*}
%	Note that
%	\(
%		1 - \cos(2\pi/\mm)
%	\eqsim_{\toinf \mm}
%		2 \pi^2 / \mm^2
%	=
%		2 \pi^2 \mm^{-2r/k} e^{2\zeta_\mm}.
%	\)
\end{prop}

\begin{Proof}[Proof of \cref{res-p2:intro:comp:cor:special}]
The lower bound argument is exactly the same as for \cref{res-p2:intro:comp:cor:gen}.

For the upper bound, we slightly refine the argument used to prove \cref{res-p2:intro:comp:nil-abe}.
First, we claim that $\gab \cong \mbz_p^r$.
Indeed, the Frattini subgroup $\Phi(G)$ satisfies $\Phi(G) = [G,G] G^p$ when $G$ is a $p$-group, where $G^p \cq \langle g^p \mid g \in G \rangle$.
By definition of being special, $\Phi(G) = [G,G]$.
Thus, $G^p \le [G,G]$.
In particular, the Abelianisation is of exponent $p$, as required.
Thus, $\widebar G \cong \mbz_p^\ell$ as $\gcom \cong \mbz_p^s$.

The general arguments for Abelian groups used to prove \cref{res-p2:intro:tv} can be specialised to the group $\mbz_p^\ell$.
We do not include the details here, but rather defer them to \cite[\cref{res-p5:intro:p}]{HOt:rcg:abe:extra}.
The conditions for this approach are only $k \ge \ell$, rather than $k - \ell \gg 1$ as previously.

We turn to the entropic time.
We have $\mm \mbz_p^r = \mbz_p^r$ unless $p \wr \mm$, as $p$ is prime.
Thus, the worst-case $\mm$ in
\(
	\tent_*(k, \mbz_p^r)
=
	\maxt{\mm \wr p^r}
	\tent_\mm(k, \mbz_p^r)
\)
is $\mm = p$,
and so
\(
	\tent_*(k, \mbz_p^r)
=
	\tent_p(k, \mbz_p^r).
\)

For a nilpotent group $G$, for $Z$ to generate $G$ it suffices that $Z$ generates $\gab$; here, $\gab = \mbz_p^r$.
We analyse Abelian group generation in \cite[\cref{res-p5:gen:dichotomy}]{HOt:rcg:abe:extra}.
That lemma shows that $k - r \gg 1$ is always sufficient to generate the group \whp, but if $p \gg 1$ then merely $k - r > 0$ is sufficient.

Finally, if $k - r \asymp k$ and $p \gg 1$,
%then \cite[\cref{res-p0:re:p:t0}]{HOt:rcg:supp}
\cref{res-p2:cutoff1:ent-times,res-p2:comp:cor:t*-eval}
gives $\tent_p(k, \mbz_p^r) \eqsim \tent_\infty(k, \mbz_p^r) = \tpro_0(\abs{\mbz_p^r})$.
Indeed, in the notation there, $\zeta_p = \tfrac1k (k - r) \log p$ is the relative entropy; so, $\zeta_p \gg 1$.
\end{Proof}

For the Heisenberg group $\HH_{m,d}$, we have the explicit expression $\HH_{m,d}^\ab \cong \mbz_m^{2d-4}$, even for $m$ not prime.
This allows us to evaluate $\tent_*(k, \HH_{m,d}^\ab)$ even when $m$ is not prime, provided $m \gg 1$.
%When $m = p$ is prime, we can relax the condition $k - 2d \gg 1$ to $k \ge 2d-3$.

\begin{Proof}[Proof of \cref{res-p2:intro:comp:cor:heis}]
Let $r \cq 2d-4$.
We have $\hab \cong \mbz_m^r$ and $\hcom \cong \mbz_m$.
The conditions of \cref{res-p2:intro:comp:cor:heis} are precisely those required to apply \cref{res-p2:intro:comp:cor:gen}.
Hence, there is cutoff at $\tent_*(k, \hab)$ \whp.
It remains to evaluate $\tent_*(k, \hab) = \tent_*(k, \mbz_m^r)$ when $k - r \asymp k$ and $m \gg 1$.

First, observe that $\mm \mbz_m^r = \mbz_{m/\gcd(\mm, m)}^r$ for all $\mm$.
Hence, by replacing $\mm$ with $\gcd(\mm, m) \le \mm$, we need only consider $\mm$ with $\mm \wr m$.
Next, $\abs{\mbz_m^r / \mbz_{m/\mm}^r} = \mm^r$.
Hence,
\[
	\tent_*(k, \mbz_m^r)
=
	\maxt{\mm \wr m}
	\tent_\mm(k, \mm^r).
\]

In the notation of \cref{res-p2:comp:cor:t*-eval},
$\zeta_\mm = \tfrac1k (k - \mm) \log \mm \asymp \log \mm$.
Hence, $\zeta_\mm \gg 1$ if (and only if) $\mm \gg 1$.
If the maximising $\mm$, call it $\mm_*$, satisfies $\mm_* \gg 1$, then
\cref{res-p2:cutoff1:ent-times,res-p2:comp:cor:t*-eval}
give
\[
	\maxt{\mm \wr m}
	\tent_\mm(k, \mm^r)
=
	\tent_{\mm_*}(k, \mm_*^r)
\eqsim
	k \mm_*^{2r/k} / (2 \pi e)
\eqsim
	\tent_\infty(k, \mm_*^r).
%	\tpro(\mm_*^r).
\]

%We claim that $\mm_* = m$.
%We need to show that $\mm_* = m$.
It remains to show that $\mm_* = m$.
The map $\mm \mapsto k \mm^{2r/k} / (2 \pi e)$ is increasing.
So, if $\mm_* \gg 1$, then in fact $\mm_* = m$.
Similarly, \cref{res-p2:comp:cor:t*-eval} show that $\tent_\mm(k, \mm)/k \lesssim \mm^{2r/k}$
as $\zeta_\mm \asymp \log \mm \gtrsim 1$.
%whenever $\zeta_\mm \gtrsim 1$.
Hence,
\[
	\tent_\mm(k, \mm^r)/k
\lesssim
	\mm^{2r/k}
\ll
	m^{2r/k}
\asymp
	\tent_m(k, \mm^r)/k;
\Quad{thus,}
	\mm_* = m.
\qedhere
\]
\end{Proof}

Finally, after \cref{res-p2:intro:comp:cor:special} we mentioned that special groups are ubiquitous amongst $p$-groups of a given size.
We elaborate on this claim in the following remark.
%Finally, we expand on \cref{rmk-p2:intro:comp:p-groups}, which gave bounds on the number of certain $p$-groups.

\begin{rmkt}
\label{rmk-p2:comp:cor:p-groups:hist}
In their classical work \cite{H:enumerating-p-groups}, \citeauthor{H:enumerating-p-groups} gave upper and lower bounds on the number groups of size $p^\ell$ for a prime $p$.
The upper bound was later refined by \textcite{S:enumerating-p-groups}.
Together, they show that this number is $p^{(2/27) \ell^3 \pm \Oh{\ell^{8/3}}}$.
The lower bound $p^{(2/27) \ell^3 - \Oh{\ell^2}}$ is obtained from \textcite[Theorem~2.1]{H:enumerating-p-groups} by counting step-2 groups whose Frattini group is equal to the centre and is elementary Abelian of size $p^s$ and of index $p^r$, where $r = \ell - s$.
It is classical that such a group is special if and only if it has exponent $p$, ie every element other than the identity has order $p$.

\textcite{H:enumerating-p-groups} showed that the number of such groups of size $p^\ell$ is between $p^{(1/2)sr(r+1) - r^2 - s^2}$ and $p^{(1/2)sr(r+1) - s(s-1)}$ if $s \le \tfrac12 r(r+1)$ and 0 otherwise.
A small variant of his argument shows that the number of special groups of size $p^{s+r}$ whose commutator subgroup is of size $p^s$ is between $p^{(1/2) sr(r-1) - r^2 - s^2}$ and $p^{(1/2) sr(r-1) - s(s-1)}$ for $s \le \tfrac12 r(r-1)$.
(In \cite{H:enumerating-p-groups}, one includes in $F$ all elements of order $p$.
See also the short argument in \textcite[Page~152]{S:enumerating-p-groups}; there, the change is considering the case that $b(i,j) = 0$ for all $1 \le i \le r$ and $1 \le j \le s$.)
Taking $r$ and $s$ such that $\abs{ r - 2s } < 3$, we see that the logarithm of the number of groups of size $p^\ell$ is dominated by special groups.
\end{rmkt}

%\begin{rmkt}
%\label{rmk-p2:comp:cor:p-groups:ext}
%	%
%For special $p$-groups, the commutator subgroup is elementary Abelian (by definition).
%Thus $\widebar{[G,G]} = [G,G]$ and $d(\widebar{[G,G]}) = \log_p \abs \gcom$.
%(\textit{Extra special} groups have $\gcom \cong \mbz_p$.)
%It is noted in \cref{rmk-p2:intro:comp:p-groups} that, for special groups, if $\abs G = p^\ell$ then $\widebar G = \mbz_p^\ell$.
%Indeed, the Frattini subgroup $\Phi(G)$ satisfies $\Phi(G) = [G,G] G^p$ when $G$ is a $p$-group where $G^p \cq \langle g^p \mid g \in G \rangle$. By definition of being special, $\Phi(G) = [G,G]$, thus $G^p \le [G,G]$.
%In particular, the Abelianisation is of exponent $p$.
%Thus, as for Heisenberg groups, instead of applying the general \cref{res-p2:intro:tv} in \cref{res-p2:intro:comp:nil-abe}, we can apply \cite[\cref{res-p5:intro:p}]{HOt:rcg:abe:extra} which is specialised to $\mbz_p^\ell$.
%	%
%\end{rmkt}

\subsection{Expander Graphs of Nilpotent Groups with $k \gtrsim \log \abs G$}
\label{sec-p2:gap}

The isoperimetric constant $\Phi_*$ can be defined for general Markov chains; see \cite[\S 7.2]{LPW:markov-mixing}. The isoperimetric constant of a graph is that of the nearest-neighbour RW on the graph.
It is easy to see that the isoperimetric constant of a Markov chain, of its time reversal and of its additive symmetrisation are all equal.
But, for any generators $Z$, the additive symmetrisation of the RW on $G^+(Z)$ is the RW on $G^-(Z)$;
so, $\Phi_*(G^+(Z)) = \Phi_*(G^-(Z))$.
It thus suffices to work with undirected~Cayley~graphs.

We analyse the spectral gap via considering the $1/n^c$-mixing time for some $c > 0$.

\begin{prop}
\label{res-p2:gap:mix}
	Let $G$ be a nilpotent group.
	Suppose that $k - d(\widebar G) \asymp k \asymp \log \abs G$.
	Let $\eps > 0$ and set $t \cq (1 + \eps) \tent_*^-(k,G)$.
	Then, there exists a constant $c > 0$ so that $d_{G_k}^-(t) \le \abs G^{-c}$ \whp.
\end{prop}

\begin{Proof}
Consider first Abelian $G$; here, $G = \widebar G$.
Since \cref{hyp-p2:sep} is satisfied, $d_{G_k}(t) \le e^{-c'(k-d(G))}$ \whp for some constant $c' > 0$, by \cref{res-p2:sep:tv-quant}.
The claim now follows as $k - d(G) \asymp \log \abs G$ here.

Consider now nilpotent $G$; here, $G \ne \widebar G$.
We apply our nilpotent-to-Abelian method.
There, we upper bounded the modified $L_2$ distance for the RW on $G$ (at time $t$) by the modified $L_2$ distance for the RW on $\widebar G$ (at time $t$); see specifically \cref{res-p2:comp:prod-decomp,res-p2:comp:unif-gcd,res-p2:comp:l2-gcd}.
For Abelian groups, we used the modified $L_2$ calculation (in \S\ref{sec-p2:cutoff1}--\S\ref{sec-p2:cutoff3}).
Thus, the nilpotent case is an immediate application of the nilpotent-to-Abelian method and Abelian case.
\end{Proof}

We apply \cref{res-p2:gap:mix} along with standard mixing-type results.

\begin{Proof}[Proof of \cref{res-p2:gap:res}]
	%
%Let $n \cq \abs G$.
As noted in \cref{rmk-p2:intro:gap}, it suffices to consider only $k \asymp \log n$.

First, we use the discrete Cheeger inequality, for reversible Markov chains:
	\(
		\gamma \lesssim \Phi_* \lesssim \sqrt \gamma ,
%	\Quad{where}
%		\text{$\gamma$ is the \textit{spectral gap}};
	\)
	where
		$\gamma$ is the spectral gap;
	and
%		$\Phi_*$ the isoperimetric constant;
	see, eg, \cite[Theorem~13.10]{LPW:markov-mixing}.
Thus, it suffices to show that $\gamma \asymp 1$.

Next, we use a standard relation between the mixing time and spectral gap:
	for a reversible Markov chain
	with
		invariant distribution uniform on $n$ states,
%	and
		mixing time $\tmix$
	and
		spectral gap~$\gamma$,
%	we~have
	\[
		\tmix(1/n^c) \asymp \gamma^{-1} \log n
	\Quad{for any constant}
		c > 0;
	\]
	see, eg, \cite[Theorem~20.6 and Lemma~20.11]{LPW:markov-mixing}.
	Thus, $\gamma \asymp 1$ if $\tmix(1/n^c) \lesssim \log n$ for such $c$.

%Combining these, to show that the graph is an expander, ie $\Phi_* \asymp 1$, it suffices to show that $\tmix(1/n^c) \lesssim \log n$ for some constant~$c > 0$.
%(For any lazy or continuous-time chain, we always have $\tmix(1/n^c) \gtrsim \log n$.)
%
The mixing claim follows immediately from \cref{res-p2:gap:mix} and the fact that $k \asymp \log \abs G$.
%The proof is complete.
	%
\end{Proof}

\section{Concluding Remarks and Open Questions}
\label{sec-p2:conc-rmks}

%In this section, we give some concluding remarks.
\begin{itemize}[itemsep = 0pt, topsep = \smallskipamount, label = $\bcdot$]
	\item [\S\ref{sec-p2:conc-rmks:const-k}]
	We discuss lack of cutoff in the regime where $k$ is a fixed constant.
	
	\item [\S\ref{sec-p2:conc-rmks:roichman}]
	We give a very short proof, which is a small variant on Roichman's argument \cite[Theorem~2]{R:random-random-walks}, establishing an upper bound on mixing, for arbitrary groups and any $k \gg \log \abs G$.
	
	\item [\S\ref{sec-p2:conc-rmks:l2-mix}]
	We briefly discuss cutoff in other metrics, namely $L_2$ and relative entropy.
		
	\item [\S\ref{sec-p2:open-conj}]
	To conclude, we discuss some questions which remain open and gives some conjectures.
\end{itemize}
Throughout this section, we only sketch details.

\subsection{Lack of Cutoff When $k$ Is Constant}
\label{sec-p2:conc-rmks:const-k}

%\input{anc/short_sections/CRC_Constk.tex}
%
%\addtocounter{thm}{-1}
%\refstepcounter{thm}
%\label{res-p2:conc-rmks:const-k}

Throughout the paper, we have always been assuming that \toinf k \asinf{\abs G}.
It is natural to ask what happens when $k$ does not diverge.
This case has actually already been covered by \textcite{DSc:growth-rw}, using their concept of \textit{moderate growth}.
In short, there is no cutoff.

\citeauthor{DSc:growth-rw} establish this not only for Abelian groups, but for nilpotent groups.
Recall that a group $G$ is called \textit{nilpotent of step at most $L+1$} if its lower central series terminates in the trivial group after at most $L$ steps:
	$G_{(0)} \cq G$ and $G_{(\ell)} \cq [G_{(\ell-1)}, G]$ for $\ell \in \mbn$ with $G_L = \bra{\id}$.

\smallskip

For a Cayley graph $G(Z)$, use the following notation.
	Write $\Delta \cq \diam G(Z)$ for its diameter.
	For the lazy simple random walk on $G(Z)$,
	write
		$\trel \cq \trel\rbr{ G(Z) }$ for the relaxation time (ie, the inverse of the spectral gap)
	and
		$\tmix \cq \tmix\rbr{ \tfrac14; G(Z) }$ for the (TV) mixing time.
%		, for $\eps \in (0,1)$.
When considering sequences $(G_N(Z_{(N)}))_\Ninn$, add an $N$-sub/superscript, analogously to before.

We phrase the result of \textcite{DSc:growth-rw} in our language.

\begin{thm}[{cf \cite[Corollary~5.3]{DSc:growth-rw}}]
\label{res-p2:conc-rmks:const-k}
	Let $(G_N)_\Ninn$ be a sequence of finite, nilpotent groups.
	For each $\Ninn$, let $Z_{(N)}$ be a symmetric generating set for $G_N$ and write $L_N$ for the step of $G_N$.
	Suppose that $\sup_N \abs{Z_{(N)}} < \infty$ and $\sup_N L_N < \infty$.
	Then,
	\[
		\tmix^N / |Z_{(N)}|
	\lesssim
		\Delta_N^2
	\lesssim
		\trel^N
	\lesssim
		\tmix^N
	\quad
		\asinf N.
	\]
	In particular, $(\tmix^N)_\Ninn$ does not exhibit the cutoff phenomenon.
\end{thm}

We give a very brief exposition of the results of \textcite{DSc:growth-rw}, including the definition of moderate growth, leading to this conclusion in \cite[\S\ref{sec-p5:const-k}]{HOt:rcg:abe:extra}.

\subsection{A Variant on Roichman's Argument}
\label{sec-p2:conc-rmks:roichman}

In this subsection we give a very short argument upper bounding the mixing time for arbitrary groups and $k \gg \log \abs G$; it is a small modification of Roichman's argument \cite[Theorem~2]{R:random-random-walks}, but it applies in both the undirected and directed cases.
(Roichman \cite[Theorem~1]{R:random-random-walks} deals with the directed case, but requires additional matrix algebra machinery.)

\begin{thm}%[Cutoff Whenever $k \gg \log \abs G$]
\label{res-p2:conc-rmks:roichman}
	Let $\eps > 0$.
	Let $G$ be a finite group and $k$ and integer with $k \gg \log \abs G$ and $\log k \ll \log \abs G$.
	Then, the RW on $G_k^\pm$ is mixed \whp at time $(1 + \eps) T(k, \abs G)$,
%	Then, the mixing time of the RW on $G_k^\pm$ is at most $(1 + \eps) T(k, \abs G)$ \whp,
	where
	\[
		T(k, n)
	\cq
		\log n / \log(k / n)
	\Qfor
		n, k \in \mbn.
	\]
	In particular, this upper bound $T(k, \abs G)$ does not depend on the algebraic structure of the group.
\end{thm}

\begin{Proof}
	%
%The proof proceeds as follows.
%Assume that $k \gg \log \abs G$ and $\log k \ll \log \abs G$.
Let $\eps > 0$ and set $t \cq (1 + \eps) \log \abs G / \log(k/\log \abs G)$.
Note that $1 \ll t \ll k$.
Choose some $\omega \gg 1$, diverging arbitrarily slowly;
set $t_\pm \cq \floor{t \rbr{ 1 \pm \omega/\sqrt t}}$ and $L \cq \omega \floor{t^2/k}$.
The number of generators picked at most once is at least $k - L$ \whp; of these, the number picked exactly once lies in $[t_-, t_+]$ \whp.
Take $\typ$ to be the event that these two conditions hold for two independent copies, $W$ and $W'$.
We use a modified $L_2$ calculation (see, eg, \cref{res-p2:cutoff1:mod-l2}) meaning that we need~to~control
\[
	D
\cq
	\abs G \, \pr{ S = S' \mid W = W', \, \typ } - 1.
\]

Let $\mce$ be the event that some generator is used precisely once in $W$ and never in $W'$, or~vice~versa:
\[ \textstyle
	\mce
\cq
	\bigcup_{i \in [k]} \rbb{
		\bra{ \abs{W_i} = 1, \: \abs{W'_i} = 0 }
	\cup
		\bra{ \abs{W'_i} = 1, \: \abs{W_i} = 0 }
	}.
\]
Then, $S' \cdot S^{-1} \sim \Unif(G)$ on $\mce$.
Indeed, if $Z \sim \Unif(G)$ and $X,Y \in G$ are independent of $Z$, then $XZY \sim \Unif(G)$;
here, $Z$ corresponds to one of the generators used precisely once in $W$ and not in $W'$ or vice versa, with the obvious choice of $X$ and $Y$ so that $XZY = S' S^{-1}$.
%Hence
%\[
%	\pr{ S = S' \mid W = W', \, \mce, \, \typ }
%=
%	1/\abs G.
%\]
Off $\mce$, every generator picked once in $W$ must be picked at least once in $W'$, and vice versa.
There are at most $L$ generators which are picked more than once in $W'$.
Thus,
\[
	\pr{\mce \mid \typ}
\le
	\maxt{a \in [t_-, t_+], b \le L}
	1/\binomt{k-b}{a-b}
=
	1/\binomt{k-L}{t_- - L}.
\]
An application of Stirling's approximation shows that this probability is $\oh{1/\abs G}$ when $\omega$ diverges sufficiently slowly.
Combined with the modified $L_2$ calculation, this proves the upper bound.
\end{Proof}

\medskip

Finally, consider the case $k = \abs G^\alpha$ for some fixed $\alpha \in (0,1)$.
The discrete-time chain cannot be mixed at time $\ceil{1/\alpha} - 1$ by considering the size of its support,
but noting that $\binom kt \gg \abs G$ for $t \cq \floor{1/\alpha} + 1$, by the above argument we see that the walk is mixed \whp after $t$ steps.

\citeauthor{D:phd} proves a more general statement than this which allows the generators to be picked from a distribution other than the uniform distribution; see \cite[Theorems~3.3.1 and~3.4.7]{D:phd}.
%(In \cite[Theorem~3.4.7]{D:phd}, take $\beta \cq 2$ to recover the case of uniform generators.)

\subsection{Mixing in Different Metrics}
\label{sec-p2:conc-rmks:l2-mix}

One can also consider cutoff in the $L_2$ distance.
Recall that, for a time $t \ge 0$,
\[
	d_{G_k}^{(2)}(t)
=
	\normb{ \pr[G_k]{S(t) \in \cdot} - \pi_G }_{2, \pi_G}
=
	\rbb{ \abs G^{-1} \sumt{g \in G} \rbb{ \abs G \, \prt[G_k]{ S(t) = g } - 1 }^2 }^{1/2}.
\]
%One can then define mixing and cutoff with respect to $L_2$ analogously to TV ($L_1$) distance.

%\begin{defn*}
%	Let $\mm \in \mbz \cup \bra{\infty}$.
%	Let $\tilde \tent_\mm(k,G)$ be the time $t$ at which the return probability for SRW on $\mbz_\mm^k$ at time $2t$ is $\abs{G/\mm G}^{-1}$.
%%	By independence of coordinates of the SRW on $\mbz_\mm^k$,
%	Equivalently,
%	\(
%		\tilde \tent_\mm(k,G) \cq k s
%	\)
%	where $s$ is the unique solution to
%	\(
%		\pr{ X_{2s} = 0 } = \abs{G/\mm G}^{-1/k}
%	\)
%	where $(X_s)_{s\ge0}$ is a rate-$1$ SRW on $\mbz_\mm$.
%	Set $\tilde \tent_*(k,G) \cq \max_{\mm \in \mbn} \tilde \tent_\mm(k,G)$.
%\end{defn*}

For reasons explained below, we \emph{strongly believe} the following is true---and can be proved in the framework which we have developed in this article.
It was stated as \cref{conj-p2:intro:l2} in \S\ref{sec-p2:intro:res:tv}.

\begin{customconj}{\ref{conj-p2:intro:l2}}
	For $\mm \in \mbz \cup \bra{\infty}$,
	let $\tilde \tent^\pm_\mm \cq \tilde \tent^\pm_\mm(k,G)$ be the time $t$ at which the return probability for RW on $\mbz_\mm^k$ at time $2t$ is $\abs G^{-1}$.
	Set $\tilde \tent^\pm_*(k,G) \cq \max_{\mm \in \mbn} \tilde \tent^\pm_\mm(k,G)$.
	Then, under similar conditions to those of \cref{res-p2:intro:tv}, \whp, the RW on $G_k$ exhibits cutoff in the $L_2$ metric at time $\tilde \tent^\pm_*(k,G)$.
\end{customconj}

%\begin{customconj}{\ref{conj-p2:intro:l2}}
%	Let $G$ be an Abelian group and suppose that $1 \ll k \lesssim \log \abs G$.
%	Suppose that $k - d(G) \gg 1$.
%	Then, \whp, the RW on $G^\pm_k$ exhibits cutoff in $L_2$ at time $\tilde \tent_*(k,G)$.
%\end{customconj}

The main change is that we now no longer perform a modified $L_2$ calculation.
Replacing $r_* \cq \tfrac12 \abs G^{1/k} (\log k)^2$ with $r_* \cq \tfrac12 \abs G^{1/k} \log \abs G$, local typicality then holds with probability $1 - \oh{1/\abs G}$; cf \cref{res-p2:gap:mix}.
Thus, we may condition on local typicality as this can only change the $L_2$ distance by at most an $\oh1$ additive term.
On the other hand, we no longer condition on global typicality.
Instead, we must handle directly terms like
\(
	\pr{ W = W' }
\text{ or }
	\pr{ W_\mm = W'_\mm }.
\)

For Approach \#1, we must handle a gcd. Under the assumption that $1 \ll k \lesssim \sqrt{\log \abs G}$, increasing $r_*$ as we have has little effect on the proof, in essence because $(\log n)^d = n^{\oh1}$.
In Approach \#2, we replace $\abs \mch$ by $\abs \mch \log \abs G$, but still $k \gg \sqrt{\log \abs G}$
implies that $k \gg \log(\abs \mch \log \abs G)$.
Lastly, the combination of the two approaches works when 
\(
	\sqrt{\log \abs G / \log \log \abs G}
\ll
	k
\lesssim
	\sqrt{\log \abs G}.
\)

\medskip

Using somewhat similar adaptations, we believe that cutoff in the relative entropy (abbreviated \textit{RE}) distance can be established. In this case, we quantify the probability with which global typicality holds:
	the maximal relative entropy of a measure on $G$ with respect to $\pi_G$ is $\log \abs G$;
	thus, naively at least, to condition on global typicality we desire it to hold with probability $1 - \oh{1/\log \abs G}$---for $L_2$ we had $1 - \oh{1/\abs G}$.
Also, one should modify local typicality as previously.
This gives conditions on $k$ and $d(G)$.
Under such conditions, the RE and TV cutoff times should then be the same.

We believe that with more effort these conditions can be improved via obtaining some estimates on the relative entropy given that global typicality fails.

\subsection{Open Questions and Conjectures}
\label{sec-p2:open-conj}

We close the paper with some questions which are left open.

\newcounter{oq}

\refstepcounter{oq}
\label{oq-p2:product-condition}
\subsubsection*{\theoq: Does the Product Condition Imply Cutoff?}

The problem of singling out abstract conditions under which the cutoff phenomenon occurs
%without necessarily pinpointing its precise location,
has drawn considerable attention.
For a reversible Markov chain $X$, write $\tmix(X)$ for its mixing time and $\gamma_\textrm{gap}(X)$ for its spectral gap.
In \citeyear{P:conj-product-condition}, \textcite{P:conj-product-condition} proposed a simple spectral criterion for a sequence $(X^N)_\Ninn$ of reversible Markov chains, known as the \textit{product condition}:
\begin{center}
\(
	\textrm{cutoff}
\Quad{is equivalent to}
%	\textrm{the \textit{product condition}},
%\Quad{ie}
	\tmix(X^N) \gamma_\textrm{gap}(X^N) \to \infty \ \asinf N.
\)
\end{center}

%The Peres conjecture asks
%%for a sequence $(\mcg_N)_\Ninn$ of graphs,
%if
%\[
%	\textrm{cutoff}
%\Quad{is equivalent to}
%	\textrm{the product condition,} \
%	\textrm{ie $\tmix(\mcg) \gamma_\textrm{gap}(\mcg) \to \infty$ as $\abs \mcg \to \infty$}.
%%	\textrm{ie $\tmix(\mcg_N) \gamma_\textrm{gap}(\mcg_N) \to \infty$ as $N \to \infty$}.
%\]
%%(Implicitly this is for a sequence of graphs.)

It is well-known that the product condition is a necessary condition for cutoff; see, eg, \cite[Proposition~18.4]{LPW:markov-mixing}.
It is relatively easy to artificially create counter-examples, but these are not `natural'; see, eg,
\cite[\S 18]{LPW:markov-mixing}
%\cite[\S 6]{CSc:cutoff-ergodic}
where constructions due to Aldous and due to Pak are described.
The product condition is widely believed to be sufficient for most `natural' chains.

We conjecture that the product condition implies cutoff for random Cayley graph of Abelian groups.
In fact, we conjecture this whenever $G$ is \textit{nilpotent} of bounded \textit{step} (denoted $\step G$), ie has lower central series terminating at the trivial group and this sequence is of bounded length.
%The \textit{step} (or \textit{class}) of a nilpotent group is the length of its lower central series.

\begin{conj-ind}\theoq
	Let $(G_N)_\Ninn$ be a sequence of finite, nilpotent group and $(Z_{(N)})_\Ninn$ a sequence of subsets with $Z_{(N)} \subseteq G_N$ for all $\Ninn$.
	For each $\Ninn$,
	write $\tmix^N$, respectively $\gamma_\textrm{gap}^N$, for the mixing time, respectively spectral gap, of the SRW on $G^-_N(Z_{(N)})$.
	
	Suppose that $\limsup_\Ninf \step G_N < \infty$ and that the product condition holds,
	ie
	\(
		\tmix^N \gamma_\textrm{gap}^N \to \infty
	\)
	\asinf N.
	Then, the sequence of SRWs exhibits cutoff.
\end{conj-ind}

An equivalence between the product condition and cutoff has been established for birth-and-death chains by \textcite{DLP:cutoff-birth-death} and, more generally, for RWs on trees by \textcite{BHP:cutoff-rev}.
It is believed to imply cutoff for the SRW on transitive expanders of bounded degree, but this is known only in the case of Ramanujan graphs, due to \textcite{LP:ramanujan}.

\refstepcounter{oq}
\label{oq-p2:choice-of-gen}
\subsubsection*{\theoq: An Explicit Choice of Generators}

We have shown that choosing the generators $Z$ uniformly gives cutoff \whp at a time which does not depend on $Z$, in many regimes.
In particular, this means that there is cutoff for almost all choices of generators at a time independent of the choice of generators.
This `almost universal' mixing time is given by $\tent_*(k,G)$ from \cref{def-p2:cutoff2:ent-orig}.
A question raised to us by \textcite{D:cayley:private} is to find \emph{explicit} sets of generators for which cutoff occurs; see also \cite[Chapter~4G, Question~2]{D:group-rep}.

%\vspace{-\baselineskip}
\begin{openproblem-ind}\theoq
	Let $G$ be an Abelian group and $1 \ll k \lesssim \log \abs G$.
	Find an explicit choice of generators $Z$ (implicitly a sequence) so that the RW on $G(Z)$ exhibits cutoff.
	Further, find generators so that the cutoff time is $\tent_*(k,G)$ asymptotically.
\end{openproblem-ind}

\textcite[Theorem~1.11]{H:cutoff-cayley-<} shows for the cyclic group $\mbz_p$ with $p$ prime that the choice $Z \cq \sbr{ 0, \pm1, \pm2, ..., \pm2^{\ceil{\log_2 p}-1} }$,
which he describes as ``an approximate embedding of the classical hypercube walk into the cycle'',
gives rise to a random walk on $\mbz_p$ which has cutoff.
The cutoff time is not the entropic time, however.
Although the entropic time is the mixing time for `most' choice of generators,
finding an explicit choice of generators which gives rise to cutoff at the entropic time is still open---even for the cyclic group of prime order.

%\vspace*{3\bigskipamount}
%{\color{red} \huge check references, particularly ``customeprint'' parts}

\renewcommand{\bibfont}{\sffamily}
%\renewcommand{\bibfont}{\sffamily\small}
%\renewcommand{\bibfont}{\small}
%\appto\bibfont{\sffamily\small}
\printbibliography[heading=bibintoc]

\end{document}